\documentclass[11pt,a4paper]{article}
\usepackage{latexsym,amssymb,amsmath,mathrsfs}
\usepackage{sectsty}
\usepackage{color}
\usepackage{bm}
\usepackage[T1]{fontenc}
\usepackage[anchorcolor=blue,%
 bookmarks=true,%
 bookmarksnumbered=true,%
]{hyperref}
\usepackage{cite}
\usepackage{pst-all}
\usepackage{graphicx}

\def\Ch{\mathrm{Ch}}
\allsectionsfont{\normalsize\sc\center}

\newtheorem{thm}{Theorem}[section]
\newtheorem{pr}[thm]{Proposition}
\newtheorem{cor}[thm]{Corollary}
\newtheorem{lem}[thm]{Lemma}
\newtheorem{df}[thm]{Definition}
\newtheorem{rem}[thm]{Remark}

\newenvironment{pf}{\par\begin{trivlist}%
\item[]{\bf Proof.}\ }{\hfill $\square$ \end{trivlist}\par}
\newenvironment{apf}[1]{\par\begin{trivlist}%
\item[]{\bf Proof of #1.}\ }{\hfill $\square$ \end{trivlist}\par}

\makeatletter \@addtoreset{equation}{section} \makeatother

\newcommand{\E}{\mathbb{E}}
\newcommand{\Q}{\mathbb{Q}}
\newcommand{\R}{\mathbb{R}}

\newcommand{\N}{\mathbb{N}}

\renewcommand{\d}{\mathrm{d}}
\newcommand{\m}{\mathfrak{m} }

\newcommand{\Exp}{\mbox{\rm Exp}}

\def\.{{\hspace{0.2mm}}}
\def\PP{{\bf P}}
\def\EE{{\bf E}}
\def\FF{{\bf F}}

\def\bJ{{\bf J}}
\def\rf{{\rm ref}}

\def\E{{\mathscr E}}
\def\F{{\mathscr F}}
\def\C{{\mathscr C}}

\def\H{{\mathscr H}}

\def\Q{{\mathscr Q}}

\def\1{{\bf 1}}
\def\eps{\varepsilon}
\def\loc{{\rm loc}}

\def\wh{\widehat}
\def\wt{\widetilde}

\def\<{{\langle}}
\def\>{{\rangle}}
\def\d{\rm d}

\newcommand{\Rd}{\mathbb{R}^d}

\renewcommand{\R}{\mathbb{R}}
\renewcommand{\N}{\mathbb{N}}

\renewcommand{\m}{\frak{m}}

\makeatletter

\@addtoreset{equation}{section}
\makeatother

\date{}

\title
{\large\bf Stability of estimates for fundamental solutions under Feynman-Kac perturbations for symmetric Markov processes}

\author{Daehong Kim\thanks{Department of Mathematics and Engineering,
Graduate School of Science and Technology,
Kumamoto University,
Kumamoto, 860-8555 Japan ({\sf daehong@gpo.kumamoto-u.ac.jp}). Supported in part by JSPS Grant-in-Aid for 
Scientific Research (KAKENHI) 17K05304.
},\ \  
Panki Kim\thanks{Department of Mathematical Sciences and Research Institute of Mathematics, Seoul National University, Building 27, 1 Gwanak-ro, Gwanak-gu
Seoul 08826, Republic of Korea ({\sf pkim@snu.ac.kr}). 
Supported in part by 
the  National  Research  Foundation  of  Korea  (NRF)  grant  funded  by  the  Korea government (MSIP)  (No.  2016R1E1A1A01941893). 
} \ \ and \ \
Kazuhiro Kuwae\thanks{Department of Applied Mathematics, Fukuoka University,
Fukuoka 814-0180, Japan ({\sf kuwae@fukuoka-u.ac.jp}). Supported in part by JSPS Grant-in-Aid for Scientific Research (KAKENHI) 17H02846 
and by fund (No.:185001) from the Central Research Institute of Fukuoka University.
}}

\date{\today}
\begin{document}

\maketitle 
\begin{abstract}  
In this paper, when a given symmetric Markov process {\bf X} satisfies  the stability of global 
heat kernel two-sided (upper) estimates by Markov perturbations (See Definition \ref{d:stHKE}),  
we give a necessary and sufficient condition on the stability of global two-sided (upper) estimates for fundamental solution of Feynman-Kac semigroup of {\bf X}.
As a corollary, under the same assumptions, a weak type global two-sided (upper) estimates holds for 
the fundamental solution of 
 Feynman-Kac semigroup with (extended) Kato class conditions for measures. 
This generalizes all known results on the stability of global integral kernel estimates by symmetric Feynman-Kac
perturbations with Kato class conditions in the framework of symmetric Markov processes. 
\end{abstract}



\section{Introduction and results}\label{sec:Result} 

Let $(E,d)$ be a locally compact separable metric space and $\m$ a Radon measure on $E$ with full topological support. Denote by $E_{\partial}$ the one point compactification of $E$. Let 
${\bf X}=(\Omega,\mathscr{F}_{\infty},\mathscr{F}_t,X_t, $ $\zeta, \PP_x, x\in E_{\partial}, t \in [0, +\infty[)$
 be an $\m$-symmetric 
 Markov process on $E_\partial$ 
with the life time  $\zeta:=
\inf\{t>0\mid X_t=\partial\}$. Let $(\E,\F)$ be a symmetric regular Dirichlet form on $L^2(E;\m)$ associated to {\bf X}. 

{\it Throughout the paper, we assume that $(\E,\F)$ is irreducible {\rm (}{\bf (I)} in short\,{\rm)}, that is, any $(T_t)$-invariant set is $\m$-trivial {\rm(}see \cite[\S 1.6]{FOT}\,{\rm)}, {\bf X} has 
 no killing inside $E$,
that is, $\PP_x(X_{\zeta-}\!\in E, \zeta < \infty)=0$ for any $x \in E$, and that {\bf X} satisfies 
 the absolute continuity condition {\rm(}{\bf (AC)} in short {\rm)}, that is, $P_t(x,{\rm d}y):=\PP_x(X_t\in{\rm d}y)\ll \m({\rm d}y)$ for any $t>0$ and $x\in E$. 
}

Let $p_t(x,y)$ be a function of $(t,x,y)\in ]\,0,\,+\infty\,[\times E\times E$. 
Throughout this paper, we say that $p_t(x,y)$ is a \emph{heat kernel} of the form $(\E,\F)$ if it satisfies the following (\cite[Section~3.3]{Grigo2014.2}): 
\begin{enumerate}
\item for any $t>0$, $p_t(x,y)$ is 
$\m\times \m$-measurable in $(x,y)\in E\times E$; 
\item for any $t>0$, $p_t(x,y)\geq0$ for $\m$-a.e.~$x,y\in E$, and $f\in\mathscr{B}_b(E)$ 
$$
\int_Ep_t(x,y)f(y)\m({\rm d}y)=\int_E f(y)\PP_x(X_t \in {\rm d}y)=:\EE_x[f(X_t)],
$$
for $\m$-a.e.~$x\in E$;
\item for any $t>0$, $p_t(x,y)=p_t(y,x)$ $\m$-a.e.~$x,y\in E$; 
\item  for all $t,s>0$, 
$$
p_{t+s}(x,y)=\int_Ep_t(x,z)p_s(z,y)\m({\rm d}z)
$$ 
for $\m$-a.e.~$x,y\in E$.
\end{enumerate}
Moreover, we say that $p_t(x,y)$ is a \emph{heat kernel in the strict sense} 
associated with the form $(\E,\F)$ if it is Borel measurable in $(t,x,y)$ and 
(2)--(4) above hold for all $x,y\in E$.  Under {\bf (AC)}, {\bf X} admits a heat kernel $p_t(x,y)$ in the strict sense {\rm(}see \cite[(2.1),(2.2)]{Mori:pKato}, \cite[Theorem~2]{Yan:1988}{\rm)}.  

Let $(\E,\F_e)$ be the extended Dirichlet space of $(\E,\F)$. Note that any element $f\in \F_e$ admits a strictly $\E$-quasi continuous version $\tilde{f}$ with $\tilde{f}(\partial)=0$ (\cite[Lemma~2.1]{Kw:stochI}).  
Throughout this paper, we always take a strictly $\E$-quasi continuous version of the element of $\F_e$, that is, we omit tilde from $\tilde{f}$ for $f\in\F_e$. 
Since {\bf X} has 
 no killing inside $E$,  by the Beurling-Deny formula (\cite[Theorem 3.2.1]{FOT}), $\E$ can be decomposed as follows: for $f, g \in \F$,
\begin{align*}
\E(f,g)&=\E^{(c)}(f,g) + \E^{(j)}(f,g) \\
&:=\E^{(c)}(f,g)+\int_{E\times E \setminus {\sf diag}}(f(x)-f(y))(g(x)-g(y))J({\rm d}x{\rm d}y)
\end{align*}
where $ \E^{(c)}$ is the strongly local part of $(\E, \F)$ and 
${\sf diag}:=\{(x,x)\mid x\in E\}$ is the diagonal set of $E \times E$.

This contains the particular case that $\E^{(c)}=0$, i.e., 
{\bf X} is a pure jump process, 
or that $J=0$, i.e., {\bf X} is a diffusion process.

In what follows, we consider $k\geq0$ and a subset $\mathbb{T}$ of $
[\,0,\,+\infty\,[^{\,2}:=[\,0,\,+\infty\,[\times[\,0,\,+\infty\,[$. 
For measurable functions $f$ on $[\,0,\,+\infty\,[\times E$ and $g$ on $[0,+\infty[^{\,2}$, 
we consider the function $\Psi$ on $[\,0,\,+\infty\,[\times E\times E$ by 
\begin{align}
\Psi_t(x,y):=f(t,x)g(t,d(x,y)).\label{eq:TypeFunc}
\end{align}
Let $\Phi$ be another measurable function defined on $[\,0,\,+\infty\,[\times E\times E$. 
We then write $\Phi_t(x,y)\precsim_k\Psi_t(x,y)$ 
(resp.~$\Phi_t(x,y)\succsim_k\Psi_t(x,y)$)
for $\m$-a.e.~$x,y$ with $d(x,y)\in\mathbb{T}_t$ if there exist $c_2,C_2>0$ (resp.~$c_1,C_1>0$) independent of $t\geq0$ such that 
\begin{align*}
\Phi_t(x,y)\leq C_2e^{kt}f(t,x)g(t,c_2d(x,y))\qquad\text{(resp.~$C_1e^{-kt}f(t,x)g(t,c_1d(x,y))\leq\Phi_t(x,y)$}\text{)}
\end{align*}
for $\m$-a.e.~$x,y\in E$ with $d(x,y)\in \mathbb{T}_t$. 
Here $\mathbb{T}_t:=\{s\in[\,0,\,+\infty\,[\; \mid (t,s)\in\mathbb{T}\}$ is the section of $\mathbb{T}$ at $t\in[\,0,\,+\infty\,[$. 
Moreover,  we write $\Phi_t(x,y)\asymp_k\Psi_t(x,y)$ for $\m$-a.e. $d(x,y)\in\mathbb{T}_t$ if $\Phi_t(x,y)\precsim_k\Psi_t(x,y)$ for $\m$-a.e. $d(x,y)\in\mathbb{T}_t$ and 
$\Phi_t(x,y)\succsim_k\Psi_t(x,y)$ for $\m$-a.e. $d(x,y)\in\mathbb{T}_t$. Similarly, $\Phi_t(x,y)\asymp_k\Psi_t(x,y)$ for $d(x,y)\in\mathbb{T}_t$ can be defined.
If we omit \lq\lq$\m$-a.e.\rq\rq, it means that each statement holds for all $x,y\in E$ with $d(x,y)\in \mathbb{T}_t$. 

When $k=0$, we simply write $\Phi_t(x,y)\precsim\Psi_t(x,y)$,  $\Phi_t(x,y)\succsim \Psi_t(x,y)$ and $\Phi_t(x,y)\asymp\Psi_t(x,y)$) instead of  $\Phi_t(x,y)\precsim_0\Psi_t(x,y)$, $\Phi_t(x,y)\succsim_0 \Psi_t(x,y)$ and  $\Phi_t(x,y)\asymp_0\Psi_t(x,y)$, respectively.

For measurable functions $G$ on $E\times E$ and $H$ on $[\,0,\,+\infty\,[$, we write $G(x,y)\precsim H(d(x,y))$ (resp.~$H(d(x,y))\precsim G(x,y)$) 
$\m$-a.e.~$x,y\in E$ if there exist $c_2,C_2>0$ (resp.~$c_1,C_1>0$) such that 
\begin{align*}
G(x,y)\leq C_2 H(c_2 d(x,y))\quad \text{(resp.~$C_1H(c_1d(x,y))\leq G(x,y)$}\text{)}\;\; \m\text{-a.e.}~x,y\in E.
\end{align*}
We write $G(x,y)\asymp H(d(x,y))$ $\m$-a.e.~$x,y\in E$ if $G(x,y)\precsim H(d(x,y))$ $\m$-a.e.~$x,y\in E$ and 
$ G(x,y)\succsim H(d(x,y))$ $\m$-a.e.~$x,y\in E$. 

Let $\phi$, $\psi$ be measurable functions defined on $[\,0,\,+\infty\,[\times E\times E$. We then write $\phi_t(x,y)\lesssim_k 
\psi_t(x,y)$ (resp,~$\phi_t(x,y)\gtrsim_k 
\psi_t(x,y)$) if there exists $C_2>0$ (resp.~$C_1>0$) independent of $t>0$ such that 
\begin{align*}
\phi_t(x,y)\leq C_2e^{kt}\psi_t(x,y)\qquad\text{(resp.~$C_1e^{-kt}\psi_t(x,y)\leq\phi_t(x,y)$}\text{)},\quad  \text{ for all }\quad t>0,\quad x,y\in E.
\end{align*}
We write $\phi_t(x,y)\simeq_k \psi_t(x,y)$ if both $\phi_t(x,y)\lesssim_k 
\psi_t(x,y)$  and $\phi_t(x,y)\gtrsim_k 
\psi_t(x,y)$ hold. When $k=0$, we simply write $\phi_t(x,y)\lesssim\psi_t(x,y)$ 
(resp. $\phi_t(x,y)\gtrsim \psi_t(x,y)$) instead of  $\phi_t(x,y)\lesssim_0\psi_t(x,y)$ (resp. $\phi_t(x,y)\gtrsim_0 \psi_t(x,y)$). $\phi_t(x,y)\simeq \psi_t(x,y)$ also means $\phi_t(x,y)\simeq_0 \psi_t(x,y)$. 
The reader should be aware of the difference between $\precsim_k$ and $\lesssim_k$.

We will consider
the following assumptions on the  
stability of the two-sided (upper) estimates for heat kernel $p_t(x,y)$ of ${\bf X}$. For this, let 
$(\wt{\E},{\F})$ be another symmetric Dirichlet form on 
$L^2(E;\wt{\m})$ with the same domain $\F$. 
Denote by 
$\wt{J}$ the jumping measure with respect to  
$(\wt{\E},{\F})$.  
We prepare the following conditions: 
There exists a constant $C_E>0$ such that 
\begin{align}
C_E^{-1}\m&\leq \wt{\m}\leq C_E\m, \label{eq:comparisonmeasure}\\
C_E^{-1}\E^{(c)}(f,f)&\leq \wt{\E}^{(c)}(f,f)\leq C_E\E^{(c)}(f,f)\quad\text{ for }f\in\F\label{eq:comparisondiffusion}
\end{align}
and,  there exist measurable functions $J$ and $\wt{J}$ on $E\times E$ such that
$J({\d}x{\d} y)=J(x,y)\m({\d} x)\m({\d} y)$ and $\wt{J}({\d}x{\d} y)=\wt{J}(x,y)\m({\d} x)\m({\d} y)$ satisfying 
\begin{align}
C_E^{-1}J(x,y)\leq\wt{J}(x,y)\leq C_EJ(x,y)\quad \m\text{-a.e.~}x,y\in E.\label{eq:comparisonjumpdensity}
\end{align}
We further remark the following:  
\begin{itemize}
\item 
\eqref{eq:comparisondiffusion} is equivalent to 
$C_E^{-1}\mu_{\<f\>}^c\leq \wt{\mu}_{\<f\>}^c\leq C_E\mu_{\<f\>}^c$ for $f\in \F$
where $\mu_{\<f\>}^c$ (resp.~$\wt{\mu}_{\<f\>}^c$) 
 is the energy measure of continuous part for $f\in{\F}$  with respect to  $({\E},{\F})$ (resp.~  $(\wt{\E},{\F})$)
(see \cite[Proposition~1.5.5(b)]{LeJan:Measure}, \cite{Ms:compmedia}). 
\item \eqref{eq:comparisonjumpdensity} 
implies that
\begin{align}
J(x,y)\asymp J^0(d(x,y))\quad\m\text{-a.e.~}x,y\in E\;\Longrightarrow\;
\wt{J}(x,y)\asymp J^0(d(x,y))\quad\m\text{-a.e.~}x,y\in E.  \label{eq:comparisonjumpasymp}
\end{align}
In fact, all results in this paper hold true under \eqref{eq:comparisonjumpasymp} instead of \eqref{eq:comparisonjumpdensity}.
\end{itemize}

Throughout this paper, we assume the existence of densities $J$ and $\wt{J}$. {Further, we make the following assumptions}. 
\begin{enumerate} 
\item[{{\bf (A.1)}}] 
{\rm (with $\mathbb{T}$)}
{Suppose that the heat kernel $p_t(x,y)$ on $]\,0,\,+\infty\,[\times E\times E$ in the strict sense associated with $(\E,\F)$ satisfies the following: there exists a function $\phi_2$ of the type \eqref{eq:TypeFunc} whose second factor defined on a subset $\mathbb{T}$ of $[\,0,\,+\infty\,[^{\,2}$ such that}  
\begin{align}
{p}_t(x,y)\precsim \phi_2(t,x,y)\quad\m\text{-a.e.}\quad 
d(x,y)\in\mathbb{T}_t.
\label{caparabilityUpperOriginalAsymp}
\end{align}
Then any Dirichlet form
$(\wt{\E},{\F})$
satisfying  \eqref{eq:comparisonmeasure}, \eqref{eq:comparisondiffusion} and \eqref{eq:comparisonjumpdensity}, 
 also admits a heat kernel 
$\wt{p}_t(x,y)$ on $]\,0,\,+\infty\,[\times E\times E$ such that 
\begin{align}
\wt{p}_t(x,y)\precsim \phi_2(t,x,y)\quad\m\text{-a.e.}\quad 
d(x,y)\in\mathbb{T}_t.
\label{caparabilityUpper}
\end{align}

\item[{\bf (A.2)}]{\rm (with $\mathbb{T}$)}
{Suppose that the heat kernel $p_t(x,y)$ on $]\,0,\,+\infty\,[\times E\times E$ in the strict sense associated with $(\E,\F)$ satisfies the following: there exist functions $\phi_1, \phi_2$ of the type \eqref{eq:TypeFunc} whose second factor defined on a subset $\mathbb{T}$ of $[\,0,\,+\infty\,[^{\,2}$ such that }
\begin{align}
\phi_1(t,x,y)\precsim p_t(x,y)\precsim \phi_2(t,x,y)\quad\text{ for }\quad d(x,y)\in\mathbb{T}_t.
\label{caparabilityOriginalasymp}
\end{align}
Then any Dirichlet form $(\wt{\E},{\F})$ satistying \eqref{eq:comparisonmeasure}, \eqref{eq:comparisondiffusion} and \eqref{eq:comparisonjumpdensity}
 also admits a heat kernel 
$\wt{p}_t(x,y)$ on $]\,0,\,+\infty\,[\times E\times E$ in the strict sense  such that 
\begin{align}
\phi_1(t,x,y)\precsim \wt{p}_t(x,y)\precsim \phi_2(t,x,y)\quad\text{ for }\quad d(x,y)\in\mathbb{T}_t.
\label{caparability}
\end{align}
\end{enumerate}

\begin{rem}\label{rem:TimeSet}
{\rm 
For example, if $\mathbb{T}=\{(t,s)\mid s\leq t\}$ 
(resp.~$\mathbb{T}=\{(t,s)\mid s\geq t\}$), then 
$\mathbb{T}_t=[\,0,\,t\,]$ 
(resp.~$\mathbb{T}_t=[\,t,\,+\infty\,[$). Note here that $\mathbb{T}_t=[\,0,\,+\infty\,[$ for any $t\in\,[\,0,\,+\infty\,[$ 
provided $\mathbb{T}=[\,0,\,+\infty\,[^{\,2}$, in this case, $d(x,y)\in\mathbb{T}_t$ automatically holds.
} 
\end{rem}

\begin{df}\label{d:stHKE}
{\rm 
We say that 
\emph{the stability of heat kernel estimate (resp.~heat kernel upper estimate) of $(\E,\F)$ (or the associated process ${\bf X}$)  globally holds
 by 
Markov perturbations} provided {\bf (A.2)}
(resp.~{\bf (A.1)})
is satisfied with $\mathbb{T}=[\,0,\,+\infty\,[^{\,2}$.
}
\end{df}

{\bf (A.1)} and {\bf (A.2)} are 
satisfied for quite wide class of  
symmetric Markov processes. For instance, diffusion processes on 
smooth complete Riemannian manifold with non-negative Ricci curvature, diffusion processes on fractals, symmetric stable-like 
processes and so on (see Section~\ref{sec:example}). 
Both {\bf (A.1)} and {\bf (A.2)} are 
stable under a Girsanov transformation (see after Theorem~\ref{FeynmanKacMod3} below). 
This property plays a key role in the proof of main results. 

Let $\mu_1$ (resp.~$\mu_2$) be a smooth measure in the strict sense corresponding to a positive 
continuous additive functional (PCAF in abbreviation) $A^{\mu_1}$ (resp.~$A^{\mu_2}$) in the strict sense 
with respect to ${\bf X}$  
and let $\mu:=\mu_1-\mu_2$, which is a signed smooth measure in the strict sense.  
We denote the corresponding 
continuous additive functional of $\mu=\mu_1-\mu_2$ by $A^\mu:=A^{\mu_1}-A^{\mu_2}$ to emphasize the correspondence between $\mu$ and $A$. 
Let $F_i$ $(i=1,2)$ be a non-negative bounded function 
on $E\times E_\partial$ symmetric on $E\times E$
which is extended to a function defined on $E_\partial \times E_\partial$ vanishing 
on the diagonal set ${\sf diag}$ of $E_\partial \times E_\partial$.
 We set $F:=F_1-F_2$. Then $A^F_t:=\sum_{0<s\le t}F(X_{s-},X_s)$ (whenever it is summable) is an additive functional of ${\bf X}$. 
Let ${\F}_{\rm loc}$ be the space of functions locally in $\F$ in the ordinary sense,  
$QC(E_{\partial})$ the space of 
strictly $\E$-quasi 
continuous  functions on $E_{\partial}$ (See Section~\ref{sec:preliminaries} for the definition)  
and  $S_D^1({\bf X})$ the class of smooth measures of Dynkin class with respect to ${\bf X}$
(see $\S2$ for definitions).

We consider a bounded finely continuous (nearly) Borel function 
$u\in{\F}_{\loc}\cap QC(E_{\partial})$ satisfying $\mu_{\<u\>}\in S_D^1({\bf X})$.
In \cite[Theorem~6.2(2)]{KimKuwaeTawara}, we proved that 
the additive functional $u(X_t)-u(X_0)$ admits the following strict decomposition: for all $t \in [\,0,\,+\infty\,[$
\begin{align}
u(X_t)-u(X_0)=M_t^u+N_t^u\quad {\PP}_x\text{-a.s.~for all }x\in E, \label{eq:Fukudecstrict} 
\end{align}
where $M^u$ is a square integrable martingale additive functional in the strict sense, and 
$N^u$ is a continuous additive functional (CAF in abbreviation) in the strict sense which is locally of zero energy. 
We note that $N^u$ is not a process of finite variation in general. 
$M^u$ can be decomposed as
\begin{align}
M_{t}^{u}=M_{t}^{u,c} + M_{t}^{u,j}\label{eq:FukuDecom}
\end{align}
where $M_{t}^{u,j}$ 
and $M_{t}^{u,c}$ are the jumping and continuous part of $M^{u}$ respectively. 
Those are defined $\PP_x$-a.s.~for all~$x\in E$ by \cite[Theorem~6.2(2)]{KimKuwaeTawara}. The strict decompositions \eqref{eq:Fukudecstrict} and \eqref{eq:FukuDecom} on $[\,0,\,+\infty\,[$ guarantee the extension of the supermartingale multiplicative functional 
 on $[\![0,\zeta[\![$ up to $[\,0,\,+\infty\,[$ 
(see \cite[Proposition~3.1]{KK:AnalChara}). Let $\mu_{\< u\>}$, $\mu^{c}_{\< u\>}$ and  $\mu^{j}_{\< u\>}$ 
be the smooth Revuz measures in the strict sense associated with the quadratic variational processes (or the sharp bracket PCAFs in the strict sense) $\< M^u\>$, $\< M^{u,c}\>$ and $\< M^{u,j}\>$ 
respectively. Then
$$
\mu_{\< u\>}({\rm d}x)=\mu^{c}_{\< u\>}({\rm d}x)+\mu^{j}_{\< u\>}({\rm d}x)
$$
We now 
consider the transforms by the additive functionals $A_t:=N_t^u+A_t^{\mu}+A_t^F$ of the form
\begin{align}
e_{A}(t):=\exp (A_t), \quad t \ge 0.     
\label{generalFKtansform}
\end{align}  
The transform \eqref{generalFKtansform} defines a semigroup, namely, the generalized Feynman-Kac semigroup 
\begin{align}
P_t^A\!f(x):=\EE_x[e_A(t)f(X_t)], \quad f \in \mathscr{B}_b(E), ~t\ge 0.  
\label{FKsemigroup}
\end{align}
Under some conditions, $(P_t^A)_{t\geq0}$ forms a strongly continuous semigroup on 
$L^2(E;\m)$ (see Lemma~\ref{lem:strongcontinuity} below). 
Since $N_t^u$ is time reversible, and $F_1$ and $F_2$ 
are symmetric on $E\times E$, $(P_t^A)_{t>0}$ is 
symmetric in the sense that 
$$
(P_t^Af,g)_{\m}=(f,P_t^Ag)_{\m}\quad\text{ for }f,g\in  \mathscr{B}_+(E)
$$
where $(f,g)_{\m}:=\int_Ef(x)g(y)\m({\rm d}x).$
We say that $(P_t^A)_{t>0}$ admits an \emph{integral kernel} $p_t^A(x,y)$ on $]\,0,\,+\infty\,[\times E \times E$ if it satisfies 
the following: 
\begin{enumerate}
\item for any $t>0$, $p_t^A(x,y)$ is $\m\times \m$-measurable in $(x,y)\in E\times E$; 
\item for any $t>0$, $p_t^A(x,y)\geq0$ for $\m$-a.e.~$x,y\in E$, and 
for $f\in\mathscr{B}_b(E)$
$$
\int_Ep_t^A(x,y)f(y)\m({\rm d}y)=\EE_x[e_A(t)f(X_t)],
$$
for $\m$-a.e.~$x\in E$;
\item for any $t>0$, $p_t^A(x,y)=p_t^A(y,x)$ $\m$-a.e.~$x,y\in E$; 
\item  for all $t,s>0$, 
$$
p_{t+s}^A(x,y)=\int_Ep_t^A(x,z)p_s^A(z,y)\m({\rm d}z)
$$ 
for $\m$-a.e.~$x,y\in E$.
\end{enumerate}
We say that $p_t^A(x,y)$ is an \emph{integral kernel of  $(P_t^A)_{t>0}$ in the strict sense} if it is Borel measurable in $(t,x,y)$ and 
(2)--(4) above hold for all $x,y\in E$.  

Stability of the short or global time estimate for an integral kernel of Feynman-Kac semigroup of various Markov processes has been studied by many authors (\cite{CKS:Stability, KimKuwae:stability, Song06, Tak:ultra, Tak:GaussianHK, Wada, Wang}). Takeda \cite{Tak:GaussianHK} studied the global stability of \emph{Li-Yau estimate} for the integral kernel under local Feynman-Kac perturbations and gave a necessary and sufficient condition for the global stability of Li-Yau estimate in an analytic way. More precisely, let 
${\bf B}=(B_{t}, \PP_{x})$ 
be a Brownian motion on a complete
non-compact Riemannian
 manifold $(M,g)$ and assume the heat kernel $p_t(x,y)$ associated with $\frac{\,1\,}{2}\Delta_M$, the half of Laplace-Beltrami operator on $(M,g)$, satisfies the following Li-Yau estimate (cf. \cite{Grigo2006})
\begin{equation}
\frac{C_1}{\m(B(x,\sqrt{t}))} 
\exp\left(-c_1\frac{d^2(x,y)}{t}\right)
\le p_t(x,y) \le 
\frac{C_2}{\m(B(x,\sqrt{t}))} 
\exp\left(-c_2\frac{d^2(x,y)}{t}\right), 
\label{LiYau*}
\end{equation}
where $C_1, c_1, C_2$ and $c_2$ are positive constants and $\m(B(x,r))$ denotes the volume of the geodesic ball of radius $r$ centered at $x \in M$. 
 For a signed smooth measure $\hat{\mu}=\hat{\mu}_+-\hat{\mu}_-$,  
 let $p_t^{\hat{\mu}}(x,y)$ be the integral kernel associated with the Schr\"odinger operator $\frac{\,1\,}{2}\Delta_M +\hat{\mu}$. Then the necessary and sufficient condition that $p_t^{\hat{\mu}}(x,y)$ also satisfies \eqref{LiYau*} is 
 $\lambda(\hat{\mu})>1$ where
\begin{align}
\lambda(\hat{\mu})=\inf\left\{\frac{\,1\,}{2}\int_M\<{\nabla f}, {\nabla f}\>\,{\rm d}\m+\int_Mf^2{\rm d}\hat{\mu}^- \;\Biggl|\; f \in C_0^{\infty}(M), \int_Mf^2{\rm d}\hat{\mu}^+=1\right\}, 
\label{Acondtion}
\end{align}
provided $\hat{\mu}_+\in S_{C\!S_{\infty}}^1({\bf B})$ and $\hat{\mu}_-\in S_{D\!S_0}^1({\bf B})$.
 Here $S_{C\!S_{\infty}}^1({\bf B})$ (resp.~$S_{D\!S_0}^1({\bf B})$) denotes 
the class of conditionally Green-tight smooth 
measures of Kato class (resp.~the class of conditionally Green bounded measures)  with respect to $B$
(see \cite[Definition~6.1]{KimKurKuwae:Gauge}). 
Note that the analytic condition such as \eqref{Acondtion} is very useful in many concrete cases for confirming the global stability of 
an integral kernel under Feynman-Kac perturbations. 

The purpose of this paper is to give the analytic condition on $u, \mu$ and $F$ under which the integral kernel $p_t^A(x,y)$ of the Feynman-Kac semigroup $(P_t^A)_{t>0}$ satisfies 
\eqref{caparabilityUpper} (resp.~\eqref{caparability}) under 
{\bf (A.1)}
(resp.~{\bf (A.2)} 
). 
Before stating our main results, let us introduce a spectral function determined by a quadratic form and a measure. 
Let $ \F_b:=\F\cap L^{\infty}(E;\m)$.
For given $u, \mu$ and $F$ in the above defintion of the additive functionals $A_t=N_t^u+A_t^{\mu}+A_t^F$, 
let $(\Q,\F_b)$ be a quadratic form on $L^2(E;\m)$ defined by \begin{align}
\Q(f,g):={\E}(f,g)+{\E}(u,fg) - \H(f,g), \quad 
f, g \in \F_b,
\label{eq:quadraticform}
\end{align}
where
\begin{align*}
{\H(f,g)}:=\int_Ef(x)g(x)\mu({\rm d}x)+\int_E\int_Ef(x)g(y)\left(e^{F(x,y)}-1\right)N(x,{\rm d}y)\mu_H({\rm d}x).
\end{align*}
It is known that $(\Q,\F_b)$ is well-defined under $\mu_{\<u\>}+\mu_1+\mu_2+N(F_1+F_2)\mu_H\in S_{D}^1({\bf X})$ (\cite{KK:AnalChara}
). Set 
$\Q_{\alpha}(f,g):=\Q(f,g)+\alpha(f,g)_{\m}$ for $\alpha\geq0$ and $f,g\in\F_b$. Clearly,  $\Q(f,g)=\Q_0(f,g)$.
For a measure $\nu\in S_D^1({\bf X})$, define the spectral function by 
\begin{align}
\lambda^{\Q_{\alpha}}(\nu):=\inf\left\{{\Q}_{\alpha}(f,f) \;\Biggl|\; f \in \C,
 ~\int_{E}f^{2}{\rm d}\nu=1\right\}.
\label{eq:eigenvalueexpress1}
\end{align}
Here $\C$ stands for the special standard core of $(\E, \F)$ (\cite{FOT}). Note that $(\Q_{\alpha},\F_b)$ is not necessarily non-negative definite and $\nu\mapsto 
\lambda^{\Q}(\nu)$ is {non-increasing} in the sense that $\nu_1\leq\nu_2$ implies 
$\lambda^{\Q}(\nu_1)\geq \lambda^{\Q}(\nu_2)$.
Define the signed measures $\overline{\mu}:=\overline{\mu}_1-\overline{\mu}_2$ and 
$\overline{\mu}^*:=\overline{\mu}_1^*-\overline{\mu}_2^*$
 by 
\begin{align}
&\overline{\mu}_1:=N\left(e^{F+[u]} \!\!- (F+[u]) \!-1+F_1\right)\mu_H+\mu_{1}+\frac{\,1\,}{2}\mu_{\<u\>}^c, \quad \overline{\mu}_2:=N(F_{2})\mu_H+\mu_{2}, \label{e:mu12} \\
&\overline{\mu}_1^*:=N(e^{F_1+[u]}-(F_1+[u])-1+F_1)\mu_H+\mu_{1}+\frac{\,1\,}{2}\mu_{\<u\>}^c, \qquad \overline{\mu}_2^*:=\mu_{2}, \label{e:mu12*} 
\end{align}
where $[u](x,y):=u(x)-u(y)$. 
Note that $\overline{\mu}_1^*=\mu_1+N(e^{[u]}(e^{F_1}-1))\mu_H+
N(e^{[u]}-[u]-1)\mu_H+\frac{\,1\,}{2}\mu_{\<u\>}^c$. 

Let us denote by $S_{N\!K_1}^1({\bf X})$ the class of natural semi-Green-tight smooth measures of extended Kato class,
$S_{N\!K_{\infty}}^1({\bf X})$ the class of natural Green-tight smooth measures of Kato class, and 
$S_{D_0}^1({\bf X})$ the class of Green-bounded smooth measures in the strict sense 
(see Definitions \ref{d:S_1} and \ref{df:newKatoclass}  for precise definitions of these). 

Now, we are ready to state our results:

\begin{thm}\label{thm:mainresult1}
Suppose that {\bf X} is transient. 
Let $u\in {\F}_{\loc}\cap QC(E_{\partial})$ be a bounded finely continuous (nearly) Borel function on $E$. Assume that $\mu_1+N(e^{F_1}-1)\mu_H
\in S_{N\!K_1}^1({\bf X})$, $\mu_{\<u\>}\in S_{N\!K_{\infty}}^1({\bf X})$ and   
$\mu_2+N(F_2)\mu_H\in S_{D_0}^1({\bf X})$ hold. Then we have the following: 
\begin{enumerate}
\item\label{item:mainresult1} 
If $\lambda^{\Q}(\overline{\mu}_1) >0$ and {\bf (A.1)} 
holds with $\mathbb{T}$,
then there exists an integral kernel $p_t^A(x,y)$ of the Feynman-Kac semigroup $(P_t^A)_{t>0}$ such that $p_t^A(x,y)\precsim \phi_2(t,x,y)$ $\m$-a.e. $d(x,y)\in\mathbb{T}_t$.

\item\label{item:mainresult1*} 
If $\lambda^{\Q}(\overline{\mu}_1) >0$ and {\bf (A.2)}
holds with $\mathbb{T}$, then there exists an integral kernel $p_t^A(x,y)$ of the Feynman-Kac semigroup $(P_t^A)_{t>0}$ in the strict sense such that 
\begin{align}
\phi_1(t,x,y) \precsim
p_t^A(x,y)\precsim \phi_2(t,x,y),\quad d(x,y)\in\mathbb{T}_t.\label{eq:FinalULestimate}
\end{align}

\item\label{item:mainresult2}  
If there exist an integral kernel $p_t^A(x,y)$ of the Feynman-Kac semigroup $(P_t^A)_{t>0}$ in the strict sense and a heat kernel $p_t(x,y)$ of {\bf X} in the strict sense 
such that $p_t^A(x,y) \lesssim p_t(x,y)$, then $\lambda^{\Q}(\overline{\mu}_1) >0$, in particular, under {\bf (A.2)}
with $\mathbb{T}=[\,0,\,+\infty\,[^{\,2}$ and $\phi_1(t,x,y)\simeq\phi_2(t,x,y)$, 
$p_t^A(x,y) \lesssim \phi_2(t,x,y)$ implies $\lambda^{\Q}(\overline{\mu}_1) >0$. 
\end{enumerate}
\end{thm}


\begin{rem}\label{rem:mainresult1}
{\rm \begin{enumerate}
\item\label{item:mainresult(i)} 
Under $\mu_{\<u\>}\in S_{N\!K_{\infty}}^1({\bf X})$, 
the condition $\mu_1+N(e^{F_1}-1)\mu_H\in 
S_{N\!K_1}^1({\bf X})$  in Theorem~\ref{thm:mainresult1} is equivalent to $\overline{\mu}_1^*\in  S_{N\!K_1}^1({\bf X})$ 
(see the proof of \cite[Lemma~3.1]{KimKuwae:general}), consequently, this is also equivalent to 
$\mu_1+N(e^{[u]}(e^{F_1}-1))\mu_H\in S_{N\!K_1}^1({\bf X})$, because 
$\overline{\mu}_1^*=\mu_1+N(e^{[u]}(e^{F_1}-1))\mu_H+N(e^{[u]}-[u]-1)\mu_H+\frac{\,1\,}{2}\mu_{\<u\>}^c$. 
\end{enumerate}
}
\end{rem} 

Theorem~\ref{thm:mainresult1}
extends the result on the stability of the Li-Yau estimate for the 
heat kernel of a Brownian motion on Riemannian manifold under Feynman-Kac perturbations proved by Takeda~\cite{Tak:GaussianHK}. 
The first main
contribution in our result is to add the perturbation by continuous additive functional of locally of zero energy. 
The global integral kernel estimate
 under such a perturbation was firstly discussed by Glover-Rao-Song~\cite[Theorem~2.9]{GRS1992} (see also Glover-Rao-$\check{\rm S}$iki\'c-Song 
\cite[Proposition~1.4]{GRSS1994}) 
in the framework of Brownian motion. But the stability  of global heat kernel estimate has not been treated in this direction. Indeed, the global estimates shown in \cite[Theorem~2.9]{GRS1992} is weaker than 
the usual Gaussian estimate like Brownian motion. 
The second main contribution  in our result
is the relaxation of the class of measures in the creation part. 
In Theorem~\ref{thm:mainresult1}\eqref{item:mainresult1}, we adopt the class of natural semi-Green-tight measures of extended Kato class in creation part 
instead of the 
Green-tight measures of Kato class. 
So 
Theorem~\ref{thm:mainresult1}\eqref{item:mainresult1} under 
{\bf (A.1)}
 (more precisely Theorem~~\ref{applthm:mainresult1} with Remark~\ref{rem:Devyver}) 
also covers the result by Devyver \cite[Theorem~4.1]{Devyver:HKRiesz}. 
The third main contribution 
is the relaxation 
 from the conditionally (semi-)Green-tight measures of (extended) Kato class (resp.~the class of conditionally Green bounded measures) into  (semi-)Green-tight measures of (extended) Kato class (resp.~the class of Green bounded measures) in Theorem~\ref{thm:mainresult1}\eqref{item:mainresult2}. 
Our Theorem~\ref{thm:mainresult1} improves Takeda's result \cite[Theorem~2]{Tak:GaussianHK} so that 
$p_t^{\mu}(x,y)$ satisfies Li-Yau estimate if and only if 
\eqref{Acondtion} holds for $\mu_+\in S_{N\!K_1}^1({\bf X})$ and 
$\mu_- \in S_{D_0}^1({\bf X})$, which
is weaker than the assertion in  Devyver~\cite[Theorem~4.1]{Devyver:HKRiesz}. 

The proof of Theorem~\ref{thm:mainresult1}  \eqref{item:mainresult1}\eqref{item:mainresult1*} is constituted of several steps: 
First, by using a Girsanov transform 
in terms of $u$, we reduce the proof to the case $u=0$.  
Second, by using another Girsanov transform 
in terms of $F_1,F_2$, we reduce the proof to the local 
perturbation under the transformed process ${\bf Y}$.

The proof of Theorem~\ref{thm:mainresult1}\eqref{item:mainresult2} is based on Theorem~\ref{thm:analgauge1} below, the equivalences among gaugeability, positivity 
of the bottom of spectrum and subcriticality. These equivalences had been regarded 
to hold by way of conditional gaugeability for Feynman-Kac functionals (see \cite{Chen:gaugeability2002,Tak:SubConGaugeabilityJFA2002}). 
In \cite{KK:AnalChara, KimKurKuwae:Gauge}, we prove that these equivalences hold without showing the conditional gaugeability, but showing the semi-conditional gaugeability under conditional (semi-)Green-tightness of 
related measures. Theorem~\ref{thm:analgauge1} does not require   the conditional (semi-)Green-tightness of measures and
 it is based on the characterization of the gaugeability under the 
Green-boundedness of underlying symmetrizing measure  (\cite[Lemma~4.3(5)]{KimKurKuwae:Gauge}) and 
a time change method.        

\medskip

The analogous statements for Theorem \ref{thm:mainresult1} without assuming the transience of {\bf X} are given 
in the next corollary, which plays a crucial role in the proof of Theorem \ref{thm:shortstability} below. For $\alpha >0$, let us denote ${\bf X}^{(\alpha)}$ by the $\alpha$-subprocess killed at rate $\alpha \m$.

\begin{cor}\label{cor:mainresult}
Let $u\in {\F}_{\loc}\cap QC(E_{\partial})$ be a bounded finely continuous (nearly) Borel function on $E$. Assume 
$\mu_1+N(e^{F_1}-1)\mu_H
\in S_{N\!K_1}^1({\bf X}^{(\alpha)})$, $\mu_{\<u\>}\in S_{N\!K_{\infty}}^1({\bf X}^{(\alpha)})$ and   
$\mu_2+N(F_2)\mu_H\in S_{D}^1({\bf X})$ 
for a fixed $\alpha>0$.  Then we have the following: 
\begin{enumerate}
\item\label{item:cormainresult1} 
If $\lambda^{\Q_{\alpha}}(\overline{\mu}_1) >0$ and {\bf (A.1)}
holds with $\mathbb{T}$,
then there exists an integral kernel $p_t^A(x,y)$ of the Feynman-Kac semigroup $(P_t^A)_{t>0}$ such that $p_t^A(x,y)\precsim_{k} \phi_2(t,x,y)$ $\m$-a.e. $d(x,y)\in\mathbb{T}_t$ 
for some constant $k:=k(\alpha)\geq0$ depending on $\alpha$.

\item\label{item:cormainresult1*}
If $\lambda^{\Q_{\alpha}}(\overline{\mu}_1) >0$ and 
{\bf (A.2)}
holds with $\mathbb{T}$, then there exists an integral kernel $p_t^A(x,y)$ of the Feynman-Kac semigroup $(P_t^A)_{t>0}$ in the strict sense such that 
\begin{align*}
\phi_1(t,x,y)
\precsim_{k}
p_t^A(x,y)\precsim_{k} \phi_2(t,x,y), \quad d(x,y)\in\mathbb{T}_t
\end{align*} 
for some constant $k:=k(\alpha)\geq0$ depending on $\alpha$.

\item\label{item:cormainresult2}  
If there exist a kernel $p_t^A(x,y)$ of the Feynman-Kac semigroup $(P_t^A)_{t>0}$ in the strict sense 
and a heat kernel $p_t(x,y)$ of {\bf X} in the strict sense  
such that $p_t^A(x,y)\lesssim_{k} p_t(x,y)$, then $\lambda^{\Q_{\alpha}}(\overline{\mu}_1) >0$ holds for $\alpha > k$, in particular, under 
{\bf (A.2)}
with $\mathbb{T}=[\,0,\,+\infty\,[^{\,2}$ and  
$\phi_1(t,x,y)\simeq_{\ell}\phi_2(t,x,y)$ for some $\ell\geq0$, 
$p_t^A(x,y) \lesssim_k \phi_2(t,x,y)$ implies $\lambda^{\Q_{\alpha}}(\overline{\mu}_1) >0$ for $\alpha > k+\ell$. 
\end{enumerate}
\end{cor}

\begin{rem}\label{rem:mainresult}
{\rm 
\begin{enumerate}
\item 
Under $\mu_{\<u\>}\in S_{N\!K_{\infty}}^1({\bf X}^{(\alpha)})$, 
the condition $\mu_1+N(e^{F_1}-1)\mu_H\in S_{N\!K_1}^1({\bf X}^{(\alpha)})$ in Corollary~\ref{cor:mainresult} is equivalent to 
$\overline{\mu}_1^*\in S_{N\!K_1}^1({\bf X}^{(\alpha)})$
(see the proof of \cite[Lemma~3.1]{KimKuwae:general}), consequently, this is also equivalent to 
$\mu_1+N(e^{[u]}(e^{F_1}-1))\mu_H\in S_{N\!K_1}^1({\bf X}^{(\alpha)})$, because $\overline{\mu}_1^*=\mu_1+N(e^{[u]}(e^{F_1}-1))\mu_H+N(e^{[u]}-[u]-1)\mu_H+\frac{\,1\,}{2}\mu_{\<u\>}^c$. 
\end{enumerate}
}
\end{rem}

The proof of Corollary~\ref{cor:mainresult}\eqref{item:cormainresult1}\eqref{item:cormainresult1*} (resp.~Corollary~\ref{cor:mainresult}\eqref{item:cormainresult2}) is quite similar to the proof of Theorem~\ref{thm:mainresult1}\eqref{item:mainresult1}\eqref{item:mainresult1*} (resp.~Theorem~\ref{thm:mainresult1}\eqref{item:mainresult2}).   However, the calculation in the proof of Corollary~\ref{cor:mainresult}\eqref{item:cormainresult1}\eqref{item:cormainresult1*} contains  extra factors yielding exponential growth and decay (see \eqref{epsrelation}). 

\medskip

We remark that from any Green-bounded measures a natural 
semi-Green-tight measure of extended 
Kato class  for $\alpha$-subprocess
can be constructed by  scaling (see Lemma~\ref{lem:extended} below). From this  observation, we can obtain the following result as a consequence of Corollary~\ref{cor:mainresult}.
See Definition \ref{d:S_1} for the definitions of $S_{E\!K}^1({\bf X})$ and  $S_{K}^1({\bf X})$.

\begin{thm}\label{thm:shortstability}
Let $u\in  {\F}_{\loc}\cap QC(E_{\partial})$ be a bounded finely continuous (nearly) Borel function on $E$. 
Assume $\mu_1+N(e^{F_1}-1)\mu_H\in S_{E\!K}^1({\bf X})$, $\mu_{\<u\>} \in S_{K}^1({\bf X})$ 
and $\mu_2 +N(F_2)\mu_H \in S_{D}^1({\bf X})$. 
Then we have the following: 
\begin{enumerate}
\item\label{item:shortstability1} Under 
{\bf (A.1)}
with $\mathbb{T}$, 
there exists an integral kernel $p_t^A(x,y)$ of the Feynman-Kac semigroup $(P_t^A)_{t>0}$ such that $p_t^A(x,y) \precsim_{k} \phi_2(t,x,y)$ $\m$-a.e. $d(x,y)\in\mathbb{T}_t$ for some $k\geq0$. 
\item\label{item:shortstability1*} 
Under {\bf (A.2)}
with $\mathbb{T}$, 
there exists an integral  kernel $p_t^A(x,y)$ of the Feynman-Kac semigroup $(P_t^A)_{t>0}$ 
in the strict sense such that 
$$
\phi_1(t,x,y)\precsim_{k} p_t^A(x,y) \precsim_{k}\phi_2(t,x,y),\quad d(x,y)\in\mathbb{T}_t
$$ for 
some 
$k\geq0$.
\end{enumerate}
\end{thm}

\begin{rem}\label{rem:shortstability}
{\rm 
\begin{enumerate}
\item 
Under $\mu_{\<u\>}\in S_{K}^1({\bf X})$, 
the condition $\mu_1+N(e^{F_1}-1)\mu_H\in S_{E\!K}^1({\bf X})$ 
in Theorem~\ref{thm:shortstability} is equivalent to 
$\overline{\mu}_1^*\in S_{E\!K}^1({\bf X}^{(\alpha)})$
(see the proof of \cite[Lemma~3.1]{KimKuwae:general}), consequently, this is also equivalent to 
$\mu_1+N(e^{[u]}(e^{F_1}-1))\mu_H\in S_{E\!K}^1({\bf X})$, 
because 
$\overline{\mu}_1^*=\mu_1+N(e^{[u]}(e^{F_1}-1))\mu_H+N(e^{[u]}-[u]-1)\mu_H+\frac{\,1\,}{2}\mu_{\<u\>}^c$. 
\end{enumerate}
}
\end{rem}

Theorem~\ref{thm:shortstability} also extends the previous known results on the integral kernel estimates under perturbation by 
measures of Kato classes in the framework of symmetric Markov processes.  
The conditions for measures 
in Theorem~\ref{thm:shortstability} are milder than those for    
known results (cf. \cite{Song06} and reference therein).  

\medskip

The rest of the paper is organized as follows. 

In Section~\ref{sec:preliminaries}, we collect basic terminologies and fundamental facts, and define several classes on Green-tight measures of Kato class.
In Section~\ref{sec:AofGandSubcriticality}, we summarize and develop the results in \cite{KK:AnalChara, KimKurKuwae:Gauge} in the framework of this paper. Section~\ref{sec:Stability results} is devoted to prove our main results. The main ingredients used in the proofs are in the front of this section.  In Section~\ref{sec:example}, we look our results more closely in the cases of symmetric diffusion process and symmetric jump process. 
\medskip 

In this paper, we use $c$ and $C$ as positive constants which may be different at different occurrences. 
In this paper, we use the following notations: For $a,b\in\R$, $a\lor b:=\max\{a,b\}$, $a\land b:=\min\{a,b\}$. $E_{\partial}$ is the 
one point compactification of $E$ where $\partial$ is the cemetery point for processes.
We use  $QC(E_{\partial})$ to denote the family of all strictly $\E$-quasi continuous functions on $E_{\partial}$, $C_0(E)$ 
(resp.~$C_\infty(E)$)
to denote the family of continuous functions on $E$ with compact support 
(resp.~ vanishing at infinity).
$\mathscr{B}_+(E)$ (resp.~$\mathscr{B}_b(E)$) to denote the space of non-negative (resp.~bounded) Borel functions on $E$, and $C_b(E)$ to denote
the space of bounded continuous functions on $E$.
For a non-negative Borel measure $\nu$  and a kernel $K$ on $E$, we write $K\nu(x):=\int_{E}K(x,y)\nu({\rm d}y)$. 
We also use the notation $Kf(x)=K\nu(x)$ when $\nu({\rm d}x)=f(x){\rm d}x$ for any $f\in\mathscr{B}_+(E)$ or $f\in\mathscr{B}_b(E)$. 
For a Borel subset $F$ of $E$, 
 $\tau_{F}:=\inf\{t>0\mid X_t\notin F\}$ (resp.~
$\sigma_{F}:=\inf\{t>0\mid X_t\in F\}$) is the first 
exit time of $X_t$ from $F$ (resp.~hitting time of $X_t$ to 
$F$). 
$\F_{F}:=\{u\in\F\mid u=0\;\m\text{-a.e.~on }E\setminus F\}$.

\section{Preliminaries and Green-tight measures of Kato class}\label{sec:preliminaries}

Under {\bf (AC)}, 
for $\alpha >0$, we can define the $\alpha$-order resolvent kernel
$R_{\alpha}(x,y)(<\infty)$ for $x,y \in E$ (see \cite[Lemma~4.2.4]{FOT}). 
We can also define $0$-order resolvent kernel $R(x,y):=R_0(x,y):=\lim_{\alpha\to0}R_{\alpha}(x,y)\leq\infty$ for $x,y \in E$, 
which is also called the \emph{Green kernel} of ${\bf X}$ provided ${\bf X}$ is transient. 
Recall that we use notations $R_{\alpha}\nu(x)=\int_ER_{\alpha}(x,y)\nu({\d}y)$ for $\alpha>0$ and $R\nu(x)=\int_ER(x,y)\nu({\d}y)$ where 
 $\nu$ is a Borel measure on $E$.

 An increasing sequence $\{F_n\}$ of closed sets is said to be an \emph{$\E$-nest} (resp.~\emph{strict $\E$-nest}) if $\PP_x(\lim_{k\to\infty}\tau_{F_n}=\zeta)=1$ 
(resp.~$\PP_x(\lim_{k\to\infty}\sigma_{F_n^c}=\infty)=1$) $\m$-a.e.~$x\in E$. It is shown in \cite{MR} that for an increasing sequence $\{F_n\}$ of closed sets, $\{F_n\}$ is an $\E$-nest if and only if $\bigcup_{n=1}^{\infty}\F_{F_n}$ is 
$\E_1^{1/2}$-dense in $\F$
where $\E_{1}(\cdot,\cdot):=\E(\cdot,\cdot)+(\cdot,\cdot)_\m$, 
equivalently 
$\lim_{n\to\infty}\text{\rm Cap}(K\setminus F_n)=0$ for any compact set $K$ in view of 
\cite[Lemma~5.1.6]{FOT}.
The regularity of the 
given Dirichlet form $(\E,\F)$ for ${\bf X}$ tells us that there always exists an $\E$-nest of compact sets (\cite[Chapter IV 4 (a), Chapter V, Proposition~2.12]{MR}). 
Hence any $\E$-nest can be taken to be a sequence of compact sets. 
A function $f$ defined on $E$ (resp.~$E_{\partial}$) is said to be \emph{$\E$-quasi continuous} (resp.~\emph{strictly $\E$-quasi continuous}) if 
there exists an $\E$-nest (resp.~a strict $\E$-nest) $\{F_n\}$ of closed sets such that $f|_{F_n}$ (resp.~$f|_{F_n\cup\{\partial\}}$) is continuous for each $k\in\mathbb{N}$.  
Recall that $QC(E_{\partial})$ is the family of all strictly $\E$-quasi continuous functions on $E_{\partial}$.

A positive Radon measure $\nu$ on $E$ is said to be \emph{of finite energy integral} if there exists $C>0$ such that 
\begin{equation}
\int_{E}|v(x)|\nu({\rm d}x)\leq C\sqrt{\E_{1}(v,v)}\quad \text{ for all }v\in\F\cap C_0(E), 
\label{S0}
\end{equation}
Denote by 
$S_{0}({\bf X})$  the family of measures of finite energy integrals and 
$S_{00}({\bf X}):=\{ \nu \in S_{0}({\bf X})\mid \nu(E)<\infty, \|R_1\nu\|_\infty<\infty
\}.$
A positive Radon measure $\nu$ on $E$ is said to be \emph{of $0$-order finite energy integral} if {\bf X} is transient and there exists 
$C>0$ satisfying (\ref{S0}) with $\E(\cdot,\cdot)$ in place of $\E_{1}(\cdot,\cdot)$, and we denote by $S_0^{(0)}({\bf X})$ the family of measures of $0$-order finite energy integrals. 
A  positive Borel measure $\nu$ on $E$ is said to be \emph{smooth} if $\nu$ charges no 
exceptional set and 
there exists an $\E$-nest $\{F_n\}_{n \ge 1}$ of compact sets satisfying $\nu(F_n)<\infty$ 
for each $n \ge 1$. 
It is known that $\nu$ is smooth if and only if $\nu$ charges no 
exceptional set and there exists an 
$\E$-nest of compact sets $\{F_n\}$ such that  $\1_{F_n} \nu \in S_{00}({\bf X})$ 
for each $n \ge 1$ (\cite[cf.~Theorem~2.2.4]{FOT}). 
We denote by $S({\bf X})$ the family of all smooth measures. $S({\bf X})$ contains all positive Radon measures on $E$ charging no set of zero capacity. Let $S_1({\bf X})$ be the family of positive smooth measures in the strict sense (\cite[Page 238]{FOT}).

\begin{df}\label{d:S_1}{\rm
(1)
A measure $\nu \in S_1({\bf X})$ is said to be \emph{of Dynkin class} (resp.~\emph{Green-bounded}) \emph{with respect to ${\bf X}$} if
  $\sup_{x\in E}R_{\alpha}\nu(x)<\infty$ for some $\alpha>0$ (resp.~$\sup_{x\in E}R\nu(x)<\infty$). 
 Denote by $S_D^1({\bf X})$ (resp.~$S_{D_0}^1({\bf X})$) the family of measures of Dynkin class (resp.~of Green-bounded measures)

\noindent (2)
 A measure $\nu \in S_1({\bf X})$ is said to be \emph{of Kato class (resp.~of extended Kato class) with respect to  ${\bf X}$} if 
  $\lim_{\beta\to\infty} \sup_{x\in E}R_{\beta}\nu(x)=0$ (resp.~$\lim_{\beta\to\infty} \sup_{x\in E}R_{\beta}\nu(x)<1$). 
 Denote by  $S_K^1({\bf X})$ (resp.~$S_{E\!K}^1({\bf X})$) the family of measures of Kato class (resp.~of extended Kato class).
 }
\end{df}
  
When {\bf X} is transient, denote by 
$S_{00}^{(0)}({\bf X}):=\{\nu\in S_0^{(0)}({\bf X})\mid 
  \nu(E)<\infty\text{ and }U\nu\in L^{\infty}(E;\m)\}$, the family of finite measures of $0$-order energy integrals with bounded $0$-order potentials. Here 
  $R\nu\in\F_e$ is the $0$-order potential for $\nu\in S_0^{(0)}({\bf X})$ defined by 
  $\E(R\nu,v)=\int_Ev{\rm d}\nu$ for $v\in \F\cap C_0(E)$. 
  By definition, we see $S_{00}^{(0)}({\bf X})\subset S_{D_0}^1({\bf X})$ and $S_{00}({\bf X})\subset S_D^1({\bf X})$ (\cite[cf.~Theorems~5.1.6 and 5.1.7]{FOT}). 
  Note that 
 any measure $\nu\in S_D^1({\bf X})$ is a positive Radon measure in view of Stollmann-Voigt's inequality: $
 \int_{E}u^2 d\nu\leq\|R_{\alpha}\nu\|_{\infty}\E_{\alpha}(u,u)$, $u\in\F,\alpha\geq0$ 
 (\cite[Theorem~3.1]{SV:potential}).   
 Conversely, any positive Radon measure $\nu$ satisfying $\sup_{x\in E}R_{\alpha}\nu(x)<\infty$ for some $\alpha>0$ 
  always belongs to $S_1({\bf X})$ in view of \cite[Proposition 3.1]{KwTak:Kato}.

We say that a positive continuous additive functional (PCAF in abbreviation) in the strict sense $A^{\nu}$ of ${\bf X}$ and a positive measure $\nu \in S_1({\bf X})$ are in the Revuz correspondence if they satisfy for any bounded $f \in \mathscr{B}_{+}(E)$,  
\begin{equation*}
\int_{E}f(x)\nu({\rm d}x) = \uparrow \lim_{t \downarrow 0}\frac{1}{t}\int_{E}\EE_{x}\left[\int_0^tf(X_s){\rm d}A_s^{\nu}\right]\m({\rm d}x). 
\label{Revuz corresp}
\end{equation*}
It is known that the family of equivalence classes of the set of PCAFs in the strict sense and the family of positive measures belonging to $S_1({\bf X})$ are in one to one correspondence under the Revuz correspondence (\cite[Theorem 5.1.4]{FOT}). 

A function $f$ on $E$ is said to be \emph{locally in $\F$ in the broad sense} (denoted as $f \in \dot{\F}_{\loc}$) if there is an increasing sequence of finely open Borel sets $\{E_n\}$ with $\bigcup_{n=1}^{\infty}E_n = E$ q.e. and for every $n \ge 1$, there is $f_{n} \in \F$ such that $f=f_{n}$ $\m$-a.e. on $E_n$. 
A function $f$ on $E$ is said to be \emph{locally in $\F$ in the ordinary sense} (denoted as $f \in {\F}_{\loc}$) if for any relatively compact open set $G$, there exists an element $f_G\in\F$ such that $f=f_G$ $\m$-a.e. on $G$. Clearly $\F_{\loc}\subset \dot\F_{\loc}$. It is shown in \cite[Theorem~4.1]{Kw:func}, $\F_e\subset \dot{\F}_{\loc}$.

For a signed measure $\nu=\nu_1-\nu_2$ with $\nu_i\in S_{D}^1({\bf X})$ (resp.~$\nu_i\in S_{D_0}^1({\bf X})$) ~$(i=1,2)$, $R_{\alpha}\nu$ (resp.~$R\nu$) is a difference of bounded $\alpha$-excessive (resp.~excessive) functions belonging to $\F_{\loc}$. 
Let $N_t^{R_{\alpha}\nu}$ be the CAF locally of zero energy appeared in 
the generalized Fukushima decomposition (see \cite[Theorem~6.1]{KimKuwaeTawara}). 
The following lemma is a variant of \cite[Lemma~5.4.1]{FOT}.
See \cite[Exercise~ 2.2.4]{FOT}.
\begin{lem}\label{lem:bddvar}
For a signed measure $\nu=\nu_1-\nu_2$ with $\nu_i\in S_{D}^1({\bf X})$ $(i=1,2)$, we have 
\begin{align}
N_t^{R_{\alpha}\nu}=\alpha\int_0^tR_{\alpha}\nu(X_s){\rm d}s-A_t^{\nu}, \quad t\in[0,\zeta[\label{eq:bddvarpositive}
\end{align}
$\PP_x$-a.s. for q.e.~$x\in E$, 
and in the transient case with $\nu_i\in S_{D_0}^1({\bf X})$ $(i=1,2)$, we have 
\begin{align}
N_t^{R\nu}=-A_t^{\nu}, \quad t\in[0,\zeta[\label{eq:bddvarzero}
\end{align}
$\PP_x$-a.s. for q.e.~$x\in E$. 
\end{lem}
 Next lemma is needed for the proof of main theorems.
 \begin{lem}
 \label{lem:Greenbdd}
Suppose that {\bf X} is transient and $\nu\in S({\bf X})$. 
Then there exists an $\E$-nest $\{F_n\}$ of compact sets 
 such that $\1_{F_n}\nu \in S_{00}^{(0)}({\bf X})$ for each $n\in\mathbb{N}$. 
 \end{lem}
 \begin{pf}
 The proof follows from \cite[Exercise~2.2.4]{FOT}. We omit the details. 
\end{pf}
 
 Now, let us introduce the notions of natural Green-tight measures of (extended) Kato class in the strict sense (\cite{KK:AnalChara}). First, we explain 
 the notion of weighted capacity of the Dirichlet form associated with the time changed process: 
 
 Let $\nu \in S_1({\bf X})$ and denote by $A_t^{\nu}$ the  
PCAF in the strict sense associated to $\nu$ in Revuz correspondence. Denote by 
${\sf S}_{o}^{\nu}$ the support of $A^{\nu}$ defined by 
${\sf S}_o^{\nu}:=\{x\in E\mid \PP_x(R=0)=1\}$, where $R(\omega):=\inf\{t>0\mid A_t^{\nu}(\omega)>0\}$.  ${\sf S}_{o}^{\nu}$ is nothing but the fine support of $\nu$, i.e., the topological support of $\nu$ with respect to the fine topology of {\bf X}. 
Let $(\check{\bf X}, \nu)$ be the time changed process of {\bf X}
 by $A_t^{\nu}$ and $(\check{\E},\check{\F})$ the associated Dirichlet form on 
$L^2({\sf S}^{\nu};\nu)$, where ${\sf S}^{\nu}$ is the support of $\nu$. It is known that 
$(\check{\E},\check{\F})$ is a regular Dirichlet form having $\C|_{{\sf S}^{\nu}}$ as its core and ${\sf S}^{\nu}\setminus {\sf S}_o^{\nu}$ is $\check{\E}$-polar, i.e., $1$-capacity $0$ set with respect to $(\check{\E},\check{\F})$. 
The life time of $(\check{\bf X}, \nu)$ is given by $A_{\zeta}^{\nu}$. 
Let $\check{\E}_1(f,f):=\check{\E}(f,f)+\int_Ef^2{\rm d}\nu$, and 
$\check{R}_1\varphi(x):=\EE_x[\int_0^{\infty}e^{-A_t^{\nu}}\varphi(X_t){\rm d}A_t^{\nu}]$,
the $1$-order resolvent of a $\nu$-a.e.~strictly positive bounded function $\varphi\in L^1(E;\nu)$ under $(\check{\bf X}, \nu)$.
Let $C^{\nu}:2^E\to[\,0,\,+\infty\,]$ be the 
weighted $1$-capacity with respect to $(\check{\E},\check{\F})$, i.e., for an open subset $G$ of $E$, 
we define 
$$
C^{\nu}(G):=\inf\left\{\check{\E}_1(f,f)\mid f\in\check{\F},\; f\geq \check{R}_1\varphi\quad\nu\text{-a.e.~on }G\right\}
$$
and for arbitrary subset $A$ of $E$
$$
C^{\nu}(A):=\inf\left\{C^{\nu}(G)\mid A\subset G,\;G\text{ is  an open subset of }E\right\}.
$$ We emphasize that $C^{\nu}$ is defined to be an outer capacity on $E$. By definition, $C^{\nu}(E\setminus {\sf S}^{\nu})=0$.  
Note that $C^{\nu}(E)\leq \check{\E}_1(\check{R}_1\varphi,\check{R}_1\varphi)=
\int_E\varphi(x) \check{R}_1\varphi(x)\nu({\rm d}x)
<\infty$ always holds. Note also that $C^{\nu}$ is tight in the sense that 
there exists an increasing sequence $\{K_n\}$ of compact subsets of ${\sf S}^{\nu}$ such that 
$\lim_{n\to\infty}C^{\nu}({\sf S}^{\nu}\setminus K_n)=0$ equivalently $\lim_{n\to\infty}C^{\nu}(E\setminus K_n)=0$. 

 \begin{df}\label{df:newKatoclass}  
{\rm Let 
$\nu\in S_1({\bf X})$.
\begin{enumerate}
\item\label{item:newKatoclass1}   $\nu$ is said to be a \emph{natural Green-tight measure of Kato class with respect to ${\bf X}$} if $\nu\in S_{D_0}^1({\bf X})$ and for any $\varepsilon>0$ 
there exist a closed subset $K=K(\varepsilon)$ of $E$ and a constant $\delta >0$ such that for all Borel set $B \subset K$ with $C^{\nu}(B) < \delta$, 
$$
\sup_{x\in E}\EE_x\left[A_{\tau_{B\cup K^c}}^{\nu}
\right]< \varepsilon.  
$$
Denote by $S_{N\!K_{\infty}}^1({\bf X})$ the family of natural Green-tight measures of Kato class with respect to ${\bf X}$.
\item\label{item:newKatoclass2}    $\nu$ is said to be a \emph{natural semi-Green-tight measure of extended Kato class with respect to ${\bf X}$} if $\nu\in S_{D_0}^1({\bf X})$ 
and
there exist a closed subset $K$ of $E$ and a constant $\delta >0$ such that for all Borel set $B \subset K$ with $C^{\nu}(B) < \delta$, 
$$
\sup_{x\in E}\EE_x\left[A_{\tau_{B\cup{K^c}}}^{\nu}
\right]< 1.  
$$
Denote by  $S_{N\!K_{1}}^1({\bf X})$  the family of  natural semi-Green-tight measures of extended Kato class with respect to ${\bf X}$.
\end{enumerate}  
}
\end{df}
\noindent
 We note that the closed set $K$ appeared in Definition~\ref{df:newKatoclass} can be taken to be compact, because the weighted $1$-capacity $C^{\nu}$ is tight. 

\begin{rem}\label{AdvantageofNV} 
{\rm \quad
\begin{enumerate}
\item\label{advanitem1} Let us denote by $S_{C\!K_{\infty}}^1({\bf X})$ (resp. $S_{C\!K_1}^1({\bf X})$) the family of Green-tight Kato class measures (resp. the family of semi-Green-tight Kato class measures) in the sense of Chen with respect to ${\bf X}$ (see~\cite{Chen:gaugeability2002, KK:AnalChara} for precise definitions). 
It is proved in \cite[Lemma~4.4]{KK:AnalChara} that $S_{C\!K_1}^1({\bf X})\subset S_{N\!K_1}^1({\bf X})\subset S_{E\!K}^1({\bf X})\cap S_{D_0}^1({\bf X})$ and $S_{C\!K_{\infty}}^1({\bf X})\subset S_{N\!K_{\infty}}^1({\bf X})\subset S_{K}^1({\bf X})\cap S_{D_0}^1({\bf X})$. 

\item\label{advanitem3} It is proved in \cite[Proposition 4.1]{KimKurKuwae:Gauge} that the class $S_{N\!K_{\infty}}^1({\bf X})$ is not so wide in the following sense: $S_{C\!K_{\infty}}^1({\bf X})=S_{N\!K_{\infty}}^1({\bf X})$ whenever ${\bf X}$ has the doubly Feller property of resolvent.
 
\item\label{advanitem2}
The advantage of the natural versions of semi-Green-tight measure is that 
they are stable under some Girsanov transform (cf.~\cite[Corollary~5.1 and Corollary~5.2]{KK:AnalChara}).
\end{enumerate}
}
\end{rem} 
If {\bf X} is transient, we can consider the Green kernel $R(x,y)$ for $x,y\in E$ under ${\bf (AC)}$. We set ${\sf d}:=\{(x,y)\mid R(x,y)=0\text{ or }+\infty\}$ and $E^z:=\{x\in E\mid (x,z)\in E\times E\setminus{\sf d}\}$. Since $R(x,y)>0$ for all $x,y\in E$ (see \cite[Lemma~6.1]{KimKurKuwae:Gauge}), we see 
${\sf d}:=\{(x,y)\mid R(x,y)=+\infty\}$. By Getoor~\cite{Get:TranRec}, there exists 
 $g\in L^1(E;\m)$ with $0<g\leq 1$ $\m$-a.e. such that $Rg\in\mathscr{B}_b(E)$, we see $\m(E\setminus E^z)=0$ for all $z\in E$. More strongly we have the following:
\begin{pr}\label{pr:DiagonalVanish}
Suppose that {\bf X} is transient. Then
$E\setminus E^z$ is exceptional for all $z\in E$. 
\end{pr} 
\begin{pf}
For any $\nu\in S_{00}^{(0)}({\bf X})$, $R\nu\in\mathscr{B}_b(E)$. 
This implies that for given $z\in E$, $\nu(E\setminus E^z)=0$ for all 
$\nu\in S_{00}^{(0)}({\bf X})$. By Lemma~\ref{lem:Greenbdd}, 
for given $z\in E$, $\nu(E\setminus E^z)=0$ for any smooth measure $\nu$. Therefore we obtain the assertion. 
\end{pf}

\section{Analytic characterization of gaugeability and subcriticality}\label{sec:AofGandSubcriticality}
Recall that  
$F=F_1-F_2$ where $F_1$ and $F_2$ are non-negative bounded functions on $E\times E_\partial$ symmetric on $E\times E$
which is extended to a function defined on $E_\partial \times E_\partial$ vanishing 
on the diagonal set ${\sf diag}$ of $E_\partial \times E_\partial$. Let
$$
A_{t}^{F}=A^{F_1}_t-A^{F_2}_t,\quad A_t^{F_i}:=\sum_{0<s\le t}F_i(X_{s-},X_{s})\quad  (i=1,2).
$$
 Note that if $N(F_1+F_2)\mu_H \in S_1({\bf X})$, then $N(|F|)\mu_H \in S_1({\bf X})$. In this case, $A^F$ can be defined as an additive functional in the strict sense.  
Note that $A_t^F=\sum_{0<s\le t}\1_{\{s<\zeta\}}F(X_{s-},X_{s})$ provided $F(x,\partial)=0$ for $x\in E_{\partial}$. Hereafter, we always assume $F(x,\partial)=0$, $x\in E_{\partial}$. 
For a bounded finely continuous (nearly) Borel function 
$u\in{\F}_{\loc}\cap QC(E_{\partial})$ satisfying $\mu_{\<u\>}\in S_1({\bf X})$ and $N(|F|)\mu_H\in S_1({\bf X})$, we set     
$$
F^{u}(x,y):=F(x,y)+u(x)-u(y)
\quad \text{
and } \quad G^{u}=e^{F^{u}}-1$$
with identifying $F^0=F$ and $G^0=G:=e^{F}-1$. 

Let $M^{F}_{t}=\sum_{0<s \le t}F(X_{s-},X_{s}) - \int_{0}^{t}N(F)(X_{s}){\rm d}H_{s}$.
By \cite[(3.1)--(3.4)]{KK:AnalChara},
 there exist  purely discontinuous locally square integrable local martingale additive functionals $M^{F^u}$ and $M^{G^u}$
 on $[\![0,\zeta[\![$ such that 
 $\Delta M^{F^u}_{t}=F^{u}(X_{t-},X_{t})$ and 
 $\Delta M^{G^u}_{t}=G^{u}(X_{t-},X_{t})$, 
 $t\in[0,\zeta[$ 
 $\PP_x$-a.s.~for all $x\in E$ 
$M^{F^u}$ and $M^{G^u}$ are given by 
\begin{equation}
M^{F^u}_{t}=M_{t}^{F} + M_{t}^{-u,j} + M_{t}^{-u,\kappa}, \quad t <\zeta
\label{MF1}
\end{equation}
and \begin{equation}
M_{t}^{G^u}=M_{t}^{F^u}+\sum_{0<s\le t}(G^{u}-F^{u})(X_{s-},X_{s}) - \int_{0}^{t}N(G^{u}-F^{u})(X_{s}){\rm d}H_{s}, \quad t <\zeta. 
\label{MG1}
\end{equation}

Let $Y_{t}:={\rm Exp} (M^{G^u}+M^{-u,c})_t$ be the Dol\'eans-Dade exponential of $M^{G^u}_t+M^{-u,c}_t$, that is, $Y_{t}$ is the unique solution of  
\begin{equation}
Y_{t}=1+\int_{0}^{t}Y_{s-}{\rm d}(M_{s}^{G^u}+M^{-u,c}_s), \quad t < \zeta,~~{\PP}_{x}\text{-a.s.} 
\label{SDE}
\end{equation}
It then follows from Dol\'eans-Dade formula and (\ref{MG1}) that 
\begin{align}
Y_{t}=\exp \left(M_{t}^{F^u}+M_t^{-u,c}-\int_{0}^{t}N(G^{u}-F^{u})(X_{s}){\rm d}H_{s}-\frac{\,1\,}{2}\<M^{u,c}\>_t\right)
\label{eqY}
\end{align}
(\cite[Theorem 3.1]{KK:AnalChara}). 
Note that $Y_{t}$ is a 
positive and local martingale, therefore supermartingale on $[\![0,\zeta[\![$. 
It is inconvenient to treat additive functionals on $[\![0,\zeta[\![$ for our purpose. 
We see that $Y_t$ 
can be extended for all $t\in[\,0,\,+\infty\,[$ provided $\mu_{\<u\>} + N(F_1+F_2)\mu_H\in S_D^1({\bf X})$ (\cite[Proposition 3.1]{KK:AnalChara}). 
\vskip 0.1cm

Let us denote by ${\bf Y}=(X_{t}, {\PP}_{x}^{Y})$ 
the transformed process of ${\bf X}$ by $Y_{t}$. 
The transition semigroup $({P}_{t}^{Y})_{t\ge 0}$ of ${\bf Y}$ 
is defined by 
\begin{align*}
{P}_{t}^{Y}\!f(x):={\EE}_{x}[Y_{t}f(X_t)]. 
\end{align*}

\begin{thm}[{\cite[Theorem~3.2]{KK:AnalChara}}]\label{FeynmanKacMod3} 
Assume that a bounded function 
$u\in{\F}_{\loc}\cap QC(E_{\partial})$ admits a Fukushima's decomposition holding up to infinity under ${\PP}_x$ for q.e.~$x\in E$ 
and $N(F_1+F_2)\mu_H \in S_1({\bf X})$. Let $({\E}^{Y},{\F}^{Y})$ be the Dirichlet form of ${\bf Y}$ on $L^{2}(E;e^{-2u}\m)$. Then $\F={\F}^{Y}$ and for any $f\in \F^{Y}$
\begin{equation*}
\begin{split}
{\E}^{Y}(f,f)&=\frac{\,1\,}{2}\int_{E}e^{-2u(x)}\mu_{\< f\>}^{c}({\rm d}x) +\int_{E\times E}(f(x)-f(y))^{2}e^{F(x,y)-u(x)-u(y)}J({\rm d}x{\rm d}y). 
\end{split}
\end{equation*} 
\end{thm}

By Theorem~\ref{FeynmanKacMod3} combined with the boundedness of $u$ and $F$, there exists a constant $C_E >0$ such that 
$$
C_E^{-1}\E^{(c)}(f,f) \le \E^{Y,(c)}(f,f) \le C_E\E^{(c)}(f,f)
$$
for $f \in \F^Y$ and 
$$
C_E^{-1}J({\rm d}x{\rm d}y)\leq J^Y({\rm d}x{\rm d}y)\leq C_EJ({\rm d}x{\rm d}y).
$$ 
Here $J^Y({\rm d}x{\rm d}y):=e^{F(x,y)-u(x)-u(y)}J({\rm d}x{\rm d}y)$ is the jumping measure of $(\E^Y,\F^Y)$. 
Therefore we see under {\bf (A.1)}
(resp.~{\bf (A.2)}) 
that $({\E}^{Y},{\F}^{Y})$ admits a heat kernel $p_t^Y(x,y)$ (resp.~heat kernel $p_t^Y(x,y)$ in the strict sense)  on $]\,0,\,+\infty\,[\times E \times E$ such that 
\begin{align*}
p_t^Y(x,y)&\precsim \phi_2(t,x,y)\quad \m\text{-a.e.}~~ d(x,y)\in \mathbb{T}_t\quad\text{(resp.~$\phi_1(t,x,y)\precsim p_t^Y(x,y)\precsim \phi_2(t,x,y)$,\quad $d(x,y)\in\mathbb{T}_t$)}
\end{align*}
for $u$ and $F$ appeared in Theorem~\ref{FeynmanKacMod3}.

Recall that $(\Q,\F_b)$ is the quadratic form on $L^2(E;\m)$ defined in 
\eqref{eq:quadraticform}. 
Note that for bounded finely continuous (nearly) Borel function $u\in{\F}_{\loc}\cap QC(E_{\partial})$ and $f \in \F\cap C_0(E)$, we see that $fe^u=f(e^u-1)+f=f(e^{u_O}-1)+f \in \F$, where 
$O$ is a relatively compact open set with ${\rm supp} [f] \subset O$ and 
$u_O \in \F_b$ satisfying $u=u_O$ $\m$-a.e. on $O$. 
Moreover, it follows from \eqref{eqY}, Theorem~\ref{FeynmanKacMod3} and the Feynman-Kac formula that for $f \in \F\cap C_0(E)$,   
$$
\E^Y(f,f)=\Q(fe^{-u},fe^{-u}) + \int_{E}f^2e^{-2u}{\rm d}\overline{\mu}.
$$

Consider the non-local Feynman-Kac transforms by the additive functionals $A_t:=N_t^u+A_t^{\mu}+A_t^F$ of the form 
\begin{equation} 
e_A(t):=\exp (A_t), \quad t \ge 0. 
\label{genFKtansform}
\end{equation}
We see that for $\mu_{\<u\>} + N(F_1+F_2)\mu_H\in S_D^1({\bf X})$  
\begin{align}
e_A(t)
=e^{u({X}_t)-u(X_0)}Y_t\exp \left(A_t^{\overline{\mu}}\right)
\quad t\in[\,0,\,+\infty\,[,
\label{relationeAandY}
\end{align}
which implies that for $x\in E$ and $f\in\mathscr{B}_+(E)$,
\begin{align}
P_t^A\!f(x)&=e^{-u(x)}\EE_x^{Y}\left[\exp \left(A_t^{\overline{\mu}}\right)(e^uf)(X_t)\right].
\label{eq:relationsemi2}
\end{align}
We note that if $\nu\in S_{D}^1({\bf X})$ (resp.~$\nu\in S_{D_0}^1({\bf X})$), then 
Stollmann-Voigt's inequality (\cite[Theorem~3.1]{SV:potential}) tells us 
$\int_{E}f^2{\rm d}\nu\leq\|R_1\nu\|_{\infty}\E_1(f,f)$ for $f\in\F$ (resp.~$\int_{E}f^2{\rm d}\nu\leq\|R\nu\|_{\infty}\E(f,f)$ for $f\in\F_e$), hence 
$\F\subset 
L^2(E;\nu)$ (resp.~$\F_e\subset 
L^2(E;\nu)$). 
It is shown in \cite{KimKurKuwae:Gauge} that under ${\mu}_1+N(e^{F_1}-1)\mu_H\in S_{E\!K}^1({\bf X})$, $\mu_{\<u\>}\in S_K^1({\bf X})$ and 
$\mu_2+N(F_2)\mu_H\in S_1({\bf X})$ and $\alpha\geq0$, the symmetric 
$\alpha$-order 
resolvent kernel 
 $R^{A}_{\alpha}(x,y)$ of the Feynman-Kac semigroup $P_t^{A}$ can be defined for all 
 $x,y\in E$ (possibly infinite), it is $\alpha$-excessive with respect to $(P_t^{A})_{t \ge 0}$ in $x$ (and in $y$), and it satisfies a resolvent equation as in \cite[(4.2.12)]{FOT}. 
 Moreover, for each $x\in E$, $y\mapsto R_{\alpha}^A(x,y)$ is finely continuous (see \cite[Lemma~6.1]{KK:AnalChara}). We write $R^A(x,y):=R_0^A(x,y)$ for $x,y\in E$.

We consider the following conditions: 
\begin{enumerate}
\item[{\bf (A)}] $\mu_1+N(e^{F_1}-1)\mu_H\in S_{E\!K}^1({\bf X})$, $\mu_{\<u\>}\in S_K^1({\bf X})$ and $\mu_2+N(F_2)\mu_H\in S_D^1({\bf X})$. 
\item[${\bf (A)}^*$\!\!] $\mu_1+N(e^{F_1}-1)\mu_H\in S_{E\!K}^1({\bf X}^*)$, $\mu_{\<u\>}\in S_K^1({\bf X}^*)$.
\end{enumerate}
Here ${\bf X}^*$ is the subprocess killed by $e^{-A_t^{\mu_2}-A_t^{F_2}}$. It is easy to see that ${\bf (A)}$ implies ${\bf (A)^*}$. 
Let us recall the quadratic form $(\Q,\F)$ on $L^2(E;\m)$ defined in \eqref{eq:quadraticform}.  

\begin{lem}[{cf.~\cite[Lemma~2.1]{KimKuwae:general}}]
\label{lem:strongcontinuity}
Suppose ${\bf (A)}^*$. 
Then for any sufficiently large $\alpha>0$ 
there exists $C>0$ such that 
\begin{align}
C^{-1}\E_1(f,f)\leq\Q_{\alpha}(f,f)\leq C\E_1(f,f),\quad \text{ for }\quad f\in\F.
\label{eq:Eqiuvalent}
\end{align}
Consequently, 
$(P^A_t)_{t\geq0}$
 defines a strongly continuous semigroup on 
$L^2(E;\m)$ associated with $(\Q,\F)$ on $L^2(E;\m)$. 
\end{lem}

\vskip 0.1cm

Next theorem extends \cite[Theorem~1.1]{KimKurKuwae:Gauge} 
(see also \cite[Theorem~1.1]{KK:AnalChara}). 
\begin{thm}
\label{thm:analgauge1} 
Suppose that ${\bf X}$ is transient. Let $u\in {\F}_{\loc}\cap QC(E_{\partial})$ be a bounded finely continuous (nearly) Borel function on $E$. 
Assume $\mu_1+N(e^{F_1}-1)\mu_H\in S_{N\!K_1}^1({\bf X})$, $\mu_{\<u\>}\in S_{N\!K_{\infty}}^1({\bf X})$ and   
$\mu_2+N(F_2)\mu_H\in S_{D_0}^1({\bf X})$. 
Then the following are equivalent:
\begin{enumerate}
\item\label{item:gaugeTh} The functional \eqref{genFKtansform} is gaugeable, that is, $\sup_{x \in E}\EE_{x}[e_{A}(\zeta)] <\infty$.
\item\label{item:analyticTh} $\lambda^{\Q}(\overline{\mu}_1) > 0$.
\item\label{item:subcriticalTh} For each $x\in E$, $R^{A}(x,y) < \infty$ for $\m$-a.e.~$y\in E$. 
\item\label{item:subcriticalTh*} 
$R^{A}(x,y) < \infty$ for $\m$-a.e.~$x,y\in E$.
\end{enumerate}
\end{thm}
\begin{pf}
As noted in Remark~\ref{rem:mainresult1}(1), 
the condition $\mu_1+N(e^{F_1}-1)\mu_H\in S_{N\!K_1}^1({\bf X})$ is equivalent to $\mu_1+N(e^{[u]}(e^{F_1}-1))\mu_H\in S_{N\!K_1}^1({\bf X})$
under $\mu_{\<u\>}\in S_{N\!K_{\infty}}^1({\bf X})$. 
The equivalence 
\eqref{item:gaugeTh}$\Longleftrightarrow$\eqref{item:analyticTh} is proved in \cite[Theorem~1.1]{KimKurKuwae:Gauge} and 
the implication \eqref{item:subcriticalTh} $\Longrightarrow$\eqref{item:subcriticalTh*} is easy. So it suffices to prove the implications 
\eqref{item:subcriticalTh*}$\Longrightarrow$\eqref{item:gaugeTh}$\Longrightarrow$\eqref{item:subcriticalTh}. 

Recall that $[u](x,y):=u(x)-u(y)$. Let ${\bf U}$ be the transfomred process of ${\bf X}$  
by the supermartingale multiplicative 
functional $U_t:=\Exp(M^{e^{[u]}-1}+M^{-u,c})_t$. Let $\overline{\nu}:=\overline{\nu}_1-\overline{\nu}_2$ be the measure given by $\overline{\nu}_1=\mu_1+N(e^{[u]}-[u]-1)\mu_H+\frac{\,1\,}{2}\mu_{\<u\>}^c$ and $\overline{\nu}_2=\mu_2$. 
Then we see by \cite[Corollary~5.1(1),(5)]{KimKurKuwae:Gauge} that $e^{-2u}(\overline{\nu}_1+N(e^{[u]}(e^{F_1}-1))\mu_H)\in S_{N\!K_1}^1({\bf U})$ and $e^{-2u}(\overline{\nu}_2+N(e^{[u]}F_2)\mu_H)\in S_{D_0}^1({\bf U})$. From these facts with \cite[Lemma 4.4 and (4.10)]{KimKurKuwae:Gauge}, one can easily show that \eqref{item:gaugeTh} is equivalent to 
$\sup_{x\in E}\EE_x^U[\exp(A_{\zeta}^{\overline{\nu}}+A_{\zeta}^F)]<\infty$ and $R^A(x,y)$ can be written by $R^A(x,y)=e^{-u(x)-u(y)}(R^U)^{e^{-2u}\overline{\nu},F}(x,y)$
where 
$R^{U}$ is 
$0$-order 
resolvent kernel of  ${\bf U}$.
In view of these, we can and do assume $u=0$. 

Since {\bf X} is transient, by Getoor~\cite{Get:TranRec}, 
there exists $g\in L^1(E;\m)$ with $0<g\leq 1$ such that 
$Rg\in\mathscr{B}_b(E)$. 
In particular, $g\m\in S_{D_0}^1({\bf X})$. 
 Let $(\check{\bf X}, g\m)$ be the time changed process by $A_t^{g\m}:=\int_0^tg(X_s){\rm d}s$, i.e., 
$(\check{\bf X}, g\m):=(\Omega,X_{\tau_t^{g\m}},\PP_x)_{x\in E}$, where $\tau_t^{g\m}:=\inf\{s>0\mid A_s^{g\m}>t\}$ is the right continuous inverse of $A_t^{g\m}$.  
Suppose that \eqref{item:gaugeTh} holds. Then by \cite[Lemma~5.1(2)]{KimKurKuwae:Gauge} 
\begin{align}
\sup_{x\in E}\EE_x[e_{A_{\tau^{g\m}}}(A^{g\m}_{\zeta})]=\sup_{x\in E}\EE_x[e_A(\zeta)]<\infty.\label{eq:gaugetimechange}
\end{align} 
Here $A_{\tau^{g\m}}$ is the 
CAF under $(\check{\bf X}, g\m)$ and 
$A^{g\m}_{\zeta}$ is the life time of $(\check{\bf X}, g\m)$. 
Applying \cite[Lemma~4.3]{KimKurKuwae:Gauge} to the time changed process $(\check{\bf X},g\m)$  with $g\m\in S_{D_0}^1(\check{\bf X},g\m)$, 
\eqref{eq:gaugetimechange} is equivalent to 
\begin{align*}
\sup_{x\in E}R^Ag(x)=\sup_{x\in E}\EE_x\left[\int_0^{\infty}e_A(t)g(X_t){\rm d}t\right]
=\sup_{x\in E}\EE_x\left[\int_0^{\infty}e_A(\tau_t^{g\m}){\rm d}t\right]
<\infty.
\end{align*} 
Then we have \eqref{item:gaugeTh}$\Longrightarrow$\eqref{item:subcriticalTh}. 
Next suppose that \eqref{item:subcriticalTh*} holds. 
Then there exists $N\in\mathscr{B}(E)$ such that 
for $x\in E\setminus N$, $R^A(x,y)<\infty$ $\m$-a.e.~$y\in E$. 
For each fixed $x\in E\setminus N$, we set 
\begin{align*}
g_x(y):=\sum_{n=1}^{\infty}\frac{1}{2^nn\m(G_n)}\1_{E_n^x\setminus E_{n-1}^x}(y),
\end{align*}
where $\{G_n\}$ is an increasing sequence of relatively compact  
open sets such that $E=\bigcup_{n=1}^{\infty}G_n$ and $\m(G_n)\geq1$, and 
$E_n^x:=\{y\in G_n\mid R^A(x,y)\leq n\}$ for $n\in\mathbb{N}$ with $E_0^x:=\emptyset$. It is easy to see that $0< g_x\leq 1$ $\m$-a.e. and $g_x\in \mathscr{B}_b(E)$, because $R^A(x,y)<\infty$ $\m$-a.e.~$y\in E$ yields $\m(\bigcap_{n=1}^{\infty}(E\setminus E_n^x))=0$.
Then for any $z\in E$
\begin{align}
R^Ag_x(z)\leq\sum_{n=1}^{\infty}\frac{1}{2^n}=1.\label{eq:gaugeGx}
\end{align}
Take any $g\in \mathscr{B}(E)$ with $0<g\leq1$ on $E$. Then $(g_x\land g)\m\in S_{D_0}^1({\bf X})$. 
Applying \cite[Lemma~4.3{{(1)$\Longleftrightarrow$(4)}}]{KimKurKuwae:Gauge} to the time changed process $(\check{\bf X}, (g_x\land g)\m)$ with 
$(g_x\land g)\m\in S_{D_0}^1(\check{\bf X}, (g_x\land g)\m)$, this  
is equivalent to the gaugeability of AF $A_{\tau^{(g_x\land g)\m}}$ under $(\check{\bf X}, (g_x\land g)\m)$ (see \cite[Lemma~5.1(2)]{KimKurKuwae:Gauge}):
\begin{align*}
\sup_{y\in E}\EE_y[e_A(\zeta)]=\sup_{y\in E}\EE_y[e_{A_{\tau^{(g_x\land g)\m}}}(A^{(g_x\land g)\m}_{\zeta})]<\infty.
\end{align*} 
Therefore we obtain \eqref{item:gaugeTh}.  
\end{pf}

\section{Proofs of Theorem~\ref{thm:mainresult1}, Corollary~\ref{cor:mainresult} and Theorem~\ref{thm:shortstability}}\label{sec:Stability results}

Recall that 
$\overline{\mu}:=\overline{\mu}_1-\overline{\mu}_2$ and 
$\overline{\mu}^*:=\overline{\mu}_1^*-\overline{\mu}_2^*$ are the signed measures 
defined in  \eqref{e:mu12} and  \eqref{e:mu12*} respectively.
In this section, we will prove our main results stated in Section \ref{sec:Result}. Before proving them, we give several lemmas on the measure $
\overline{\mu}
$ and its gauge function $h$ with respect to the process ${\bf Y}$. Throughout this section, 
$u$ is a bounded finely continuous (nearly) Borel function in ${\F}_{\loc}\cap QC(E_{\partial})$.
 
\begin{lem}\label{lem;barmusemiGreen}
Assume that {\bf X} is transient. Suppose that 
$\mu_1+N(e^{F_1}-1)\mu_H\in S_{N\!K_1}^1({\bf X})$, $\mu_{\<u\>}\in S_{N\!K_{\infty}}^1({\bf X})$ and $\mu_2+N(F_2)\mu_H\in S_{D_0}^1({\bf X})$. Then, the following are equivalent:
\begin{enumerate}
\item\label{item:gaugeTh1} The functional $\exp (A_t^{\overline{\mu}})$ is gaugeable under ${\bf Y}$, that is, $\sup_{x \in E}\EE_{x}^Y[\exp (A_{\zeta}^{\overline{\mu}})] <\infty$.
\item\label{item:analyticTh1} $\lambda^{\Q}(\overline{\mu}_1) > 0$.
\end{enumerate}
\end{lem} 
\begin{pf}
By assumption, we see $\mu_{\<u\>}+N(F_1+F_2)\mu_H\in S_{D_0}^1({\bf X})$. Then we have from \cite[Lemma~4.9(2)]{KK:AnalChara} that 
\begin{align}
\EE_x^Y\left[
\exp(A_{\zeta}^{\overline{\mu}})\right]=e^{u(x)}\EE_x[e^{-u(X_{\zeta-})}e_A(\zeta)].\label{eq:equality}
\end{align}
 Hence $\sup_{x\in E}\EE_x[e_A(\zeta)]<\infty$ is equivalent to 
$\sup_{x\in E}\EE_x^Y[\exp(A_{\zeta}^{\overline{\mu}})]<\infty$. 
Therefore the assertion follows from Theorem~\ref{thm:analgauge1}. 
\end{pf}
\vskip 0.2cm
We define the gauge function $h$ with respect to ${\bf Y}$ by 
$$
h(x):=\EE_x^{Y}\left[
\exp(A_{\zeta}^{\overline{\mu}})\right]
,\quad x\in E_{\partial}.
$$
Since any AF $C_t$ of ${\bf Y}$ satisfies $\PP_{\partial}^{Y}(C_t\equiv 0)=1$, we see $h(\partial)=1$.
Then, by \eqref{eq:equality} 
\begin{align}
0< e^{-\|R^Y \overline{\mu}_2
\|_{\infty}}=e^{-\sup_{y\in E}\EE^Y_y[A_{\zeta}^{\mu_2}+A_{\zeta}^{F_2}]}
\leq h(x)\leq\sup_{y\in E_{\partial}}\EE_y^{Y}\left[e^{A_{\zeta}^{\overline{\mu}}}\right]<\infty
\label{eq:gaugeable}
\end{align}
whenever $u=0$ and $\lambda^{\Q}(\overline{\mu}_1)>0$. 
We have the following lemma in a similar way of \cite[Lemma~2]{Tak:GaussianHK}:

\begin{lem}\label{htoR}
Assume that {\bf X} is transient and $u=0$. Suppose that 
$\overline{\mu}_1^*=\mu_1+N(e^{F_1}-1)\mu_H\in S_{N\!K_1}^1({\bf X})$ and $\overline{\mu}_2=\mu_2+N(F_2)\mu_H\in S_{D_0}^1({\bf X})$, and 
 $\lambda^{\Q}(\overline{\mu}_1)>0$. Then we have $\overline{\mu}_2\in S_{D_0}^1({\bf Y})$ and  
\begin{align}
h(x)=R^{Y}(h\overline{\mu})(x)+1,
\quad x\in E.\label{eq:identity}
\end{align}
In particular, we have 
$\overline{\mu}_1=\mu_1+N(e^{F_1-F_2}+F_2-1)\mu_H\in S_{D_0}^1({\bf Y})$, hence $\overline{\mu}_1^*\in S_{D_0}^1({\bf Y})$. 
\end{lem}
\begin{pf}
Since $\lambda^{\Q}(\overline{\mu}_1)>0$, we know $\sup_{x\in E_{\partial}}h(x)<\infty$ by Lemma~\ref{lem;barmusemiGreen}.
Recall that, since we assume that $u=0$, we have  $\overline{\mu}=\overline{\mu}_1-\overline{\mu}_2$ with $\overline{\mu}_1=\mu_1+N(e^{F_1-F_2}+F_2-1)\mu_H$ and $\overline{\mu}_2=\mu_2+N(F_2)\mu_H$. 
Thus,  we see $\overline{\mu}=\overline{\mu}_1^*-({\mu}_2+N(e^{F_1}(1-e^{-F_2})\mu_H)$.  
Thanks to the boundedness of $F_1$ and $F_2$,
$(\E^Y,\F^Y)$ is equivalent to $(\E,\F)$ so that
$\nu\in S_0^{(0)}({\bf X})$ is equivalent to $\nu\in S_0^{(0)}({\bf Y})$, 
consequently $\nu\in S({\bf X})$ is equivalent to $\nu\in S({\bf Y})$. Since 
$\overline{\mu}_1^*\in S({\bf Y})$, 
there exists an $\E^{Y}$-nest $\{K_n\}$ of compact sets such that $\1_{K_n}\overline{\mu}_1^* \in S_{00}^{(0)}({\bf Y})$ for each $n \in \mathbb{N}$ by Lemma~\ref{lem:Greenbdd}. 
Set $\overline{\mu}^n:=\1_{K_n}\overline{\mu}_1^*-({\mu}_2+N(e^{F_1}(1-e^{-F_2})\mu_H)$, 
$M_t^n:=\exp(A_t^{\overline{\mu}^n})$ and $h_n(x):=\EE_x^Y[e^{A_{\zeta}^{\overline{\mu}^n}}]$. 
Then  
\begin{align*}
\wh{h}(x):= \EE_x^Y\left[\exp\left(-A_{\zeta}^{\mu_2}-A_{\zeta}^{N(e^{F_1}(1-e^{-F_2}))\mu_H}\right)\right]\leq h_n(x).
\end{align*} 
In the same way of the proof of \cite[Lemma~4.9(2)]{KK:AnalChara}, we can deduce 
\begin{align*}
\hat{h}(x)&=\EE_x\left[\exp\left(A_{\zeta}^F-A_{\zeta}^{N(e^F-1)\mu_H}-A_{\zeta}^{\mu_2}-A_{\zeta}^{N(e^{F_1}(1-e^{-F_2}))\mu_H}\right)\right]\\
&\geq \EE_x\left[\exp\left(-A_{\zeta}^{\mu_2}-A_{\zeta}^{F_2}-A_{\zeta}^{N(e^{F_1}-1)\mu_H} \right)\right]\\
&\geq \exp\left(-\|R(\overline{\mu}_2+N(e^{F_1}-1)\mu_H)\|_{\infty} \right)>0. 
\end{align*}
Hence 
\begin{align}
0<\exp\left(-\|R(\overline{\mu}_2+N(e^{F_1}-1)\mu_H)\|_{\infty} \right)\leq h_n(x)\leq h(x)\leq \sup_{x\in E_{\partial}}h(x)<\infty.\label{eq:uniformLowBDD}
\end{align}
The process $(h_n(X_t)M_t^n)_{t\geq0}$ is a closed $\PP_x^{Y}$-martingale for all $x\in E$. Indeed, by the Markov property, we see
\begin{align*}
\EE_x^{Y}[M_{\infty}^n\;|\;\mathscr{F}_t]=\EE_x^{Y}[M_t^n\cdot M_{\infty}^n\circ\theta_t\;|\;\mathscr{F}_t]
=M_t^n\EE_{X_t}^{Y}[M_{\infty}^n]=M_t^n h_n(X_t).
\end{align*}
Thus for any $(\mathscr{F}_t)$-stopping time $T$, we have  
$$
\EE_x^{Y}\left[h_n(X_T)\1_{\{T<\infty\}}\right]=\EE_x^{Y}\left[M_{\infty}^n(M_{T}^n)^{-1}\1_{\{T<\infty\}}\right].
$$
Hence, applying \cite[(1.13) Exercise p.~186]{RevuzYor}, we have 
\begin{align*}
\EE_x^Y\left[\int_0^{\infty}h_n(X_s){\rm d}A_s^{\overline{\mu}^n}\right]&=\EE_x^Y\left[e^{A_{\zeta}^{\overline{\mu}^n}}\int_0^{\infty}e^{-A_s^{\overline{\mu}^n}}
{\rm d}A_s^{\overline{\mu}^n}\right]=\EE_x^Y\left[e^{A_{\zeta}^{\overline{\mu}^n}}\int_0^{\infty}{\rm d}\left(-e^{-A_s^{\overline{\mu}^n}}\right) \right] \\
&=\EE_x^Y\left[e^{A_{\zeta}^{\overline{\mu}^n}}\!\!\left(1-e^{-A_{\zeta}^{\overline{\mu}^n}}\right)\right]=h_n(x)-1.\label{eq:equality}
\end{align*}
From this and the inequality $e^{-\|F_2\|_{\infty}}F_2\leq e^{-F_2}(e^{F_2}-1)$, 
\begin{align*}
e^{-\|F_2\|_{\infty}}
R^Y(h_n\overline{\mu}_2)(x)&\leq 
R^Yh_n({\mu}_2+N(1-e^{-F_2})\mu_H)(x)\\
&\leq
R^Yh_n({\mu}_2+N(e^{F_1}(1-e^{-F_2})\mu_H)(x)\\
&= R^Y(h_n\1_{K_n}\overline{\mu}_1^*)(x)-h_n(x)+1.
\end{align*}
Thus, we see $h_n\overline{\mu}_2\in S_{D_0}^1({\bf Y})$. 
By \eqref{eq:uniformLowBDD}, 
we have $\overline{\mu}_2\in S_{D_0}^1({\bf Y})$. 
By replacing $M_t^n$ (resp.~$h_n$) with $M_t:=\exp(A_t^{\overline{\mu}})$ (resp.~ $h$), we can deduce 
 \eqref{eq:identity} under $\overline{\mu}_2\in S_{D_0}^1({\bf Y})$ by the same manner as 
shown above (see \eqref{eq:equality}). 

Finally, we prove the last assertion.  
From \eqref{eq:identity}, 
\begin{align*}
R^Y(h\overline{\mu}_1)=R^Y(h\overline{\mu}_2)+h-1
\end{align*} 
is bounded above. 
By \eqref{eq:uniformLowBDD}, $\overline{\mu}_1\in 
 S_{D_0}^1({\bf Y})$ and $N(F_2)\mu_H\leq\overline{\mu}_2\in S_{D_0}^1({\bf Y})$ implies
\begin{align*}
\overline{\mu}_1^*&=\overline{\mu}_1+N(e^{F_1}-e^{F_1-F_2}-F_2)\mu_H\\&\leq 
\overline{\mu}_1+N\left(e^{F_1}(1-e^{-F_2})
\right)\mu_H
\leq \overline{\mu}_1+e^{\|F_1\|_{\infty}}
N(F_2)\mu_H\in 
S_{D_0}^1({\bf Y}).
\end{align*}
\end{pf}
The expression \eqref{eq:identity} 
yields the fine continuity of $h$ on $E$ with respect to ${\bf Y}$ (when $u=0$).
 From this observation, $h$ is $\E^{Y}$-quasi continuous on $E$. 
Moreover, since $R^{Y}(h\overline{\mu})$ is a difference of bounded excessive functions on $E$ with respect to ${\bf Y}$, we have  
\begin{equation}
h-1=R^{Y}(h\overline{\mu})\in{\F}_{\rm loc}^{Y}
\label{h-1Floc}
\end{equation}
by use of \cite[Lemma 2.3.2]{FOT}. Indeed, 
let $\{G_n\}$ be an increasing sequence of relatively compact open 
sets and $e_{G_n}^{Y}$ its $1$-equilibrium potential 
with respect to $(\E^{Y},\F^{Y})$ on $L^2(E;e^{-2u}\m)$.  
Take any bounded excessive function $f$ with respect to  ${\bf Y}$ and set $f_n:=f\wedge \|f\|_{\infty}e_{G_n}^{Y}$. Then $f_n \le \|f\|_{\infty}e_{G_n}^{Y}$ with \cite[Lemma~2.3.2]{FOT} yields 
$f_n\in\F^{Y}$ and 
$f=f_n$ on $G_n$ for each $n\in\mathbb{N}$. 
Thus we have 
$f\in {\F}_{\loc}^{Y}$.

\begin{lem}\label{Fukudecompforh}
Assume that {\bf X} is transient and $u=0$. Suppose that 
$\overline{\mu}_1^*=\mu_1+N(e^{F_1}-1)\mu_H\in S_{N\!K_1}^1({\bf X})$  and $\overline{\mu}_2=\mu_2+N(F_2)\mu_H\in S_{D_0}^1({\bf X})$ 
hold.  
Suppose further that  
 $\lambda^{\Q}(\overline{\mu}_1)>0$. Then the additive functional $h(X_t)-h(X_0)$ admits the following decomposition in the strict sense:
\begin{align}
\left\{\begin{array}{ll}
h(X_t)-h(X_0)=M_t^{h}+N_t^{h}, \\ 
N_t^{h}=-\int_0^t h(X_s){\rm d}A_s^{\overline{\mu}}
\end{array}\right. 
\label{FukuDecomp}
\end{align}
for all $t\in[\,0,\,+\infty\,[$ $\PP_x^{Y}$-a.s.~for~all~$x\in E$, 
where $M_t^{h}$ (resp.~$N_t^{h}$) 
is a square integrable martingale additive functional in the strict sense 
(resp.~CAF in the strict sense which is locally of zero energy) under ${\bf Y}$. 
\end{lem}
\begin{pf}
By Lemma~\ref{htoR}, we already know $\overline{\mu}_1,\overline{\mu}_2\in S_{D_0}^1({\bf Y})$. 
Recall the $\E^Y$-nest $(K_n)_{n\in\mathbb{N}}$ of compact sets satisfying $\1_{K_n}\overline{\mu}_1^*\in S_{00}^{(0)}({\bf Y})$ and let  
$\overline{\mu}^n$ 
and 
$h_n(x)$ be as defined in the proof of Lemma~\ref{htoR}. 
Note that $h_n$ increases to $h$ as $n\to\infty$.
Then 
by \eqref{eq:uniformLowBDD} and the proof of Lemma~\ref{htoR}, we have 
 $h_n(x)=R^{Y}(h_n\overline{\mu}^n)(x)+1$
  whenever $\lambda^{\Q}(\overline{\mu}_1)>0$. This expression also yields the fine continuity of $h_n$ on $E$ with respect to ${\bf Y}$. Since $\1_{K_n}\overline{\mu}_1^*\in  S_{D_0}^{1}({\bf Y})$ and 
$\mu_2+N(e^{F_1}(1-e^{-F_2}))\mu_H\in   
  S_{D_0}^{1}({\bf Y})$, we see that for each $n \in \mathbb{N}$
\begin{align*}
h_n-1=R^{Y}(h_n\overline{\mu}^n)\in \F^Y_{\loc}
\end{align*}
is a difference of bounded excessive functions on $E$ with respect to ${\bf Y}$. 
From this observation, $h_n$ is $\E^Y$-quasi continuous on $E$. 
We easily see that $\{h_n\}$ is $\E^Y$-quasi uniformly convergent to $h$ on $E$. Indeed, there exists an $\E^Y$-nest $\{A_l\}$ of compact sets such that all $h_n$ and $h$ are continuous on each $A_l$. 
Here $\E^Y$-nest $\{A_l\}$ implies $\PP_x^Y(\lim_{\ell\to\infty}\tau_{A_{\ell}}=\zeta)=1$ q.e.~$x\in E$. 
Since $h_n$ is increasing on $E_{\partial}$ 
and each $A_l$ is compact, 
$h_n$ uniformly converges to $h$ on $A_l\cup\{\partial\}$, because $h_n(\partial)=h(\partial)=1$. 
 In particular, for each $T\in\,]\,0,\,+\infty\,[$
\begin{align}
\lim_{n\to\infty}\sup_{0\le t \le T}|h_n-h|(X_t)=0 \quad \text{ $\PP^Y_x$-a.s.~on $\{T<\zeta\}$ for q.e. $x \in E$.}
\label{e:hnh}
\end{align}
We apply the generalized Fukushima's decomposition for bounded $h_n-1\in \F_{\loc}^{Y}$ under $\PP_x^{Y}$ (see \cite[Theorem~4.2]{Kw:stochI}
)  
and Lemma~\ref{lem:bddvar} (also \cite[Lemma 5.4.1]{FOT}) that  
\begin{align}
\left\{\begin{array}{ll}
h_n(X_t)-h_n(X_0)=M_t^{h_n}+N_t^{h_n}, \\ 
N_t^{h_n}=-\int_0^t h_n(X_s){\rm d}A_s^{\overline{\mu}^n}
\end{array}\right. 
\label{FukuDecomp-n}
\end{align}
for all $t\in[0,\zeta[$ $\PP_x^{Y}$-a.s.~for~q.e.~$x\in E$, 
where $M_t^{h_n}$ (resp.~$N_t^{h_n}$) 
is the martingale additive functional locally of finite energy (resp.~CAF locally of zero energy) under ${\bf Y}$. 
We define $M_t^{h_n}, N_t^{h_n}$ for $t\geq\zeta$ under 
$\PP_x^Y$-a.s.~for q.e.~$x\in E$
by 
\begin{align*}
N_t^{h_n}:=-\int_0^th_n(X_s){\rm d}A_s^{\overline{\mu}^n},\quad M_t^{h_n}:=h_n(X_t)-h_n(X_0)-N_t^{h_n}.
\end{align*}
Then $M_t^{h_n}$ is a martingale additive functional under 
$\PP_x^Y$~for q.e.~$x\in E$. Indeed, $(M_t^{h_n})_{t\geq0}$ is an AF, and 
for each $t>0$ we see $\EE_x^Y[|M_t^{h_n}|]<\infty$ and 
$\EE_x^Y[M_t^{h_n}]=0$ for q.e.~$x\in E$ under $\overline{\mu}_1,\overline{\mu}_2\in S_{D_0}^1({\bf Y})$. 
 The decomposition \eqref{FukuDecomp-n} also holds for all $t\in\,[\,0,\,+\infty\,[$ $\PP_x^Y$-a.s.~for q.e.~$x\in E$.
We then see that for $m > n>0$, 
\begin{align*}
N_t^{h_m}-N_t^{h_n}&=-
\int_0^th_m(X_s){\rm d}A_s^{\overline{\mu}^m}+
\int_0^th_n(X_s){\rm d}A_s^{\overline{\mu}^n}\\
&=-\int_0^t(h_m\1_{K_m}-h_n\1_{K_n})(X_s){\rm d} A_s^{\overline{\mu}_1^*}+
\int_0^t(h_m-h_n)(X_s){\rm d} A_s^{\mu_2+N(e^{F_1}(1-e^{-F_2})\mu_H}.
\end{align*}
Hence 
\begin{align*}
|N_t^{h_m}-N_t^{h_n}|&\leq \int_0^t(h-h_n\1_{K_n})(X_s){\rm d} A_s^{\overline{\mu}_1^*}+\int_0^t(h-h_n)(X_s){\rm d} A_s^{\mu_2+N(e^{F_1}(1-e^{-F_2}))\mu_H}\\
&\leq \int_0^t(h-h_n\1_{K_n})(X_s){\rm d} A_s^{\overline{\mu}_1^*+\mu_2+N(e^{F_1}(1-e^{-F_2}))\mu_H}.
\end{align*}
Thus, using  \eqref{e:hnh}, we have
\begin{align*}
\sup_{0\le t \le T}\left|M_t^{h_m}-M_t^{h_n}\right|&\le 
\sup_{0\le t \le T}|h_m-h_n|(X_t) + \sup_{0\le t \le T}\left|N_t^{h_m}-N_t^{h_n}\right|\\
&\le \sup_{0\le t \le T}|h-h_n|(X_t) +\int_0^T\left|h-\1_{K_n}h_n\right|(X_s){\rm d}A_s^{\overline{\mu}_1^*+{\mu}_2+e^{\|F_1\|_{\infty}}N(F_2)\mu_H} 
\\ 
&\le \sup_{0\le t \le T}|h-h_n|(X_t)(1+A_T^{\overline{\mu}_1^*+{\mu}_2
+e^{\|F_1\|_{\infty}}N(F_2)\mu_H
})
\\
&\hspace{2cm}
+\|h\|_{\infty}\int_0^T\1_{K_n^c}(X_s){\rm d}A_s^{\overline{\mu}_1^*+\mu_2+e^{\|F_1\|_{\infty}}N(F_2)\mu_H}
\\
&\longrightarrow ~0, \quad \text{as }~ m, n \to \infty,
\quad \text{$\PP^Y_x$-a.s.~on $\{T<\zeta\}$ for q.e. $x \in E$}. 
\end{align*}
Similarly, in view of the Lebeague's dominated convergence and $\overline{\mu}_1^*,\overline{\mu}_2\in S_{D_0}^1({\bf Y})$, we can deduce 
\begin{align*}
|M_{\zeta}^{h_m}-M_{\zeta}^{h_n}|&=|N_{\zeta}^{h_m}-N_{\zeta}^{h_n}|\\
&\leq 
\int_0^{\zeta}\left|h-\1_{K_n}h_n\right|(X_s){\rm d}A_s^{\overline{\mu}_1^*+{\mu}_2+e^{\|F_1\|_{\infty}}N(F_2)\mu_H}
&\to0\quad\text{ as }\quad m,n\to\infty
\end{align*}
$\PP_x^Y$-a.s.~for q.e.~$x\in E$. Therefore, for each $t\in\,[\,0,\,+\infty\,]$, $\{M_t^{h_n}\}$ and $\{N_t^{h_n}\}$ converge as $n\to\infty$ under $\PP_x^Y$ for q.e.~$x\in E$. Set $M_t^h:=\lim_{n \to \infty}M_t^{h_n}$ and $N_t^h:=\lim_{n \to \infty}N_t^{h_n}$. Then $M^h$ (resp. $N^h$) 
satisfies 
that \eqref{FukuDecomp} holds for 
all $t\in\,[\,0,\,+\infty\,[$ $\PP_x^{Y}$-a.s.~for q.e.~$x\in E$, and hence it 
is a square integrable martingale additive functional locally of finite energy (resp.~CAF locally of zero energy) under $\PP_x^Y$ for q.e. $x \in E$. 
Moreover, from \eqref{FukuDecomp-n},  we have for $p \in \mathbb{N}$
\begin{align*}
\EE_x^{Y}\left[\left|M_t^{h_n}\right|^p \right]&\leq 2\cdot 4^{p-1}\|h_n\|_{\infty}^p+2^{p-1}\EE_x^{Y}\left[\left(\int_0^t h_n(X_s){\rm d}
A_s^{\overline{\mu}_1^*+{\mu}_2+e^{\|F_1\|_{\infty}}N(F_2)\mu_H}\right)^p\,\right]\\
& \le 2\cdot 4^{p-1}\|h\|_{\infty}^p+2^{p-1}p!\left\|R^{Y}(h(\overline{\mu}_1^*+{\mu}_2+e^{\|F_1\|_{\infty}}N(F_2)\mu_H))\right\|_{\infty}^p=:C_p<\infty
\end{align*}
for~q.e.~$x\in E$, where we use \cite[Lemma~2.2]{DKK2010}. 
 So $\{M_t^{h_n}\}_{n \in \mathbb{N}}$ is uniformly $\PP_x^{Y}$-integrable for q.e.~$x\in E$. In particular, we have 
\begin{align*}
R^{Y}(\mu_{\<h_n\>})(x)=\EE_x^{Y}\left[\<M^{h_n}\>_{\infty}\right]=\EE_x^{Y}\left[\left|M_{\infty}^{h_n}\right|^2 \right] \leq C_2
\end{align*}
for q.e. $x\in E$, consequently $\sup_{x \in E}R^{Y}(\mu_{\<h_n\>})(x)\leq C_2$ because of the fine continuity of $x \mapsto R^{Y}(\mu_{\<h_n\>})(x)$. Then we have $\mu_{\<h_n\>} \in S_{D_0}^1({\bf Y})$. Therefore, we can redefine $M^{h_n}$ (resp. $N^{h_n}$) as the martingale additive functional locally of finite energy (resp. CAF locally of zero energy) in the strict sense under ${\bf Y}$, and hence the following Fukushima's decomposition in the strict sense 
holds (\cite[Theorem~6.2]{KimKuwaeTawara}):
$$
h_n(X_t)-h_n(X_0)=M^{h_n}_t +N^{h_n}_t
$$
for all $t \in \,[\,0,\,+\infty\,[$~ $\PP_x^Y$-a.s. for all $x \in E$.  Moreover, by combining the same way as above and Fatou's lemma, we can obtain 
$$
\sup_{x \in E}R^{Y}(\mu_{\<h\>})(x) \le  8\|h\|_{\infty}^2+4\left\|R^{Y}(h(\overline{\mu}_1^*+{\mu}_2+e^{\|F_1\|_{\infty}}N(F_2)\mu_H))\right\|_{\infty}^2<\infty,
$$
that is, $\mu_{\<h\>} \in S_{D_0}^1({\bf Y})$. Using this fact, we can redefine both $M^{h}$ and $N^{h}$ in the strict sense under ${\bf Y}$, and which leads us to the conclusion. 
\end{pf}

Under $\lambda^{\Q}(\overline{\mu}_1)>0$, $u=0$, 
$\overline{\mu}_1^*=\mu_1+N(e^{F_1}-1)\mu_H\in S_{N\!K_1}^1({\bf X})$ 
 and $\overline{\mu}_2=\mu_2+N(F_2)\mu_H\in S_{D_0}^1({\bf X})
 $, we 
define a martingale additive functional of ${\bf Y}$ by $M_t:=\int_0^th({X}_{s-})^{-1}{\rm d}M_s^{h}$ and denote by $L_t^h:=\text{\rm Exp}(M)_t$ the Dol\'eans-Dade exponential of $M_t$, that is, $L_t^h$ is the unique solution of $L_t^h=1+\int_0^tL_{s-}^h{\rm d}M_s$ for $t\in[0,\zeta[$, $\PP_x^{Y}$-a.s.~for q.e.~$x\in E$. Then we see
\begin{align*}
L_t^h=\exp \left(M_t-\frac{\,1\,}{2}\<M^c\>_t\right)
\prod_{0<s\le t}\frac{h(X_s)}{h(X_{s-})}\exp\left(1-\frac{h(X_s)}{h(X_{s-})}\right), 
\end{align*}
where $M_t^c$ is the continuous part of $M_t$. Note that 
$L_t^h$ is also a multiplicative functional under ${\bf Y}$. 
By Lemma~\ref{Fukudecompforh} and It\^o's formula applied to the semimartingale $h(X_t)$ in \eqref{FukuDecomp} with 
the function $\log x$, we see that $L_t^h$ has the following expression:
\begin{align}
L_t^h=\frac{h(X_t)}{h(X_0)}\exp\left(A_t^{\overline{\mu}} \right). \label{expressionL}
\end{align}
Let ${\bf Y}^{h}:=(\Omega, \wt{\mathscr{F}}_{\infty}, \wt{\mathscr{F}}_{t},  \wt{X}_{t}, {\PP}_{x}^{Y,h}, \zeta^{Y,h})$ the transformed process of ${\bf Y}$ by $L_t^h$: 
\begin{align*} 
\EE_x^{Y,h}[f(X_t)]=\EE_x^{Y}[L_t^h f(X_t)], \quad  f\in \mathscr{B}_b(E).
\end{align*}
Set $\ell:=-\log h=-\log(R^{Y}(h{\overline{\mu}})+1)$. Since $R^{Y}(h{\overline{\mu}})\in {\F}_{\loc}^{Y}$, we have 
$\ell\in \dot{\F}_{\loc}^{Y}$ and $h=e^{-\ell}$. 
It is easy to see the fine continuity of $\ell$ on $E$ with respect to ${\bf Y}$. 
Then we can see that the transform $L_t^h$ belongs to the class of Girsanov transforms considered in Section~\ref{sec:AofGandSubcriticality} (or in \cite{CZ2002}) under ${\bf Y}$. In particular,  
 ${\bf Y}^{h}$ is an $h^2
 \m$-symmetric Hunt process on $E$.  
 From Theorem \ref{FeynmanKacMod3} (cf. \cite[Theorem~3.2]{KK:AnalChara}) 
 the associated Dirichlet form 
 $(\E^{Y,h},\F^{Y,h})$ on $L^2(E;h^2
 \m)$ is identified by 
\begin{align*}
{\E}^{Y,h}(f,f)& =\frac{\,1\,}{2}\int_{E}h^2
{\rm d}\mu_{\< f\>}^{c} +\int_{E\times E}(f(x)-f(y))^{2}e^{F(x,y)}h(x)h(y)
J({\rm d}x{\rm d}y)  
\end{align*} 
for $f \in \F~ (=\F^{Y,h})$.  Let ${\bf Y}^h$ be the transformed process of ${\bf Y}$ by $L_t^h$. Then $(\E^{Y,h},\F^{Y,h})$ is nothing but the Dirichlet form on $L^2(E;h^2\m)$ associated with ${\bf Y}^h$. Moreover we see that $(\E^{Y,h}, \F^{Y,h})$ is equivalent to $(\E,\F)$ 
because of \eqref{eq:gaugeable} and 
\begin{align}
e^{-2\|R^Y\overline{\mu}_2\|_{\infty}}e^{-\|F\|_{\infty}
}\E^{(c)}(f,f)\leq (\E^{Y,h})^{(c)}(f,f)
\leq \|h\|_{\infty}^2e^{\|F\|_{\infty}
}\E^{(c)}(f,f)
\label{equivalenceofDs}
\end{align}
and 
\begin{align}
e^{-2\|R^Y\overline{\mu}_2\|_{\infty}}e^{-\|F\|_{\infty}
}J({\rm d}x{\rm d}y)\leq 
J^{Y,h}({\rm d}x{\rm d}y)\leq \|h\|_{\infty}^2e^{\|F\|_{\infty}
}J({\rm d}x{\rm d}y).\label{equivalenceofDsJ}
\end{align} 
Here $J^{Y,h}({\rm d}x{\rm d}y):=e^{F(x,y)}h(x)h(y)
J({\rm d}x{\rm d}y)$ is the jumping measure of $(\E^{Y,h},\F^{Y,h})$.
Therefore, under {\bf (A.1)}
(resp.~{\bf (A.2)}), $(\E^{Y,h},\F^{Y,h})$ admits a heat kernel $p_t^{Y,h}(x,y)$ (resp.~a heat kernel $p_t^{Y,h}(x,y)$ in the strict sense)  on $]\,0,\,+\infty\,[\times E\times E$ such that 
\begin{align*}
&p_t^{Y,h}(x,y)\precsim \phi_2(t,x,y)\quad \m\text{-a.e.~}\quad d(x,y)\in \mathbb{T}_t\\
&\text{(resp.~$\phi_1(t,x,y)\precsim p_t^{Y,h}(x,y)\precsim \phi_2(t,x,y)$, \quad $d(x,y)\in\mathbb{T}_t$)}.
\end{align*} 
\vskip 0.2cm
Now, we give the proofs of Theorem \ref{thm:mainresult1}, Corollary \ref{cor:mainresult} and  
Theorem~\ref{thm:shortstability}.

\begin{apf}{Theorem~\ref{thm:mainresult1}}   
\eqref{item:mainresult1}, \eqref{item:mainresult1*}: 
Recall that ${\bf U}$ is the Girsanov transformed process 
by the supermartingale multiplicative 
functional $U_t=\Exp(M^{e^{[u]}-1}+M^{-u,c})_t$.
Suppose $\lambda^{\Q}(\overline{\mu}_1)>0$, 
equivalently $\lambda^{\Q}(\overline{\mu}_1^*)>0$ (see 
\cite[Lemma~5.1]{KK:AnalChara}). 
Since $\mu_{\<u\>}\in S_{N\!K_{\infty}}^1({\bf X})$, we have 
$e^{-2u}\overline{\nu}_1+e^{-u}N(e^{-u}(e^{F_1}-1))\mu_H=
e^{-2u}(\overline{\nu}_1+N(e^{[u]}(e^{F_1}-1))\mu_H)=e^{-2u}\overline{\mu}_1^*\in S_{N\!K_1}^1({\bf U})$ and 
$e^{-2u}\overline{\nu}_2+e^{-u}N(e^{-u}F_2)\mu_H=e^{-2u}(\overline{\nu}_2+N(e^{[u]}F_2)\mu_H)\in S_{D_0}^1({\bf U})$ by 
\cite[Corollary~5.1(1),(5)]{KimKurKuwae:Gauge}. The Dirichlet form 
$(\E^U,\F^U)$ on $L^2(E;e^{-2u};\m)$ associated to 
 ${\bf U}$ 
satisfies \eqref{eq:comparisonmeasure}, \eqref{eq:comparisondiffusion} and \eqref{eq:comparisonjumpdensity} for some $C_E>0$.  
Under {\bf (A.1)}
(resp.~{\bf (A.2)}), {\bf U} admits  
a heat kernel $p_t^U(x,y)$ such that $p_t^U(x,y)
\precsim \phi_2(t,x,y)$ $\m$-a.e. holds 
(resp.~a heat kernel $p_t^U(x,y)$ in the strict sense such that 
$\phi_2(t,x,y)\precsim p_t^U(x,y)
\precsim \phi_2(t,x,y)$, $d(x,y)\in\mathbb{T}_t$ holds.   
So, we may and do assume $u=0$ to the end of the proof.  

We apply Lemma~\ref{htoR} and 
\eqref{eq:gaugeable}. 
Recall that the additive functional $A$ is given by $A_t:=A_t^{\mu}+A_t^F$ with $A_t^{\mu}=A_t^{\mu_1}-A_t^{\mu_2}$ and $A_t^F=A_t^{F_1}-A_t^{F_2}$. 
Since
\begin{align*}
P_t^A\!f(x)=\EE_x^Y\left[e^{A_t^{\overline{\mu}}}f(X_t)\right]=
h(x)\EE_x^{Y,h}\left[\frac{f(X_t)}{h(X_t)}\right],
\end{align*}
the integral kernel 
$p_t^A(x,y)$ can be represented as 
\begin{equation}
p_t^A(x,y) =h(x)h(y)p_t^{Y,h}(x,y).
\label{pAph}
\end{equation}
So we can conclude the assertion, because 
$0<e^{-\|R^Y\overline{\mu}_2\|_{\infty}}\leq h(x)\leq\|h\|_{\infty}$ and $p_t^{Y,h}(x,y)\precsim \phi_2(t,x,y)$ $\m$-a.e. $d(x,y)\in\mathbb{T}_t$ (resp.~$\phi_1(t,x,y)\precsim p_t^{Y,h}(x,y)\precsim \phi_2(t,x,y)$, $d(x,y)\in\mathbb{T}_t$) under 
{\bf (A.1)}
(resp.~{\bf (A.2)}).  
 Therefore, we obtain the conclusion.
 
\eqref{item:mainresult2}: 
Let ${\sf d}:=\{(x,y)\in E\times E\mid R(x,y)=\infty\}$. 
Since $p_t^A(x,y) \lesssim p_t(x,y)$, we see that there exists $C>0$ such that 
\begin{align}
R^A(x,y)\leq CR(x,y)<\infty\quad\text{ for }(x,y)\in E\times E\setminus{\sf d}.\label{eq:subcriticaloriginal}
\end{align}
Thus we have $R^A(x,y)<\infty$ for $\m$-a.e.~$y\in E^x$
where $E^x=\{y\in E\mid (x,y)\notin {\sf d}\}$. 
Since $\m(E\setminus E^x)=0$, $R^A(x,y)<\infty$ $\m$-a.e.~$y\in E$. 
By Theorem~\ref{thm:analgauge1}, we obtain $\lambda^{\Q}(\overline{\mu}_1)>0$. 
\end{apf}

Next we prove Corollary~\ref{cor:mainresult}. The proof is a mimic of the proof of Theorem~\ref{thm:mainresult1} except the representation of CAF locally of zero energy part in \eqref{FukuDecomp}.

\begin{apf}{Corollary~\ref{cor:mainresult}} 
\eqref{item:cormainresult1},\eqref{item:cormainresult1*}: 
As noted in the proof of Theorem~\ref{thm:mainresult1}, 
we may and do assume $u=0$.  
First we assume $\overline{\mu}_2=\mu_2+N(F_2)\mu_H\in S_D^1({\bf X})
$.  
Let ${\bf Y}^{(\alpha)}$ be the $\alpha$-subprocess of ${\bf Y}$ killed at rate $\alpha\m$. Set the function 
$h_{\alpha}(x)=\EE_x^{Y^{(\alpha)}}[e^{A^{\overline{\mu}}_{\zeta}}]$ for $x \in E_{\partial}$.  
Suppose $\lambda^{\Q_{\alpha}}(\overline{\mu}_1)>0$, equivalently, $\lambda^{\Q_{\alpha}}(\overline{\mu}_1^*)>0$ in view of \cite[Lemma 5.1]{KK:AnalChara}. Then $h_{\alpha}$ satisfies 
$0<e^{-\|R_{\alpha}^Y\overline{\mu}_2\|_{\infty}}
\leq  h_{\alpha}(x)\leq \sup_{x\in E}h_{\alpha}(x)<\infty$ and $h_{\alpha}-1=R^Y_{\alpha}(h_{\alpha}\overline{\mu}) \in \F_{\loc}^Y$ by virtue of Lemma \ref{htoR} and \eqref{h-1Floc}. By Lemma~\ref{htoR}, $\overline{\mu}_1\in S_D^1({\bf Y})$ and  
$\overline{\mu}_2\in S_D^1({\bf Y})$. 
Using  Lemma \ref{lem:bddvar}, by following the proof of Lemma \ref{Fukudecompforh} line by line,
we have the following strict decomposition for $h_{\alpha}-1\in \F_{\loc}^{Y}$ under $\PP_x^{Y}$:
\begin{align*}
\left\{\begin{array}{ll}
h_{\alpha}(X_t)-h_{\alpha}(X_0)=M_t^{h_{\alpha}}+N_t^{h_{\alpha}}, \\ 
N_t^{h_{\alpha}}=\alpha\int_0^t\left(h_{\alpha}-1\right)(X_s){\rm d}s-\int_0^t h_{\alpha}(X_s){\rm d}A_s^{\overline{\mu}}
\end{array}\right. 
\text{for all $t\in\,[\,0,\,+\infty\,[$ $\PP_x^{Y}$-a.s.~for~all~$x\in E$. }
\label{FukuDecompeps}
\end{align*}
Let $M_t^{(\alpha)}:=\int_0^th_{\alpha}({X}_{s-})^{-1}{\rm d}M_s^{h_{\alpha}}$ and denote by $L^{h_{\alpha}}_t:={\rm Exp} (M^{(\alpha)})_t$ 
the Dol\'eans-Dade exponential of $M^{(\alpha)}_t$.  
By using It\^o's formula applied to $h_{\alpha}(X_t)$ with the function $\log x$, we then have 
\begin{align*}
L^{h_\alpha}_t&=\exp \left(M^{(\alpha)}_t-\frac{\,1\,}{2}\<M^{(\alpha),c}\>_t\right)\prod_{0<s\leq t}\frac{h_{\alpha}(X_s)}{h_{\alpha}(X_{s-})}\exp\left(1-\frac{h_{\alpha}(X_s)}{h_{\alpha}(X_{s-})}\right)\nonumber  \\
&=\frac{h_{\alpha}(X_t)}{h_{\alpha}(X_0)}\exp \left(A_t^{\overline{\mu}}+\alpha\int_0^t\left(h_{\alpha}^{-1}-1\right)(X_s){\rm d}s\right). 
\end{align*} 
By noting again that the transform $L_t^{h_\alpha}$ belongs to the class of Girsanov transforms considered in Section~\ref{sec:AofGandSubcriticality}, the transformed process ${\bf Y}^{h_{\alpha}}:=(X_{t}, {\PP}_{x}^{Y,h_{\alpha}})$ of ${\bf Y}$ by $L^{h_\alpha}_t$ becomes an $h^2_{\alpha}\m$-symmetric Hunt process on $E$ and the associated Dirichlet form 
 $(\E^{Y,h_{\alpha}},\F^{Y,h_{\alpha}})$ is 
also equivalent to $(\E,\F)$ in view of 
\eqref{equivalenceofDs} and 
\eqref{equivalenceofDsJ}. 
Let $p_t^{Y,h_{\alpha}}(x,y)$ be the heat kernel associated to $(\E^{Y,h_{\alpha}},\F^{Y,h_{\alpha}})$. By virtue of 
the above expression for $L_t^{h_{\alpha}}$, we see 
\begin{align*}
P_t^A\!f(x)=\EE_x^Y[e^{A_t^{\overline{\mu}}}f(X_t)]=h_{\alpha}(x)
\EE_x^{Y,h_{\alpha}}\left[\exp\left(-\alpha\int_0^t(h_{\alpha}^{-1}-1)(X_s){\rm d}s\right)\frac{f(X_t)}{h_{\alpha}(X_t)}\right].
\end{align*}
Then we have 
\begin{equation}
e^{-k_1t}h_{\alpha}(x)h_{\alpha}(y)p_t^{Y,h_{\alpha}}(x,y) \le p_t^{A}(x,y) \le e^{k_2t}h_{\alpha}(x)h_{\alpha}(y)p_t^{Y,h_{\alpha}}(x,y),
\label{epsrelation}
\end{equation}
where $k_1=k_1(\alpha):=\alpha (\ell_1^{-1}-1)$ and $k_2=k_2(\alpha):=\alpha (1-\ell_2^{-1})$ for $\ell_1=\inf_{x \in E}h_{\alpha}(x)$ and $\ell_2=\sup_{x \in E}h_{\alpha}(x)$. 
Since $p_t^{Y,h_{\alpha}}(x,y) \precsim \phi_2(t,x,y)$ $\m$-a.e. $d(x,y)\in\mathbb{T}_t$ (resp.~$\phi_1(t,x,y)\precsim p_t^{Y,h_{\alpha}}(x,y) \precsim \phi_2(t,x,y)$, $d(x,y)\in\mathbb{T}_t$) holds under {\bf (A.1)}
(resp.~{\bf (A.2)}), we see from \eqref{epsrelation} that $p_t^{A}(x,y) \precsim_k \phi_2(t,x,y)$ $\m$-a.e. $d(x,y)\in\mathbb{T}_t$ (resp.~$\phi_1(t,x,y)\precsim_k p_t^{A}(x,y) \precsim_k \phi_2(t,x,y)$, $d(x,y)\in\mathbb{T}_t$) ($k:=\max\{k_1,k_2\}>0$) holds 
under {\bf (A.1)} (resp.~{\bf (A.2)}). 

\vskip 0.1cm

\eqref{item:mainresult2}: 
Suppose $p_t^A(x,y) \lesssim_k p_t(x,y)$ for $\alpha >k$. Then there exists $C>0$ such that 
for each $x\in E$ 
\begin{align*}
R_{\alpha}^A(x,y)\leq CR_{\alpha -k}(x,y)<\infty\quad \text{ for }\m\text{-a.e.~}y\in E.
\end{align*} 
Thus $R_{\alpha}^A(x,y)<\infty$ $\m$-a.e.~$y\in E$ for each $x\in E$.
Applying Theorem~\ref{thm:analgauge1} to ${\bf X}^{(\alpha)}$, we have $\lambda^{\Q_{\alpha}}(\overline{\mu}_1)>0$.  
\end{apf}
 
Finally, we give the proof of Theorem~\ref{thm:shortstability}. To do this, we need the following lemma.

\begin{lem}\label{lem:extended}
For each $\nu\in S_{E\!K}^1({\bf X})$, there exists $\alpha>0$ such that $\nu\in S_{N\!K_1}^1({\bf X}^{(\alpha)})$.
\end{lem}
\begin{pf}
By $\nu \in S_{E\!K}^1({\bf X})$, we see that $\sup_{x\in E}R_{\alpha}\nu(x)<1$ for some large $\alpha >0$. 
From this, we can see that $\sup_{x\in E}R_{\alpha}\1_{B\cup K^c}\nu(x)<1$ for any Borel set $K$ of $E$ with $\nu(K)<\infty$ and for any $\nu$-measurable set $B \subset K$ with $\nu(B)<\delta$~ $(\delta >0)$, which implies that 
$\nu$ is in the family of semi-Green-tight Kato class measures in the sense of Chen with respect to ${\bf X}^{(\alpha)}$. Therefore, 
 $\nu \in  S_{N\!K_1}^1({\bf X}^{(\alpha)})$. (See Remark \ref{AdvantageofNV}(1).)
\end{pf}
\begin{apf}{Theorem~\ref{thm:shortstability}} 
Since $\mu_{\<u\>}\in S_K^1({\bf X})$, we have 
$e^{-2u}\overline{\mu}_1^*=
e^{-2u}(\overline{\nu}_1+N(e^{[u]}(e^{F_1}-1))\mu_H)\in S_{E\!K}^1({\bf U})$ and $e^{-2u}(\mu_2+N(F_2)\mu_H)\in S_{D}^1({\bf U})$ by \cite[Lemma~4.1(3)]{KimKurKuwae:Gauge}. By Lemma~\ref{lem:extended},  there exists large $\alpha >0$ such that $e^{-2u}\overline{\mu}_1^*\in S_{N\!K_1}^1({\bf U}^{(\alpha)})$. 
Applying Stollmann-Voigt's inequality to the Girsanov transformed process ${\bf U}$, we have that for large enough $\alpha>0$
\begin{align*}
\lambda^{\Q_{\alpha}}(\overline{\mu}_1^*)&=\inf\left\{\Q_{\alpha}(f,f)\;\Biggl|\; f\in \C,\: \int_Ef^2{\rm d}\overline{\mu}_1^* =1 \right\}\\
&=\inf\left\{\E_{\alpha}^{U}(fe^u,fe^u)-\int_Ef^2{\rm d}\overline{\nu}-\int_{E\times E}
f(x)f(y)(e^{F(x,y)}-1)N(x,{\rm d}y)\mu_H({\rm d}x) \right. \\
&\hspace{7cm}\left.\;\Biggl|\; f\in \C,\: \int_Ef^2{\rm d}\overline{\mu}_1^* =1 
\right\}
\end{align*}
\begin{align*}
&\geq \inf\left\{\E_{\alpha}^{U}(fe^u,fe^u)-\int_Ef^2{\rm d}\overline{\nu}_1-\int_E
f(x)^2N(e^{}[u](e^{F_1}-1))(x)\mu_H({\rm d}x) \right. \\
&\hspace{7cm}\left.\;\Biggl|\; f\in \C,\: \int_Ef^2{\rm d}\overline{\mu}_1^* =1 
\right\}\\
&= \inf\left\{\E_{\alpha}^{U}(fe^u,fe^u)-\int_Ef^2{\rm d}\overline{\mu}_1^*\;\Biggl|\; f\in \C,\: \int_Ef^2{\rm d}\overline{\mu}_1^* =1 
\right\}\\
&\geq \left(1-\|R_{\alpha}^Ue^{-2u}\overline{\mu}_1^* \|_{\infty}\right)
\inf\left\{\E_{\alpha}^{U}(fe^u,fe^u)\;\Biggl|\; f\in \C,\: \int_Ef^2{\rm d}\overline{\mu}_1^* =1 
\right\}\\
&\geq \frac{1-\|R_{\alpha}^Ue^{-2u}\overline{\mu}_1^* \|_{\infty}}{\|R_{\alpha}^Ue^{-2u}\overline{\mu}_1^* \|_{\infty}}>0,
\end{align*}
which is equivalent to $\lambda^{\Q_{\alpha}}(\overline{\mu}_1)>0$ (\cite[Lemma~5.1]{KK:AnalChara}). 

Since {\bf (A.1)} (resp.~{\bf (A.2)}) holds for {\bf X}, we can deduce that {\bf (A.1)} 
(resp.~{\bf (A.2)}) holds for {\bf U},  
because $\mu_{\<f\>}^{U,c}({\rm d}x)=e^{-2u(x)}\mu_{\<f\>}^{c}({\rm d}x)$ is the energy measure of continuous part of $(\E^U,\F^U)$ and $J^U({\rm d}x{\rm d}y)=e^{-u(x)-u(y)}J({\rm d}x{\rm d}y)$ is the jumping measure of $(\E^U,\F^U)$.
Note that $p_t^A(x,y)=(p_t^U)^{A^{\overline{\nu},F}}(x,y)$, 
where $A^{\overline{\nu},F}_t:=A_t^{\mu}+\int_0^tN(e^{[u]}-[u]-1)(X_s){\rm d}H_s+\frac{\,1\,}{2}\<M^{-u,c}\>_t$.

Applying Corollary~\ref{cor:mainresult}\eqref{item:cormainresult1}
\eqref{item:cormainresult1*} to {\bf U}, we can construct an 
integral kernel (resp.~an integral kernel in the strict sense) $p_t^A(x,y)=(p_t^U)^{A^{\overline{\nu},F}}(x,y)$ such that   
$p_t^A(x,y)\precsim_k \phi_2(t,x,y)$ $\m$-a.e. $d(x,y)\in\mathbb{T}_t$ (resp.~$\phi_1(t,x,y)\precsim p_t^A(x,y)\precsim_k \phi_2(t,x,y)$, $d(x,y)\in\mathbb{T}_t$) under {\bf (A.1)}
(resp.~{\bf (A.2)}) for some $k>0$. 
Therefore we obtain the conclusion. 
\end{apf}

\section{Examples}\label{sec:example}

We will use the following notations in this section. 
\begin{df}
	{\rm Let $g:\,]\,0,\,+\infty\,[ \,\to\, [\,0,\,+\infty\,[$, and $a\in\,]\,0,\, +\infty\,]$,  $\beta_1, \beta_2>0$, and $c, C>0$.
		\begin{enumerate}
			\item[(1)] We say that $g$  satisfies $L(\beta_1, c)$, the (global) lower scaling condition with index  $\beta_1$, if 
			$$ \frac{g(R)}{g(r)} \geq c \left(\frac{R}{r}\right)^{\beta_1} \quad \text{for all} \quad r\leq R< \infty.$$
			\item[(2)] We say that $g$ satisfies $U(\beta_2, C)$, 
			the (global) upper scaling condition with index  $\beta_2$,
			 if
			$$ \frac{g(R)}{g(r)} \leq C \left(\frac{R}{r}\right)^{\beta_2} \quad \text{for all} \quad r\leq R< \infty.$$
					\end{enumerate}
}\end{df}

\subsection{
Brownian Motion in $\R^d$ 
}\label{ex:Brownian}
Let ${\bf X}=(\Omega, X_t,\PP_x)$ be $d$-dimensional Brownian motion
on $\R^d$. Note that {\bf X} satisfies 
{\bf (A.2)}
with $\mathbb{T}=[\,0,\,+\infty\,[^{\,2}$.
A signed Borel measure $\mu$ on $\R^d$ is said to be
\emph{of Kato class} if
\begin{align*}
\lim_{r\to0}\sup_{x\in\R^d}\int_{|x-y|<r}\frac{ |\mu
|(\d y)}{|x-y|^{d-2}}&=0
\quad \hbox{when }\quad d\geq3,\\
\hspace{-1cm}\lim_{r\to0}\sup_{x\in\R^d}\int_{|x-y|<r}(\log|x-y|^{-1})
|\mu|(\d y)&=0
\quad \hbox{when }\quad  d=2, \\
\sup_{x\in\R^d}
 |\mu| (B(x,1))&<\infty 
\quad
\hbox{when } \quad d=1.
\end{align*}
Here $| \mu|:=\mu_++\mu_-$ is the total variation measure of $\mu$.
 A signed Borel measure $\mu$ on $\R^d$ is said to be
\emph{of local Kato class} if ${\bf 1}_K\mu$ is of Kato class for every
compact subset $K$ of $\R^d$. By definition, any measure $\mu$ of
local Kato class is always a signed Radon measure. Denote by
${\bf K}_{d}$ (resp. ${\bf K}_{d}^{loc}$) the family of non-negative measures 
of Kato class (resp.
local Kato class)
on $\R^d$. It is essentially proved in
\cite{AizSim:Brown} that a 
non-negative 
 measure $\mu$ is in Kato class $
{\bf K}_{d}$ if and only if $\mu$ is a smooth measure in the strict
sense and
\begin{align*}
\lim_{t\to0}\sup_{x\in\R^d}\int_{\R^d}\left(\int_0^tp_s(x,y)ds\right)\mu(\d y)=
\lim_{t\to0}\sup_{x\in\R^d}\EE_x[A_t^{\mu}]=0,
\end{align*}
where $A^{\mu}$ is a PCAF of $X$ admitting no exceptional set
associated to $\mu$ under Revuz correspondence. That is, 
${\bf K}_{d}=S_K^1({\bf X})$. 

Take $\phi\in C_{\infty}(\R^d)\cap H^1(\R^d)_e$ with $|\nabla\phi|^2{\rm d}x\in{\bf K}_{d}$.  It is proved in \cite[Theorem~2.8]{GRS1992} that 
there exists a jointly continuous integral kernel $p_t^{\phi}(x,y)$ such that 
\begin{align*}
\EE_x[e^{N_t^{\phi}}f(X_t)]=\int_{\Rd}p_t^{\phi}(x,y)f(y){\rm d}y,\qquad f\in \mathscr{B}_b(\Rd).
\end{align*} 
Our Theorem~\ref{thm:shortstability} with 
$\phi_t^i(x,y)
=\frac{1}{(2\pi t)^{d/2}}\exp\left(-|x-y|^2/t\right)$ $(i=1,2)$ 
 tells us that 
there exist constants $C,C_1,C_2,k>0$ with $C_1\leq C_2$ such that 
\begin{align}
\frac{C^{-1}e^{-kt}}{(2\pi t)^{d/2}}e^{-\frac{|x-y|^2}{C_1t}}\leq p_t^{\phi}(x,y)\leq 
\frac{Ce^{kt}}{(2\pi t)^{d/2}}e^{-\frac{|x-y|^2}{C_2t}},\quad x,y\in\Rd, \quad t>0,
\label{eq:HKExKato}
\end{align}
which is blunter than the estimate in \cite[Theorem~2.9]{GRS1992}. 
However, Theorem~\ref{thm:mainresult1} with 
$\phi_i(t,x,y)
=\frac{1}{(2\pi t)^{d/2}}\exp\left(-|x-y|^2/t\right)$ $(i=1,2)$ tells us that 
for $d\geq3$ and $|\nabla\phi|^2{\rm d}x\in {\bf K}_{d}$, $\lambda^{\Q}\left(\frac{\,1\,}{2}|\nabla\phi|^2{\rm d}x\right)>0$ if and only if there exist $C,C_1,C_2>0$ with $C_1\leq C_2$ such that \eqref{eq:HKExKato} holds with $k=0$. 
Note that for $\phi\in C_{\infty}(\R^d)\cap H^1(\R^d)_e$, 
$|\nabla\phi|^2{\rm d}x\in {\bf K}_{d}$ is equivalent to $|\nabla\phi|^2{\rm d}x\in {\bf K}_{d}^{\infty}$
by \cite[Lemma~5.1(2)]{DKK2010} and \cite[Proposition~1]{Zhao:subcrit}, 
where  
\begin{align*}
{\bf K}_{d}^{\infty}:=\left\{\nu\in {\bf K}_d\;\left|\; \lim_{R\to\infty}\sup_{x\in\R^d}
\int_{|x|\geq R}\frac{\nu({\d} y)}{|x-y|^{d-2}}=0\right.\right\}
\end{align*}
is the class of Green-tight measures of Kato class in the sense of Zhao. 
It is easy to see ${\bf K}_d^{\infty}=S_{K_{\infty}}^1({\bf X})$, where 
\begin{align*}
S_{K_{\infty}}^1({\bf X}):=\left\{\nu\in S_K^1({\bf X})\;\left|\;\forall \eps>0, \exists K:\text{ compact such that }  \sup_{x\in E}\int_{K^c}\frac{\nu({\d} y)}{|x-y|^{d-2}}<\eps \right.\right\}.
\end{align*}
Hence ${\bf K}_{d}^{\infty}=S_{N\!K_{\infty}}^1({\bf X})$ by \cite[Proposition~4.1]{KimKurKuwae:Gauge}.  

Here $\Q$ is the quadratic form defined by 
\begin{align*}
\Q(f,g)&=\frac{\,1\,}{2}\int_{\Rd}\<\nabla f(x),\nabla g(x)\>{\rm d}x+\frac{\,1\,}{2}\int_{\Rd}
g(x)\<\nabla f(x),\nabla\phi(x)\>{\rm d}x\\
&\hspace{1cm}+\frac{\,1\,}{2}\int_{\Rd}
f(x)\<\nabla g(x),\nabla\phi(x)\>{\rm d}x\quad f,g\in H^1(\Rd)
\end{align*}
and 
$$
\lambda^{\Q}\left(\frac{\,1\,}{2}|\nabla\phi|^2{\rm d}x\right):=\inf\left\{ \Q(f,f)\;\Biggl|\; f\in C_0^{\infty}(\Rd),\quad \frac{\,1\,}{2}\int_{\Rd}f(x)^2|\nabla \phi(x)|^2{\rm d}x=1\right\}.
$$

Hereafter, we focus on the case $d=1$, that is,  ${\bf X}=(\Omega, X_t,\PP_x)$ is a $1$-dimensional 
Brownian motion. 
Let $\beta>-\frac{3}{2}$ with $\beta\ne1$. 
 Define 
\begin{align}
H_t^{\beta}:=\lim_{\varepsilon\to0}\int_0^t|X_s|^{\beta}(\text{\rm sgn}(X_s))\1_{\{|X_s|\geq\varepsilon\}}{\rm d}s.\label{eq:Cauchy*}
\end{align}
It is known that $H_t^{\beta}$ is well-defined and continuous. 
 $H_t^{\beta}$ is the so-called  
\emph{Hilbert transform of Brownian local time}. For $\beta=-1$, $H_t^{\beta}$ is called the \emph{Cauchy principal value of Brownian local time}. 
Set 
\begin{align*}
\ell_{\beta}(x):=
\left\{\begin{array}{ll} x\log|x|-x, &\text{ if }\quad \beta=-1 \\ \frac{1}{\beta+1}\int_0^x|y|^{\beta+1}{\rm d}y, &\text{ if }\quad \beta\in\,]\,-3/2,\,+\infty\,[\,\setminus\,\{\pm1\}\end{array}\right.,
\end{align*} 
$\rho_1:=\ell_{\beta}\phi_R$  and $\rho_2:=\ell_{\beta}(1-\phi_R)$
for sufficiently large $R>0$. Here $\phi_R$ is any smooth function such that 
$\phi_R=1$ on $B_R(0)$ and $\phi_R=0$ on $\Rd\setminus B_{R+1}(0)$.
It is known that $H_t^{\beta}$ is the (local) zero energy part $N_t^{\ell_{\beta}}$ in Fukushima's decomposition 
for $\ell_{\beta}\in H^1(\R)_{\rm loc}$ (see \cite{Fuk:strictdecomp}). 
Note that $|\rho_1'(x)|^2{\rm d}x$ is of Kato class and 
$|\ell_{\beta}'(x)|^2{\rm d}x$ is of local Kato class. 
If $\beta\leq 0$, then we can confirm that $\rho_2''$ is bounded and 
\begin{align*}
\ell_{\beta}''(x)(1-\phi_R(x)){\rm d}x=
\left\{\begin{array}{ll} \frac{1}{x}(1-\phi_R(x)){\rm d}x, &\text{ if }\quad  \beta=-1 \\ 
|x|^{\beta}({\rm sgn}(x))(1-\phi_R(x)){\rm d}x, &\text{ if }\quad  \beta\in ]-3/2,\,-1\,[\,\cup\,]-1,\,0\,]
\end{array}\right.
\in S_{K}^1({\bf X}).
\end{align*} 
We then see that $H_t^{\beta}=N_t^{\ell_{\beta}}=N_t^{\rho_1}-A_t^{\mu}$ for 
$\mu({\rm d}x):=-\frac{\,1\,}{2}\rho_2''(x){\rm d}x\in S_{K}^1({\bf X})-S_{K}^1({\bf X})$. 
Applying Theorem~\ref{thm:shortstability} again,  we obtain that 
there exist $C,C_1,C_2,k>0$ with $C_1\leq C_2$ such that 
\begin{align}
\frac{C^{-1}e^{-kt}}{\sqrt{2\pi t}}e^{-\frac{|x-y|^2}{C_1t}}\leq p_t^{\beta}(x,y)\leq 
\frac{Ce^{kt}}{\sqrt{2\pi t}}e^{-\frac{|x-y|^2}{C_2t}},\quad x,y\in \R, \quad t>0,\label{eq:HKExKato*}
\end{align}
where $p_t^{\beta}(x,y)$ is the integral kernel 
corresponding to the Feynman-Kac semigroup 
$P_t^{\beta}\!f(x):=
\EE_x[e^{H_t^{\beta}}f(X_t)]$ provided $\beta\in ]-3/2,\,0\,]$. 

\begin{rem}
{\rm The condition $\alpha\leq0$ in \cite[Theorem~6.2]{DKK2010} is incorrect. It should be corrected to be $\alpha<0$ for $\ell_{\alpha}''(x)(1-\phi_R(x)){\rm d}x=
|x|^{\alpha}({\rm sgn}(x))(1-\phi_R(x)){\rm d}x\in S_{K_{\infty}^+}^1({\bf X})$. 
}
\end{rem}

\subsection{Symmetric diffusion processes}\label{ex:Diffusion}
Throughout Sections \ref{ex:Diffusion}--\ref{s:rdj}, let $(E,d)$ be a locally compact separable metric space and 
$\m$ be a positive Radon measure on $E$ with full support. 
We also assume that every ball in $(E,d)$ is relatively compact.
Set $V(x,r):=\m(B(x,r))$ for $B(x,r)=\{y \in E \mid d(x,y)<r\}$. 
Throughout Sections \ref{ex:Diffusion}--\ref{s:rdj}, we 
assume that the \emph{volume doubling condition} (\emph{VD}): there exists a constant $C_D >0$ such that 
$V(x,2r) \le C_DV(x,r)$ for all $x \in E$ and $r>0$, and the \emph{reverse volume doubling condition} (\emph{RVD}): there exist constants $
c_1
>0$ and 
$d_1 >0$ such that 
\begin{align}
\frac{V(x,r_2)}{V(x,r_1)} \ge 
c_1\left(\frac{r_2}{r_1}\right)^{
d_1}\quad \text{for all}~  x \in E,~ 0< r_1 \le r_2 
 \le 2\, \text{diam}(E).
\label{eq:RVD}
\end{align}
(cf. \cite{GrigoHuLau2014}). (\emph{VD}) implies (\emph{RVD}) if $E$ is 
 connected. See  \cite[Proposition 2.1 and  a paragraph before Remark 2.1]{YZ}.

In this section, we follow the terminology and conditions in \cite{GrigoHuLau2014}.
Let $\phi$ be a fixed continuous increasing bijection on $]\,0,\,+\infty\,[$ satisfying 
$L(\beta,C_{\phi}^{-1})$ and $U(\beta',C_{\phi})$ for some $C_{\phi}\geq1$ and $1<\beta\leq\beta'$.	
	That is, 
\begin{align}
C_{\phi}^{-1}\left(\frac{R}{r}\right)^{\beta}\leq\frac{\phi(R)}{\phi(r)}
\leq C_{\phi}\left(\frac{R}{r}\right)^{\beta'}\quad\text{ for all }~0<r<R. 
\label{eq:PsiCond}
\end{align}
We consider the following condition: there exists a large $L>0$ such that
\begin{align}
\text{\rm ess-}\sup_{r \ge L}\frac{r\phi'(r)}{\phi(r)} <\infty.
\label{condPhi}
\end{align}
It is easy to see that \eqref{eq:PsiCond} and \eqref{condPhi} 
are satisfied for $\phi(r)=r^{\beta}$ with some $\beta >1$. 
Let $(\E,\F)$ be a strongly local irreducible regular Dirichlet form on $L^2(E;\m)$ and ${\bf X}=(\Omega,X_t,\PP_x, x \in E_\partial)$ an $\m$-symmetric diffusion process associated to $(\E,\F)$. 
${\bf X}$ is said to be \emph{conservative}, denoted by (\emph{C}) in short,  
if $\PP_x(X_t\in E\text{ for all }t>0)=1$ $\m$-a.e.~$x\in E$. 
We set 
\begin{align}\label{e:Phi}
\Phi(s,t)=\Phi(\phi, s,t)
:=\sup_{r>0}\left\{\frac{s}{r}-\frac{t}{\phi(r)} \right\}.
\end{align} 
\begin{df}[Green kernel bounds]
{\rm We say that \emph{the condition $(G)_{\phi}$} holds if there exist constants $\kappa\in\,]\,0,\,1\,[$ and $C\in\,]\,0,\,+\infty\,[$ such that, for any 
open ball $B=B_R(x_0)$, the Green kernel $R^B(x,y)$, $x,y\in B$ exists, is jointly continuous off diagonal, and satisfies 
\begin{align*}
R^B(x_0,y)&\leq C\int_{\kappa d(x_0,y)}^R\frac{\phi(s){\rm d}s}{s V(x_0,s)} \quad \text{ for }0<d(x_0,y)<R,
\\ R^B(x_0,y)&\geq C^{-1}\int_{\kappa d(x_0,y)}^R\frac{\phi(s){\rm d}s}{s V(x_0,s)} \quad\text{ for }0<d(x_0,y)<\kappa R.
\end{align*} 
}
\end{df}

\begin{df}[{(\!\emph{UE})\boldmath${}_{\phi}$, (\!\emph{LE})\boldmath${}_{\phi}$}]\label{df:UE}
{\rm For $c\in\,]\,0,\,1\,]$ and $k\geq0$, we say that 
(\!\emph{UE})${}_{\phi}^{c,k}$ holds if the heat kernel  $p_t(x,y)$ of 
{\bf X} exists and satisfies the following upper estimate
\begin{align}
 p_t(x,y)&\leq \frac{Ce^{kt}}{V(x,\phi^{-1}(t))}
\exp{\left(-\frac{\,1\,}{2}\Phi\left(cd(x,y),t \right) \right)} \label{eq:UE} 
\end{align}
for all $t>0$ and $\m$-a.e.~$x,y\in E$, where $C>0$ is a constant  independent of $x,y,t$. 

For $c\in\,[\,1,\,+\infty\,[$ and $k\geq0$, we say that 
(\emph{LE})${}_{\phi}^{c,k}$ holds 
if the heat kernel  $p_t(x,y)$ of 
{\bf X} exists and satisfies the following lower estimate
\begin{align}
 p_t(x,y)\geq \frac{Ce^{-kt}}{V(x,\phi^{-1}(t))}
\exp{\left(-c\Phi\left(cd(x,y),t \right) \right)} 
\label{eq:LE}
\end{align}
for all $t>0$ and $\m$-a.e.~$x,y\in E$, where $C>0$ is a constant  independent of $x,y,t$.   

We say that 
(\emph{UE})${}_{\phi}^{k}$ (resp.~(\emph{LE})${}_{\phi}^{k}$) holds if it satisfies 
(\emph{UE})${}_{\phi}^{c,k}$ (resp.~(\emph{LE})${}_{\phi}^{c,k}$) for some 
$c\in\,]\,0,\,1\,]$ (resp.~$c\in\,[\,1,\,+\infty\,[$). 
In particular, we say that 
(\!\emph{UE})${}_{\phi}$ (resp.~(\emph{LE})${}_{\phi}$) holds if it satisfies (\!\emph{UE})${}_{\phi}^{k}$ (resp.~(\emph{LE})${}_{\phi}^{k}$) with $k=0$.  
}
\end{df}

\begin{df}[{(\emph{NLE})\boldmath${}_{\phi}$}]\label{df:NLE} 
{\rm We say that 
(\emph{NLE})${}_{\phi}^{k}$ holds if the heat kernel  $p_t(x,y)$ of 
{\bf X} exists and satisfies the following lower estimate 
\begin{align}
 p_t(x,y)\geq \frac{Ce^{-kt}}{V(x,\phi^{-1}(t))}\label{eq:NLE}
 \end{align}
for all $t>0$ and $\m$-a.e.~$x,y\in E$, such that 
\begin{align}
d(x,y)\leq\eps\phi^{-1}(t),\label{eq:NLEL}
\end{align}
where $C, \eps>0$ and $k\geq0$ are constants independent of $x,y,t$. In particular, we say that 
(\emph{NLE})${}_{\phi}$ holds if 
(\emph{NLE})${}_{\phi}^{k}$ holds with $k=0$. 
}
\end{df}

\begin{rem}\label{rem:HK}
{\rm 
\begin{enumerate}
\item\label{item:HK1}
Clearly, $\Phi(s,t)=t\Phi(s/t,1)$. 
If $\phi(r)=Cr^{\beta}$ with some $C>0$ and $\beta>1$, then $\Phi(s,1)=cs^{\beta/(\beta-1)}$. Consequently, under 
\eqref{eq:PsiCond}, we always have $\Phi(s,1)\geq cs^{\beta/(\beta-1)}$ for some $c>0$. In particular, $\Phi(s,t)\geq0$ under \eqref{eq:PsiCond}.  
\item\label{item:HK2} 
It is easy to see that under \eqref{eq:NLEL} the term $\Phi(cd(x,y),t)$ in \eqref{eq:UE} is bounded by a constant, so that the upper bound 
(\emph{UE})${}_{\phi}^c$ is consistent with (\emph{NLE})${}_{\phi}$. 
 \item\label{item:HK3}  
It is known (cf.~\cite{BGK}, \cite{Grigo2014.1}) that (\emph{UE})${}_{\phi}$$+$(\emph{NLE})${}_{\phi}$ implies that the heat kernel $p_t(x,y)$ admits a locally H\"older continuous in $x,y$ version, so that \eqref{eq:UE} and 
\eqref{eq:NLE} are a posteriori true for all $x,y\in E$ under (\emph{UE})${}_{\phi}$$+$(\emph{NLE})${}_{\phi}$. 
It is known (cf.~\cite[Remark~6.12]{GrigoTelks2012}) that (\emph{LE})${}_{\phi}$ implies (\emph{NLE})${}_{\phi}$. Hence \eqref{eq:UE} and 
\eqref{eq:LE} are a posteriori true for all $x,y\in E$ under
(\emph{UE})${}_{\phi}$$+$(\emph{LE})${}_{\phi}$. 

\item\label{item:HK4} 
By \cite[Theorems~1.2 and 1.3]{GrigoHuLau2014}, (\emph{UE})${}_{\phi}$$+$(\emph{NLE})${}_{\phi}$ (resp.~(\emph{UE})${}_{\phi}$$+$(\emph{C})) is stable in the following sense; if two strongly local regular Dirichlet forms $(\E^{(1)},\F^{(1)})$ on $L^2(E;\m^{(1)})$ and $(\E^{(2)},\F^{(2)})$ on $L^2(E;\m^{(2)})$  having common domain $\F^{(1)}=\F^{(2)}$ satisfy 
$C^{-1}\m^{(1)}\leq\m^{(2)}\leq C\m^{(1)}$, 
$C^{-1}\E^{(1)}\leq \E^{(2)}\leq C\E^{(1)}$ on $\F^{(1)}\times \F^{(1)}$, and  $(\E^{(1)},\F^{(1)})$ admits (\emph{UE})${}_{\phi}$$+$(\emph{NLE})${}_{\phi}$ (resp.~(\emph{UE})${}_{\phi}$$+$(\emph{C})) for some positive constants, 
then  $(\E^{(2)},\F^{(2)})$  also does for some constants. 
This shows that {\bf (A.1)} (resp.~{\bf (A.2)})
with $\mathbb{T}:=[\,0,\,+\infty\,[^{\,2}$ 
holds under 
 (\emph{UE})${}_{\phi}$$+$(\emph{C}) (resp.~(\emph{UE})${}_{\phi}$$+$(\emph{NLE})${}_{\phi}$) with  
\begin{align*}
\hspace{-1cm}\phi_1(t,x,y):=\1_{\{d(x,y)\leq\eps\phi^{-1}(t)\}}\frac{1}{V(x,\phi^{-1}(t))}\quad\text{and}\quad \phi_2(t,x,y):=\frac{1}{V(x,\phi^{-1}(t))}
\exp{\left(-\frac{\,1\,}{2}\Phi\left(d(x,y),t \right) \right)}.
\end{align*} 

\item\label{item:HK5} 
By \cite[Theorem~1.2]{GrigoHuLau2014}, (\emph{UE})${}_{\phi}$$+$(\emph{NLE})${}_{\phi}$ is equivalent to $(G)_{\phi}$. In this case, under the transience of {\bf X}, we know the existence of global 
Green kernel $R(x,y)$ and it satisfies that there exist $\kappa\in\,]\,0,\,1\,[$ and $C\in \,]\,0,\,+\infty\,[$ such that 
\begin{align*}
C^{-1}\int_{\kappa d(x,y)}^{\infty}\frac{\phi(s){\rm d}s}{sV(x,s)}\leq R(x,y)\leq C\int_{\kappa d(x,y)}^{\infty}\frac{\phi(s){\rm d}s}{sV(x,s)}\quad \text{ for }x,y\in E\text{ with }x\ne y.
\end{align*} 
\end{enumerate}
}
\end{rem}
To apply our main results in this framework, 
we need to assume {\bf (AC)} for {\bf X}. 
Let $P_t^A\!f(x)=\EE_x[e_A(t)f(X_t)]$ be the Feynman-Kac semigroup of ${\bf X}$ defined by $e_A(t)=\exp(N_t^u+A_t^{\mu_1}-A_t^{\mu_2})$.  
Then we have the following:

\begin{thm}\label{applthm:mainresult1}
Suppose that {\bf X} {{satisfies {\bf (AC)} and it}} is transient. 
Let $u\in  {\F}_{\loc}\cap QC(E_{\partial})$ be a bounded finely continuous (nearly) Borel function on $E$. 
 Assume ${\mu}_1\in S_{N\!K_1}^1({\bf X})$, $\mu_{\<u\>}\in S_{N\!K_{\infty}}^1({\bf X})$ and   
$\mu_2\in S_{D_0}^1({\bf X})$. Then we have the following: 
\begin{enumerate}

\item\label{item:applmainresult1*} 
Suppose that {\bf X} is conservative. If 
$\lambda^{\Q}(\overline{\mu}_1) >0$ and 
{\rm(}UE\,{\rm)}${}_{\phi}$ holds, then there exists an integral kernel $p_t^A(x,y)$ of the Feynman-Kac semigroup $(P^A_t)_{t>0}$ satisfying \eqref{eq:UE} with $k=0$ and 
some $c\in]\,0,\,1\,]$ 
 (by replacing $p_t(x,y)$ with $p_t^A(x,y)$). 

\item\label{item:applmainresult1} 
If 
$\lambda^{\Q}(\overline{\mu}_1) >0$ and  
{\rm(}UE\,{\rm)}${}_{\phi}$$+${\rm(}NLE\,{\rm)}${}_{\phi}$ 
{\rm(}resp.~{\rm(}UE\,{\rm)}${}_{\phi}$$+${\rm(}LE\,{\rm)}${}_{\phi}${\rm)} holds, then there exists an integral kernel $p_t^A(x,y)$ of the Feynman-Kac semigroup $(P^A_t)_{t>0}$ in the strict sense satisfying 
\eqref{eq:UE} and \eqref{eq:NLE}
{\rm(}resp.~\eqref{eq:UE} and \eqref{eq:LE}{\rm)} 
for all $x,y\in E$ with $k=0$ and some adequate constants $c>0$ 
 (by replacing $p_t(x,y)$ with $p_t^A(x,y)$).

\item\label{item:applmainresult2}  
Assume that \eqref{condPhi} holds.  Suppose that  
{\rm(}UE\,{\rm)}${}_{\phi}$$+${\rm(}NLE\,{\rm)}${}_{\phi}$ 
for all $x,y\in E$, or 
{\rm(}LE\,{\rm)}${}_{\phi}$ for all $x,y\in E$ holds.
If there exists an integral kernel $p_t^A(x,y)$ of the Feynman-Kac semigroup $(P^A_t)_{t>0}$ in the strict sense satisfying 
\eqref{eq:UE} 
for all $x,y\in E$ with $k=0$ and $c\in]\,0,\,1\,]$ 
 (by replacing $p_t(x,y)$ with $p_t^A(x,y)$),  then $\lambda^{\Q}(\overline{\mu}_1) >0$.
\end{enumerate}
\end{thm}
\begin{pf}
\eqref{item:applmainresult1*}:  By use of  \cite[Theorem 1.3]{GrigoHuLau2014}, (\emph{UE})${}_{\phi}$$+$(\emph{C}) 
is stable under the change of Dirichlet form satisfying \eqref{eq:comparisonmeasure} and \eqref{eq:comparisondiffusion} for some $C_E>0$. 
So, {\bf (A.1)}
with 
$\mathbb{T}=[\,0,\,+\infty\,[^{\,2}$ is satisfied for    
$\phi_2(t,x,y)=\frac{1}{V(x,\phi^{-1}(t))}\exp\left(-\frac{\,1\,}{2}\Phi(d(x,y),t )\right)$. The conclusion follows from 
Theorem \ref{thm:mainresult1}\eqref{item:mainresult1}. 

\eqref{item:applmainresult1}: 
It is known that the Poincar\'e inequality and the generalized capacity inequality  are stable with respect to the quasi-isometry of Dirichlet forms (cf. \cite{GrigoHuLau2014}). From this and \cite[Theorem 1.2]{GrigoHuLau2014}, we see that  (\emph{UE})${}_{\phi}$$+$(\emph{NLE})${}_{\phi}$ and (\emph{UE})${}_{\phi}$$+$(\emph{LE})${}_{\phi}$ of the heat kernel $p_t(x,y)$ in the strict sense are stable under the change of Dirichlet form satisfying \eqref{eq:comparisonmeasure} and \eqref{eq:comparisondiffusion} for some $C_E>0$. 
So, {\bf (A.2)}
with 
$\mathbb{T}=[\,0,\,+\infty\,[^{\,2}$  is satisfied for 
$\phi_2(t,x,y)=\frac{1
}{V(x,\phi^{-1}(t))}\exp\left(-\frac{\,1\,}{2}\Phi(
d(x,y),t )\right)$ 
and 
$\phi_1(t,x,y)=\frac{1
}{V(x,\phi^{-1}(t))}\1_{\{d(x,y)<\eps\phi^{-1}(t)\}}$ (resp.~$\phi_1(t,x,y)=\frac{1
}{V(x,\phi^{-1}(t))}\exp\left(-c_1
\Phi(
d(x,y),t)\right)$ for some $c_1\in\,[\,1,\,+\infty\,[$) 
under (\emph{UE})${}_{\phi}$$+$(\emph{NLE})${}_{\phi}$ (resp.~(\emph{UE})${}_{\phi}$$+$(\emph{LE})${}_{\phi}$). 
 The conclusion now follows from  
Theorem \ref{thm:mainresult1}\eqref{item:mainresult1*}. 

\eqref{item:applmainresult2}: 
Since $\overline{\mu}_1=\overline{\mu}_1^*=\mu_1+\frac{\,1\,}{2}\mu_{\<u\>}\in S_{N\!K_1}^1({\bf X})$, we have 
$e^{-2u}\overline{\mu}_1\in S_{N\!K_1}^1({\bf Y})\subset S_{E\!K}^1({\bf Y})$ by \cite[Corollaries~5.1(5) and 5.2(2)]{KK:AnalChara}. 
One can take small $t_0>0$ such that $\sup_{x\in E}\EE_x^{Y}[A_{t_0}^{\overline{\mu}_1}]<1$, whence 
$\sup_{x\in E}\EE_x^Y[e^{A_{t_0}^{\overline{\mu}_1}}]\leq \frac{1}{1-\sup_{x\in E}\EE_x^Y[A_{t_0}^{\overline{\mu}_1}]}<\infty$ by Khasiminskii Lemma. 
Then for any $f\in \mathscr{B}_b(E)$ and $x\in E$, we have 
\begin{align*}
\int_E\left(\int_0^{t_0}p_t^A(x,y){\rm d}t\right)|f(y)|\m({\rm d}y)&=\int_0^{t_0}P_t^A|f|(x){\rm d}t\\
&\leq e^{2\|u\|_{\infty}}\|f\|_{\infty}\int_0^{t_0}\sup_{z\in E}\EE_z^Y[e^{A_t^{\overline{\mu}_1}}]{\rm d}t\\
&\leq e^{2\|u\|_{\infty}}\|f\|_{\infty}t_0\sup_{z\in E}\EE_z^Y[e^{A_{t_0}^{\overline{\mu}_1}}]<\infty,
\end{align*}
which implies that for each $x\in E$, $\int_0^{t_0}p_t^A(x,y){\rm d}t<\infty$ for $\m$-a.e.~$y\in E$. Now suppose that $p_t^A(x,y)$ satisfies (\emph{UE})${}_{\phi}$.
 for all $x,y\in E$. 

Assume first that \eqref{condPhi} 
and (\emph{UE})${}_{\phi}$
$+$(\emph{NLE})${}_{\phi}$ hold.  
By Remark~\ref{rem:HK}\eqref{item:HK5} with the 
transience of {\bf X}, for each $x\in E$, 
there exists a large $L>0$ such that $\int_{L}^{\infty}\frac{\phi(s)}{sV(x,s)}{\rm d}s<\infty$. By \eqref{condPhi}, 
\begin{align*}
\int_{L}^{\infty}\frac{\phi'(s)}{V(x,s)}{\rm d}s<\infty,
\end{align*} 
hence 
$\int_{\phi(L)}^{\infty}\frac{{\rm d}t}{V(x,\phi^{-1}(t))}<\infty$. 
Thus we have 
\begin{align*}
\int_{t_0}^{\infty}\frac{{\rm d}t}{V(x,\phi^{-1}(t))}<\infty,
\end{align*}
because $\phi$ is strictly increasing and $r\mapsto V(x,r)$ is lower semi continuous. 
Then we can conclude $\int_{t_0}^{\infty}p_t^A(x,y){\rm d}t<\infty$ for all $x,y\in E$, 
since $p_t^A(x,y)$ satisfies (\emph{UE})${}_{\phi}$ for all $x,y\in E$.  
Thus we obtain that for each $x\in E$, $R^A(x,y)<\infty$ for $\m$-a.e.~$y\in E$. 
Applying Theorem~\ref{thm:analgauge1}, 
we have $\lambda^{\Q}(\overline{\mu}_1)>0$. 

Next suppose that 
(\emph{LE})${}_{\phi}$ holds for all $x,y\in E$. 
 We recall from the proof of Theorem~\ref{thm:mainresult1}(3) that ${\sf d}=\{(x,y)\in E\times E\mid R(x,y)=\infty\}$ 
and $E^x=\{y\in E\mid (x,y)\notin{\sf d}\}$. Since for all $(x,y)\notin {\sf d}$
\begin{align*}
C^{-1}e^{-c\Phi(cd(x,y),t_0)}\int_{t_0}^{\infty}\frac{{\rm d}t}{V(x,\phi^{-1}(t))} &\leq 
C^{-1}\int_{t_0}^{\infty}\frac{1}{V(x,\phi^{-1}(t))}e^{-c\Phi(cd(x,y),t)} {\rm d}t\\
&\leq \int_{t_0}^\infty p_t(x,y){\rm d}t \leq R(x,y)<\infty,
\end{align*}
we can conclude $\int_{t_0}^{\infty}p_t^A(x,y){\rm d}t<\infty$ holds for all $(x,y)\notin{\sf d}$, 
because $p_t^A(x,y)$ satisfies (\!\emph{UE})${}_{\phi}$ for all $x,y\in E$.  
Hence we obtain that for each $x\in E$, $R^A(x,y)<\infty$ for 
$\m$-a.e.~$y\in E^x$. 
Since $\m(E\setminus E^x)=0$, we have $R^A(x,y)<\infty$ $\m$-a.e.~$y\in E$. 
Applying Theorem~\ref{thm:analgauge1}, 
we obtain $\lambda^{\Q}(\overline{\mu}_1)>0$.  
\end{pf}
\begin{rem}\label{rem:new}
{\rm
Note that under {\rm(}UE\,{\rm)}${}_{\phi}$ and $d_1>\beta$, ${\bf X}$ is transient. Indeed, by \eqref{eq:RVD} and \eqref{eq:PsiCond}, for $t\geq1$, 
\begin{align*}
\frac{V(x,\phi^{-1}(1))}{V(x,\phi^{-1}(t))}\leq c_1\left(\frac{\phi^{-1}(1)}{\phi^{-1}(t)}\right)^{d_1}\leq \frac{c_1}{t^{d_1/\beta}}
\end{align*}
implies $R(x,y):=\int_0^{\infty}p_t(x,y){\rm d}t<\infty$. 
}
\end{rem}
\begin{rem}\label{rem:Devyver}
{\rm Suppose that $\phi(r)=r^2$ and the intrinsic metric $\rho(x,y)=\sup\{f(x)-f(y)\mid f\in \F_{\loc}\cap C(E),~ \mu_{\<f\>}\leq\m\}$ is a non-degenerate complete metric compatible with the given topology. Then the (global) volume doubling condition (\emph{VD}) yields the conservativeness of {\bf X} by 
Sturm~\cite[Theorem~4]{St:DirI}. In this case, 
 (\emph{UE})${}_{\phi,\loc}$$+$(\emph{LE})${}_{\phi,\loc}$, a 
 localized version of  
 (\emph{UE})${}_{\phi}$$+$(\emph{LE})${}_{\phi}$, 
  is equivalent to (\emph{PI})${}_{\phi,\loc}$, a localized version of (weak) Poincar\'e inequality (\emph{PI})${}_{\phi}$, which is also equivalent to  
(\emph{PHI})${}_{\phi,\loc}$, a localized versions of parabolic Harnack inequality (\emph{PHI})${}_{\phi}$. 
Here (\emph{UE})${}_{\phi,\loc}$ (resp.~(\emph{LE})${}_{\phi,\loc}$) means that for any relatively compact open set $G$, \eqref{eq:UE} 
(resp.~\eqref{eq:LE}) with $k=0$ and some $c\in]\,0,\,1\,]$ (resp.~$c\in]\,1,\,+\infty\,[$) holds for all $t>0$ and for $\m$-a.e.~$x,y\in G$. 
It is easy to see that (\emph{PI})${}_{\phi,\loc}$$+$(\emph{VD}) (consequently (\emph{UE})${}_{\phi,\loc}$$+$(\emph{LE})${}_{\phi,\loc}$) is stable 
under the change of Dirichlet form satisfying \eqref{eq:comparisonmeasure} and \eqref{eq:comparisondiffusion} for some $C_E>0$ in view of the remark after  
\eqref{eq:comparisonjumpdensity}.  
Moreover, (\emph{PHI})${}_{\phi,\loc}$ 
yields the existence of 
local H\"older continuous heat kernel $p_t(x,y)$. So (\emph{UE})${}_{\phi}$ for all $x,y\in E$ is a posteori  true under  (\emph{UE})${}_{\phi}$$+$(\emph{LE})${}_{\phi,\loc}$. Then 
we can strengthen the assertion of  
Theorem~\ref{applthm:mainresult1}\eqref{item:applmainresult1*} 
so that \eqref{eq:UE} holds for all $x,y\in E$ with $k=0$ by replacing 
$p_t(x,y)$ with $p_t^A(x,y)$ provided $\lambda^{\Q}(\overline{\mu}_1) >0$ and (\emph{UE})${}_{\phi}$$+$(\emph{LE})${}_{\phi,\loc}$ holds. 
This statement extends Devyver~\cite[Theorem~4.1]{Devyver:HKRiesz}, because any complete smooth Riemannian manifold always satisfies  (\emph{UE})${}_{\phi,\loc}$$+$(\emph{LE})${}_{\phi,\loc}$ with $\phi(r)=r^2$.
}
\end{rem}

For $\alpha >0$, let ${\bf X}^{(\alpha)}$ be the $\alpha$-subprocess killed at rate $\alpha\m$. By virtue of Theorem \ref{applthm:mainresult1}, we obtain the next corollary without  assuming the transience of ${\bf X}$.    

\begin{cor}\label{applcor:mainresult1}
Suppose that {\bf X} satisfies {\bf (AC)}.
Let $u\in  {\F}_{\loc}\cap QC(E_{\partial})$ be a bounded finely continuous {\rm(}nearly{\rm)} Borel function on $E$. 
 Assume ${\mu}_1\in S_{N\!K_1}^1({\bf X}^{(\alpha)})$, $\mu_{\<u\>}\in S_{N\!K_{\infty}}^1({\bf X}^{(\alpha)})$ and $\mu_2\in S_{D}^1({\bf X})$. Then we have the following: 
\begin{enumerate}

\item\label{item:applcormainresult1*} 
Suppose that {\bf X} is conservative. 
If 
$\lambda^{\Q_{\alpha}}(\overline{\mu}_1) >0$ and  {\rm(}UE\,{\rm)}${}_{\phi}$ holds, then there exists an integral kernel $p_t^A(x,y)$ of the Feynman-Kac semigroup $(P^A_t)_{t>0}$ satisfying \eqref{eq:UE} 
with some constant $k:=k(\alpha)\geq0$ depending on $\alpha$ 
and some $c\in ]\,0,1\,]$ 
 (by replacing $p_t(x,y)$ with $p_t^A(x,y)$).

\item\label{item:applcormainresult1} 
If 
$\lambda^{\Q_{\alpha}}(\overline{\mu}_1) >0$ and {\rm(}UE\,{\rm)}${}_{\phi}$$+${\rm(}NLE\,{\rm)}${}_{\phi}$ 
(resp.~{\rm(}UE\,{\rm)}${}_{\phi}$$+${\rm(}LE\,{\rm)}${}_{\phi}$) holds, then there exists an integral kernel $p_t^A(x,y)$ of the Feynman-Kac semigroup $(P^A_t)_{t>0}$ in the strict sense satisfying 
\eqref{eq:UE} and \eqref{eq:NLE} 
{\rm(}resp.~\eqref{eq:UE} and \eqref{eq:LE}{\rm)}
for all $x,y\in E$ 
with some constant $k:=k(\alpha)\geq0$ depending on $\alpha$ 
and some adequate constants $c>0$ 
 (by replacing $p_t(x,y)$ with $p_t^A(x,y)$). 

\item\label{item:applcormainresult2}  
If there exists an integral kernel $p_t^A(x,y)$ of the Feynman-Kac semigroup $(P^A_t)_{t>0}$ in the strict sense satisfying 
\eqref{eq:UE} 
for all $x,y\in E$ with some $k\geq0$ and $c\in]\,0,\,1\,]$ 
 (by replacing $p_t(x,y)$ with $p_t^A(x,y)$), 
 then $\lambda^{\Q_{\alpha}}(\overline{\mu}_1) >0$ for $\alpha >k$.
\end{enumerate}
\end{cor}
\begin{pf}
The proof of \eqref{item:applcormainresult1} (resp.~\eqref{item:applcormainresult1*}) is similar to that of Corollary \ref{cor:mainresult}\eqref{item:cormainresult1*} 
(resp.~Corollary \ref{cor:mainresult}\eqref{item:cormainresult1})
together with Theorem~\ref{applthm:mainresult1}\eqref{item:applmainresult1} (resp.~Theorem~\ref{applthm:mainresult1}\eqref{item:applmainresult1*}). We shall prove \eqref{item:applcormainresult2}. By the same way of the proof of Theorem \ref{applthm:mainresult1}\eqref{item:applcormainresult2}, we can take a small $t_0 >0$ such that $\int_0^{t_0}e^{-\alpha t}p_t^A(x,y){\rm d}t < \infty$ for $\m$-a.e. $y \in E$ because $\overline{\mu}_1^* \in S_{E\!K}^1({\bf Y}^{(\alpha)})$. On the other hand, since $p_t^A(x,y)$ satisfies \eqref{eq:UE} 
for all $x,y\in E$ 
with some $k\geq0$ and $c\in]\,0,\,1\,]$ 
\begin{align*}
\int_{t_0}^{\infty}e^{-\alpha t}p_t^A(x,y){\rm d}t&\le \int_{t_0}^{\infty}\frac{e^{-(\alpha-k)t}}{V(x,\phi^{-1}(t))}\exp \left(-\frac{\,1\,}{2}\Phi(cd(x,y),t)\right)\!{\rm d}t \\
&\le \frac{1}{V(x,\phi^{-1}(t_0))}\int_{t_0}^{\infty}e^{-(\alpha-k)t}{\rm d}t < \infty.
\end{align*}
Thus we obtain that $R^A_{\alpha}(x,y) < \infty$ for $\m$-a.e. $y \in E$, equivalently $\lambda^{\Q_{\alpha}}(\overline{\mu}_1) >0$ by applying Theorem~\ref{thm:analgauge1} to the case $({\bf X}^{(\alpha)},\mu_{\<u\>},\mu_1,\mu_2)$.
\end{pf}

The following theorem can be proved based on Lemma~\ref{lem:extended} and Corollary~\ref{applcor:mainresult1}\eqref{item:applcormainresult1} in the same way of the proof of Theorem~\ref{thm:shortstability}. 
We omit its proof.

\begin{thm}\label{thm:ExtendedCor}
Suppose that {\bf X} satisfies {\bf (AC)}.
 Let $u\in  {\F}_{\loc}\cap QC(E_{\partial})$ be a bounded finely continuous (nearly) Borel function on $E$. 
 Assume ${\mu}_1\in S_{E\!K}^1({\bf X})$, $\mu_{\<u\>}\in S_{K}^1({\bf X})$ and $\mu_2\in S_{D}^1({\bf X})$. 
 \begin{enumerate}
 \item 
 If 
{\rm(}UE\,{\rm)}${}_{\phi}$$+${\rm(}NLE\,{\rm)}${}_{\phi}$ 
{\rm(}resp.~{\rm(}UE\,{\rm)}${}_{\phi}$$+${\rm(}LE\,{\rm)}${}_{\phi}${\rm)} holds,
  then  
there exists an integral kernel $p_t^A(x,y)$ of the Feynman-Kac semigroup $(P^A_t)_{t>0}$ in the strict sense satisfying 
\eqref{eq:UE} and \eqref{eq:NLE} 
{\rm(}resp.~\eqref{eq:UE} and \eqref{eq:LE}{\rm)}
for all $x,y\in E$ with 
some constant $k:=k(\alpha)\geq0$ depending on $\alpha$ and adequate constants $c>0$ 
 (by replacing $p_t(x,y)$ with $p_t^A(x,y)$). 
 
\item  If {\rm(}UE\,{\rm)}${}_{\phi}$$+${\rm(}C\,{\rm)} holds, then  
there exists an integral kernel $p_t^A(x,y)$ of the Feynman-Kac semigroup $(P^A_t)_{t>0}$ satisfying \eqref{eq:UE}   
with some constant $k:=k(\alpha)\geq0$ depending on $\alpha$ and 
an adequate constant $c>0$ 
 (by replacing $p_t(x,y)$ with $p_t^A(x,y)$). 
\end{enumerate}
\end{thm}

\subsection{Symmetric jump processes on metric measure spaces}\label{ex:jumpprocess}
Recall that  $(E,d)$ is  a locally compact separable metric space with  a positive Radon measure $\m$ on $E$ with full support and every ball in $(E,d)$ is relatively compact.
Recall $V(x,r)=\m(B(x,r))$.

We follow the terminology and conditions in \cite{BKKL21, CKW:Stability, CKW:H}.
We assume ${\rm diam}(E)=+\infty$ for simplicity. 
Since we have assumed (\emph{VD}) and (\emph{RVD}), we have that  there exist $c_1,c_2>0$, 
$d_2\geq d_1>0$ such that 
\begin{align}
c_1\left(\frac{R}{r}\right)^{d_1}\!\!\leq\frac{V(x,R)}{V(x,r)}\leq 
c_2\left(\frac{R}{r}\right)^{d_2}\quad\text{ for all }\quad 0<r<R<\infty,\quad x\in E.
\label{eq:Ahlfors}
\end{align}

Let $(\E,\F)$ be a regular Dirichlet form on $L^2(E;\m)$
with no killing part:
i.e., 
\begin{align*}
\E(f,g)=\E^{(c)}(f,g)+\int_{E\times E\setminus{\sf diag}}(f(x)-f(y))(g(x)-g(y))J({\rm d}x{\rm d}y), \quad \text{ for }f,g\in\F.
\end{align*}
where $ \E^{(c)}$ is the strongly local part of $ (\E, \F)$.
Let $\mu_{\<f\>}^{j}$ be the energy measure of $(\E,\F)$ defined by $\mu_{\<f\>}^{j}({\rm d}x)=2\int_E(f(x)-f(y))^2J({\rm d}x{\rm d}y)$ for 
$f\in\F$. 
Let {\bf X} be the associated $\m$-symmetric 
jump process. 
Assume further that the L\'evy system $(N,H)$ of {\bf X} satisfies 
$H_t=2t$, i.e., $J({\rm d}x{\rm d}y)=N(x,{\rm d}y)\m({\rm d}x)\, (=
\frac{\,1\,}{2}N(x,{\rm d}y)\mu_H({\rm d}x))$.

We will assume that  $\psi:\,[\,0,\,+\infty\,[\,\to\,[\,0,\,+\infty\,[$ is a 
	non-decreasing function satisfying 
	$L(\beta_1,C_L)$ and $U(\beta_2,C_U)$ 
for some $0<\beta_1 \le \beta_2$.	
	That is, 
	\begin{align}
C_L\left(\frac{R}{r}\right)^{\beta_1}\!\!\leq\frac{\psi(R)}{\psi(r)}\leq 
C_U\left(\frac{R}{r}\right)^{\beta_2}\quad\text{ for all }\quad 0<r<R<\infty.\label{eq:PhiAhlfors}
\end{align}
Note that, since $\psi$ satisfies $L(\beta_1,C_L)$,  we have $\lim_{t \to 0} \psi(t) = 0$.

Recall that $\Phi$ is defined in \eqref{e:Phi} through $\phi$.

\begin{df}[{\boldmath$\emph{HK}(\phi , \varphi, \psi)$, \boldmath$\emph{UHK}(\phi, \varphi, \psi)$, \boldmath$\emph{SHK}(\phi, \varphi, \psi)$}]\label{D:1.11}
{\rm 
Suppose  that  $\phi:\,]\,0,\,+\infty\,[\,\to\,[\,0,\,+\infty\,[$ is a non-decreasing function satisfying 
 $L(\alpha_1,c_L)$ and $U(\alpha_2, c_U)$ with some $0<\alpha_1 \le \alpha_2$ and $c_L, c_U>0$. 
Suppose  that  $\varphi:\,]\,0,\,+\infty\,[\,\to\,[\,0,\,+\infty\,[$ is a non-decreasing function satisfying that for some $c >0$,  
\begin{align}\label{comp1}
\varphi(r)\le c \psi(r), \quad \mbox{for all} \quad r>0.
\end{align}

\begin{description}
		\item{\rm (i)} We say that $\emph{HK}(\phi, \varphi, \psi)^*$ holds if there exists a heat kernel $p_t (x,y)$ of {\bf X} in the strict sense  
		which has the
		following estimates: there exist $k \ge 0$, $\eta, a_0>0$ and $c \ge1$ such that for all $t>0$ and $x,y\in E$,
		\begin{align}
		&\hspace{-0.3cm}c^{-1}e^{-kt}\left( \frac{1}{V(x, \varphi^{-1}(t))} {\bf 1}_{\{ d(x,y) \le \eta\varphi^{-1}(t)\}}+ \frac{t}{V(x, d(x,y)) \psi(d(x,y))}{\bf 1}_{\{ d(x,y) > \eta\varphi^{-1}(t)\}} \right) \notag\\
		&\le p_t(x,y) \label{eq:HKpurejump}\\
		&\le  c e^{kt}\left( \frac{1}{V(x, \varphi^{-1}(t))} \wedge 
	\Big(	\frac{t}{V(x, d(x,y))\psi(d(x,y))} + \frac{1}{V(x,\phi^{-1}(t))} e^{-a_0 \;\Phi(d(x,y), t)}\Big)
		\right).\notag
		\end{align}
	We say that \emph{HK}$(\phi, \varphi, \psi)$ holds if \emph{HK}$(\phi, \varphi, \psi)^*$ holds with $k=0$.

		\item{(ii)} 	
		We say that \emph{UHK}$(\phi,\varphi,\psi)^*$ holds if there exist $k\geq0$ and  a  
heat kernel $p_t(x,y)$ of {\bf X} such that for all $t>0$ and $\m\text{-a.e.~}x,y\in E$,
\begin{align}
\hspace{-1cm}p_t(x,y)\precsim_k \frac{1}{V(x, \varphi^{-1}(t))} \wedge 
	\left(	\frac{t}{V(x, d(x,y))\psi(d(x,y))} + \frac{1}{V(x,\phi^{-1}(t))} e^{-a_0 \;\Phi(d(x,y), t)}\right)
		.\label{eq:UHKpurejump}
\end{align} 
We say $\emph{UHK}(\phi, \varphi,\psi)$  holds if $\emph{UHK}(\phi, \varphi,\psi)^*$ holds  with $k=0$.

		\item{\rm (iii)} We say that $\emph{SHK}(\phi, \varphi, \psi)^*$ holds 
		if there exists a heat kernel $p_t (x,y)$ of {\bf X} in the strict sense such that 
		  for all $t>0$ and $x,y\in E$,
		\begin{equation*}
p(t,x,y) \asymp_{k}  \frac{1}{V(x, \varphi^{-1}(t))}\wedge 
	\left(	\frac{t}{V(x, d(x,y))\psi(d(x,y))} + \frac{1}{V(x,\phi^{-1}(t))} e^{-\;\Phi(d(x,y), t)}\right). 
		\end{equation*}		
			We say that $\emph{SHK}(\phi, \varphi, \psi)$ holds if $\emph{SHK}(\phi, \varphi, \psi)^*$ holds with $k=0$. 	
	\end{description}
	}
\end{df}

\begin{rem}[\emph{HK}(\boldmath$\phi$)]
{\rm 
When $\psi \simeq \phi\simeq \varphi$, we have
\begin{align*}
 &\frac{1}{V(x, \varphi^{-1}(t))}\wedge 
\left(	\frac{t}{V(x, d(x,y))\psi(d(x,y))} + \frac{1}{V(x,\phi^{-1}(t))} e^{-a_0 \;\Phi(d(x,y), t)}\right)
		\\
		&\asymp 
\frac{1}{V(x,\phi^{-1}(t))}\land \frac{t}{V(x,d(x,y))\phi(d(x,y))}
\end{align*} 
Thus, when $\psi \simeq \phi\simeq \varphi$, our notions \emph{HK}$(\phi, \varphi,  \psi)$ and \emph{UHK}$(\phi, \varphi, \psi)$ 
are the same as \emph{HK}$(\phi)$ and \emph{UHK}$(\phi)$  in \cite{CKW:Stability}.
}
\end{rem}
\begin{df}\label{d:chain}
{\rm	We say that a metric space $(M,d)$ satisfies the \textit{chain condition} $\Ch(A)$ if there exists a constant $A \ge 1$ such that, for any $n \in \N$ and $x,y \in M$, there is a sequence $\{z_k\}_{k=0}^n$ of points in $M$ such that $z_0 =x, z_n=y$ and
	$$ d(z_{k-1},z_k) \le A \frac{d(x,y)}{n} \quad \mbox{for all} \quad k=1,\cdots, n. $$}
\end{df}

\subsubsection{Symmetric pure jump processes}\label{ex:jumpprocessd}
{\it In this subsection, we assume  $\E^{(c)}\equiv 0$ so that $(\E,\F)$ is a regular Dirichlet form on $L^2(E;\m)$
 of pure jump type. }

It is established recently in 
\cite{CKW:Stability} and 
\cite[Theorem 2.14, Corollary 2.15, Theorem 2.17, Corollary 2.18] {BKKL21} that, if
either $\phi \simeq \psi$ or
the lower scaling index $\alpha_1$ for $\phi$ is strictly greater than 1, then 
 \emph{HK}($\phi, \phi, \psi)$ (resp.~\emph{UHK}($\phi, \phi, \psi)$$+$(\emph{C})) is  stable under the change of Dirichlet form satisfying \eqref{eq:comparisonmeasure} and \eqref{eq:comparisonjumpdensity} for some $C_E>0$. Hence, 
{\bf (A.2)} (resp.~{\bf (A.1)}) is satisfied under 
\emph{HK}($\phi, \phi, \psi)$ (resp.~\emph{UHK}($\phi, \phi, \psi)$$+$(\emph{C})) for 
$(\E,\F)$. 
Consequently, we can apply Theorem~\ref{thm:mainresult1}, Corollary~\ref{cor:mainresult} and 
Theorem~\ref{thm:shortstability}. 
We state the following only:

\begin{thm}\label{applthm;stableSJP}
Suppose that {\bf X} satisfies {\bf (AC)} and it is transient. We assume that either $\phi \simeq \psi$ or
the lower scaling index $\alpha_1$ for $\phi$ is strictly greater than 1.

Let $u\in  {\F}_{\loc}\cap QC(E_{\partial})$ be a bounded finely continuous (nearly) Borel function on $E$. 
 Assume $\mu_1+N(e^{F_1}-1)\mu_H\in S_{N\!K_1}^1({\bf X})$, $\mu_{\<u\>}\in S_{N\!K_{\infty}}^1({\bf X})$ and   
$\mu_2+N(F_2)\mu_H\in S_{D_0}^1({\bf X})$. Then we have the following: 
\begin{enumerate}
\item\label{item:applstablemainresult1SJP} 
Suppose UHK($\phi, \phi, \psi)$ and (C) hold. If  
$\lambda^{\Q}(\overline{\mu}_1) >0$, 
then there exists an integral kernel $p_t^A(x,y)$ of the Feynman-Kac semigroup $(P^A_t)_{t>0}$ 
satisfying \eqref{eq:UHKpurejump} 
with $k=0$ and $\varphi=\phi$
{(by replacing $p_t(x,y)$ with $p_t^A(x,y)$)}.

\item\label{item:applstablemainresult1*SJP}
Suppose HK($\phi, \phi, \psi)$ holds. If  
$\lambda^{\Q}(\overline{\mu}_1) >0$, 
then there exists an integral kernel $p_t^A(x,y)$ of the Feynman-Kac semigroup $(P^A_t)_{t>0}$ in the strict sense   
satisfying \eqref{eq:HKpurejump} 
with $k=0$ and $\varphi=\phi$ {(by replacing $p_t(x,y)$ with $p_t^A(x,y)$)}.   

\item\label{item:applstablemainresult2SJP}  
If there exists an integral kernel $p_t^A(x,y)$ of the Feynman-Kac semigroup $(P^A_t)_{t>0}$ in the strict sense  satisfying 
\eqref{eq:HKpurejump}  with $k=0$ and $\varphi=\phi$ {(by replacing $p_t(x,y)$ with $p_t^A(x,y)$)},  
then $\lambda^{\Q}(\overline{\mu}_1) >0$.
\end{enumerate}
\end{thm}
\begin{thm}\label{thm:shortstabilitySJP}
Suppose that {\bf X} satisfies {\bf (AC)}. 
Let $u\in  {\F}_{\loc}\cap QC(E_{\partial})$ be a bounded finely continuous (nearly) Borel function on $E$. 
Suppose that  
$\mu_1+N(e^{F_1}-1)\mu_H\in S_{E\!K}^1({\bf X})$, 
$\mu_{\<u\>}\in S_{K}^1({\bf X})$ and 
$\mu_2+N(F_2)\mu_H\in  S_{D}^1({\bf X})$ hold. 
If HK($\phi, \phi, \psi)$ holds, 
then there exists an integral kernel 
$p_t^A(x,y)$ of the Feynman-Kac semigroup $(P^A_t)_{t>0}$ 
in the strict sense  
satisfying 
\eqref{eq:HKpurejump} with $\varphi=\phi$
for some constant $k \ge 0$ depending on 
$\beta_1, \beta_2, C_L, C_U, \alpha_1,  \alpha_2, c_L, c_U$ {(by replacing $p_t(x,y)$ with $p_t^A(x,y)$)}.
\end{thm}
\begin{rem}
{\rm Theorem~\ref{applthm;stableSJP} extends \cite{Tak:ultra} and \cite{KimKuwae:stability}. 
Theorem~\ref{thm:shortstabilitySJP} also partially extends \cite{Wang}, which treats the non-symmetric $\alpha$-stable like process {\bf X} with non-symmetric bounded functions $F_1$ and $F_2$. 
}
\end{rem}

\begin{rem}
{\rm If our metric space $(E,d)$ also satisfy the \textit{chain condition} $\Ch(A)$, then by \cite[Theorem 2.17]{BKKL21}
\emph{HK}($\phi, \phi, \psi)$ is equivalent to \emph{SHK}($\phi, \phi, \psi)$. Thus, under the \textit{chain condition} $\Ch(A)$,  \emph{HK}($\phi, \phi, \psi)$ can be replaced by \emph{SHK}($\phi, \phi, \psi)$ in Theorems \ref{applthm;stableSJP} and \ref{thm:shortstabilitySJP} above.}
\end{rem}

\begin{df}[\boldmath$\bJ_{\psi}$]
{\rm We say that $\bJ_{\psi}$ holds if for $\m$-a.e.~$x\in E$, 
$N(x,{\rm d}y)\ll \m({\rm d}y)$ and for $\m$-a.e.~$x,y\in E$, the Radon-Nikodym derivative $J(x,y)$ satisfies 
\begin{align*}
\frac{C_1}{V(x,d(x,y))\psi(d(x,y))}\leq J(x,y)\leq \frac{C_2}{V(x,d(x,y))\psi(d(x,y))}\quad \text{ for } \m\text{-a.e.~}x,y\in E.
\end{align*}
}
\end{df}

We now assume that our metric measure space $(E,d)$ allows a conservative diffusion process $Z = (Z_t)_{t \ge 0 }$ which has the transition density $q(t,x,y)$  with respect to $m$ satisfying (\emph{UE})${}_{f}^{a_0,0}$ 	and (\emph{NLE})${}_{f}^{0}$ in Definitions \ref{df:UE} and \ref{df:NLE}, that is  there exist constants  $C\ge1$ and $a_0>0$ such that for all $t>0$ and $x,y\in E$,
	\begin{align}\label{hkediff}
	\frac{C^{-1}}{V(x,f^{-1}(t))} {\bf 1}_{\{f(d(x,y)) \le t\}}\le q(t,x,y)\le \frac{C}{V(x,f^{-1}(t))}\exp\big(-a_0 F(d(x,y),t) \big),
	\end{align}
	where the function $F$ is defined as	
	\begin{align*}
	F(r,t) : =F(f, r,t)= \sup_{s>0} \left[ \frac{r}{s} - \frac{t}{f(s)} \right].
	\end{align*}

Recall that 
$\psi:\,]\,0,\,+\infty\,[\,\to\,]\,0,\,+\infty\,[$ is a non-increasing function which 
satisfies $L(\beta_1,C_L)$, $U(\beta_2,C_U)$. 
\emph{For the rest of this subsection, 
we assume that $f$ is a strictly increasing function satisfying $L(\gamma_1,c_F^{-1})$ and $U(\gamma_2,c_F)$ with some $1<\gamma_1 \le \gamma_2$, and $\psi$ and $f$ satisfy  
\begin{align}\label{e:intcon}
\int_0^{1} \frac{{\d} f(s)}{\psi(s)} \,<\infty.
\end{align}
}
We now define $\phi$ as 
\begin{equation*}
\phi(r):=  \frac{f(r)}{\int_0^r \frac{1}{\psi(s)}{\d}f(s)}, \quad r>0. 
\end{equation*}
Then $\phi$ is strictly increasing function satisfying  \eqref{comp1}, 
 { $U(\alpha_2, c_U)$
and $L(\alpha_1, c_L)$ for some $\alpha_2\ge\alpha_1>0$ and $c_U, c_L>0$} (see 
\cite[Section 3]{BKKL21}).

Using \cite[Theorem 2.19(iii)]{BKKL21} and our Theorem~\ref{thm:mainresult1}, Corollary~\ref{cor:mainresult} and 
Theorem~\ref{thm:shortstability}, we have the following. 

\begin{thm}\label{applthm;stableSJP1}
Suppose that {\bf X} satisfies {\bf (AC)}. 
We assume that the lower scaling index $\alpha_1$ for $\phi$ is strictly greater than 1.
Suppose $(E,d)$ satisfies $\Ch(A)$ and allows a conservative diffusion process  whose transition density satisfies \eqref{hkediff} with a strictly increasing function $f$ satisfying  \eqref{e:intcon} and $L(\gamma_1,c_F^{-1})$ and $U(\gamma_2,c_F)$ with some $1<\gamma_1 \le \gamma_2$. 
Let $u\in  {\F}_{\loc}\cap QC(E_{\partial})$ be a bounded finely continuous (nearly) Borel function on $E$. 

\begin{enumerate}
\item Suppose that {\bf X} is transient. 
 Assume $\mu_1+N(e^{F_1}-1)\mu_H\in S_{N\!K_1}^1({\bf X})$, $\mu_{\<u\>}\in S_{N\!K_{\infty}}^1({\bf X})$ and   
$\mu_2+N(F_2)\mu_H\in S_{D_0}^1({\bf X})$. 
If 
$\lambda^{\Q}(\overline{\mu}_1) >0$ and $\bJ_{\psi}$ holds, 
then there exists an integral kernel $p_t^A(x,y)$ of the Feynman-Kac semigroup $(P^A_t)_{t>0}$ in the strict sense   
satisfying \emph{SHK}($\phi, \phi, \psi)$ (by replacing $p_t(x,y)$ with $p_t^A(x,y)$).   

\item
Suppose that  
$\mu_1+N(e^{F_1}-1)\mu_H\in S_{E\!K}^1({\bf X})$, 
$\mu_{\<u\>}\in S_{K}^1({\bf X})$ and 
$\mu_2+N(F_2)\mu_H\in  S_{D}^1({\bf X})$ hold. 
If  $\bJ_{\psi}$ holds,
then there exists an integral kernel 
$p_t^A(x,y)$ of the Feynman-Kac semigroup $(P^A_t)_{t>0}$ 
in the strict sense  
satisfying 
\emph{SHK}($\phi, \phi, \psi)^*$
for some constant $k \ge 0$
(by replacing $p_t(x,y)$ with $p_t^A(x,y)$).
 \end{enumerate}
\end{thm}

	If $\Phi$ satisfies $U(\alpha_2, c_U)$ with $\alpha_2<d_1$
	where $d_1$ is the constant in 	
then $\bJ_{\psi}$  implies that Green function $G(x,y)$ satisfies 
\begin{align}
\label{e:GE}
G(x,y)\asymp \frac{\Phi(d(x,y))}{V(x, d(x,y))}\asymp \frac{\Phi(d(x,y))}{V(y, d(x,y))}\quad\text{for any}\;\;x,y \in E.
\end{align}

\subsubsection{Symmetric diffusion with jump}\label{ex:diffusionjump}
{\it In this subsection, 
we assume that 
the lower scaling index $\alpha_1$ for $\phi$ is strictly greater than 1,
and
$$
\phi(r) \le \psi(r) \quad \text{for } r \in \,]\,0,\, 1\,] \quad \text{and} \quad 
\phi(r) \ge \psi(r) \quad \text{for } r \in \,[\,1, \,+\infty\,[. 
$$
}
Very recently it has been established  in \cite[Theorems 1.13 and 1.14] {CKW:H} that, under the above assumptions,  \emph{HK}($\phi, \phi \wedge \psi, \psi)$ (resp.~\emph{UHK}($\phi, \phi \wedge \psi, \psi)$$+$(\emph{C})) is  stable under the change of Dirichlet form satisfying \eqref{eq:comparisonmeasure} and \eqref{eq:comparisonjumpdensity} for some $C_E>0$. Hence, 
{\bf (A.2)} (resp.~{\bf (A.1)}) is satisfied under 
\emph{HK}($\phi, \phi \wedge \psi, \psi)$ (resp.~\emph{UHK}($\phi, \phi \wedge \psi, \psi)$$+$(\emph{C})) for 
$(\E,\F)$. 
Consequently, we can also apply Theorem~\ref{thm:mainresult1}, Corollary~\ref{cor:mainresult} and 
Theorem~\ref{thm:shortstability}. 

\begin{thm}\label{applthm;stableSJP2}
Suppose that {\bf X} satisfies {\bf (AC)} and it is transient. 
Let $u\in  {\F}_{\loc}\cap QC(E_{\partial})$ be a bounded finely continuous (nearly) Borel function on $E$. 
 Assume $\mu_1+N(e^{F_1}-1)\mu_H\in S_{N\!K_1}^1({\bf X})$, $\mu_{\<u\>}\in S_{N\!K_{\infty}}^1({\bf X})$ and   
$\mu_2+N(F_2)\mu_H\in S_{D_0}^1({\bf X})$. Then we have the following: 
\begin{enumerate}
\item\label{item:applstablemainresult1SJP2} 
Suppose UHK($\phi, \phi \wedge \psi, \psi)$ and (C) hold. If  
$\lambda^{\Q}(\overline{\mu}_1) >0$, 
then there exists an integral kernel $p_t^A(x,y)$ of the Feynman-Kac semigroup $(P^A_t)_{t>0}$ 
satisfying \eqref{eq:UHKpurejump} 
with $k=0$ 
and $\varphi=\phi \wedge \psi$
{(by replacing $p_t(x,y)$ with $p_t^A(x,y)$)}.

\item\label{item:applstablemainresult1*SJP2}
Suppose HK($\phi, \phi \wedge \psi, \psi)$ holds. If  
$\lambda^{\Q}(\overline{\mu}_1) >0$, 
then there exists an integral kernel $p_t^A(x,y)$ of the Feynman-Kac semigroup $(P^A_t)_{t>0}$ in the strict sense   
satisfying \eqref{eq:HKpurejump} 
with $k=0$ and $\varphi=\phi \wedge \psi$ {(by replacing $p_t(x,y)$ with $p_t^A(x,y)$)}.   

\item\label{item:applstablemainresult2SJP2}  
If there exists an integral kernel $p_t^A(x,y)$ of the Feynman-Kac semigroup $(P^A_t)_{t>0}$ in the strict sense  satisfying 
\eqref{eq:HKpurejump}  with $k=0$ and $\varphi=\phi \wedge \psi$ {(by replacing $p_t(x,y)$ with $p_t^A(x,y)$)},  
then $\lambda^{\Q}(\overline{\mu}_1) >0$.
\end{enumerate}
\end{thm}
\begin{thm}\label{thm:shortstabilitySJP2}
Suppose that {\bf X} satisfies {\bf (AC)}. 
Let $u\in  {\F}_{\loc}\cap QC(E_{\partial})$ be a bounded finely continuous (nearly) Borel function on $E$. 
Suppose that  
$\mu_1+N(e^{F_1}-1)\mu_H\in S_{E\!K}^1({\bf X})$, 
$\mu_{\<u\>}\in S_{K}^1({\bf X})$ and 
$\mu_2+N(F_2)\mu_H\in  S_{D}^1({\bf X})$ hold. 
If HK($\phi, \phi \wedge \psi, \psi)$ holds, 
then there exists an integral kernel 
$p_t^A(x,y)$ of the Feynman-Kac semigroup $(P^A_t)_{t>0}$ 
in the strict sense  
satisfying 
\eqref{eq:HKpurejump} with $\varphi=\phi \wedge \psi$
for some constant $k \ge 0$ depending on 
$\beta_1, \beta_2, C_L, C_U, \alpha_1,  \alpha_2, c_L, c_U$ {(by replacing $p_t(x,y)$ with $p_t^A(x,y)$)}.
\end{thm}

\begin{rem}
{\rm If our metric space $(E,d)$ is connected and also satisfy the \textit{chain condition} $\Ch(A)$, then by \cite[Theorem 1.13]{CKW:H}
\emph{HK}($\phi, \phi \wedge \psi, \psi)$ is equivalent to \emph{SHK}($\phi, \phi \wedge \psi, \psi)$. Thus, under the \textit{chain condition} $\Ch(A)$,  \emph{HK}($\phi, \phi \wedge \psi, \psi)$ can be replaced by \emph{SHK}($\phi, \phi \wedge \psi, \psi)$ in Theorems \ref{applthm;stableSJP2} and \ref{thm:shortstabilitySJP2} above.}
\end{rem}

\subsection{Reflected symmetric diffusions with jumps on inner uniform subdomains}\label{s:rdj}

Recall that $E$ is a locally compact separable metric space,
and $\m$ a $\sigma$-finite Radon measure with full support
on $E$. Let $(\E^0, \F^0)$ be a strongly local regular Dirichlet form
on $L^2(E; m)$, and  $\mu^0_{\<u\>}$ be the ($\E^0$-) energy measure
of $u\in \F^0$ so that
 $\E^0 (u,u)= \frac{\,1\,}{2} \mu^0_{\<u\>}(E)$.
 We assume that the intrinsic metric
$\rho$  of $(\E^0, \F^0)$
defined by
$ \rho (x, y)=\sup\big\{f(x)-f(y)\mid \ f\in \FF^0 \cap C_c(E)
\hbox{ with } \mu_{\<f\>}^0({\d} z) \leq \m({\d} z)\big\}
$
is finite
for any $x,y\in E$
and
 induces the original topology on $E$, and
  that
  $(E, \rho)$ is a complete metric space.
  
  We assume that the Dirichlet form $(\E^0, \F^0)$ is {\it Harnack-type Dirichlet space} in \cite[Chapter 2]{GS} so that Aronson-type heat kernel estimates for the diffusion
process $Z^0$ associated with $(\E^0, \F^0)$ hold, that is (\!\emph{UE})${}_{g}^{c,0}$ and (\emph{LE})${}_{g}^{c,0}$ in Definition
\ref{df:UE}
 hold with $g(r)=r^2$.

For a connected open subset $D$  of 
$(E, \rho)$, define for $x, y\in D$,
\begin{equation*}
\rho_D (x, y)=\inf\{\hbox{length}(\gamma)\;\mid\; \hbox{a continuous curve }
\gamma \hbox{ in } D \hbox{ with } \gamma (0)=x \hbox{ and }
\gamma (1)=y\}.
\end{equation*}
Denote by ${\bar  D}$ the completion of $D$ under the metric $\rho_D$.
   We extend the definition of $\m|_D$ to ${\bar D}$ by setting
$\m|_D ({\bar D}\setminus D)=0$.
For notational simplicity, we continue use $\m$ to denote this
measure $\m|_D$.

We assume that  $D$ is  an unbounded inner uniform subdomain of $E$. That is, 
there are constants
 $C_1, C_2 \in \,]\,0,\,+\infty\,[$ such
that, for any $x, y \in D$, there exists a continuous curve
$\gamma_{x, y}: [\,0,\, 1\,] \to D$ with
$\gamma_{x, y}(0) = x$, $\gamma _{x, y}(1)=y$ whose length is at most $ C_1 \rho_D (x, y)$, and 
for any $z \in \gamma_{x,y}([\,0,\, 1\,])$,
$
\rho (z, \partial D):=\inf_{w \in \partial D} \rho (z, w)   \geq C_2{\rho_D (z, x) \rho_D (z, y)}/{\rho_D(x, y)}.
$

Let $(\E^0, \F^0_D)$ be the part Dirichlet form of $(\E^0, \F^0)$ on $D$.
By \cite[Proposition 2.13]{GS}, we have  that for $x,y \in D$,
$$
 \rho_D (x, y)=\sup\left\{f(x)-f(y)\;\left|\; \ f\in \F^0_{D,{\rm loc}} \cap
 C_c(D)
\hbox{ with } \mu_{\<f\>}^0 ({\d} z) \leq \m({\d} z) \right\}\right..
$$

{Denote $B_{{\bar D}}(x, r):=\{y\in {\bar D}\mid \rho_D (x, y)<r\}$
and}  $\m(B_{{\bar D}}(x, r))$ by $V_D(x, r)$.
Define $\F_D^{0, \rf} :=\{f\in \F^0_{D,{\rm loc}}\mid \mu_{\<f\>}^0(D)<\infty\}$
and
\begin{equation*}
 \E^{0, \rf} (f, f):=\frac{\,1\,}{2} \mu_{\<f\>}^0(D) \quad \hbox{for } f\in \F_D^{0, \rf}.
\end{equation*}
In this section we consider the Markov process {\bf X} on $\bar D$ associated with the following type of non-local
symmetric Dirichlet forms $(\E, \F)$ on $L^2(D; \m)$:
 \begin{equation*}
 \F=\F_D^{0, \rf} \cap L^2(D; \m),
 \end{equation*}
 and, for $u\in \F$,
 \begin{equation*}
 \E (u, u)  = \E^{0,\rf}(u, u)+
 \frac{\,1\,}{2} \int_{D\times D}
 (u(x)-u(y))^2 J(x, y)\, \m\otimes\m({\d} x{\d} y),
 \end{equation*}
where  $J(x,y)$ is a non-negative
symmetric measurable function  on $D\times D \setminus {\sf diag}$ satisfying certain conditions to be specified below.

\begin{df}\label{D:1.3}
{\rm  Let  $\beta \in \,[\,0, \,+\infty\,]$ and $\psi$
be a
strictly increasing function  on
 $[\,0, \,+\infty\,[$  with $\psi (0)=0$ and $\psi (1)=1$
that satisfies $L(\beta_1,C_L)$ and $U(\beta_2,C_U)$ for some $0<\beta_1 \le \beta_2<2$.
For  a non-negative symmetric measurable function  $J(x, y)$ on $D\times D \setminus {\sf diag}$, we say
\begin{enumerate}

\item[(i)]   $(\bJ_{\psi, \beta, \leq})$   holds if there are positive constants $\kappa_1$ and $\kappa_2$ so that
 \begin{equation*}
J(x, y) \leq \frac{\kappa_1 }{V_D(x,  \rho_D(x, y))\psi  (\rho_D(x, y) )\exp (\kappa_2 \rho_D(x, y)^\beta)}
\quad \hbox{for   } (x, y) \in D\times D \setminus {\sf diag};
\end{equation*}

\item[(ii)]    $(\bJ_{\psi, \beta, \geq})$   holds if there are positive constants $\kappa_3$ and $\kappa_4$ so that
 \begin{equation*}
J(x, y) \geq  \frac{\kappa_3 }{V_D(x,  \rho_D(x, y))\psi  (\rho_D(x, y) )\exp (\kappa_4 \rho_D(x, y)^\beta)}
\quad \hbox{for   } (x, y) \in D\times D \setminus {\sf diag};
\end{equation*}

\item[(iii)]    $(\bJ_{\psi, \beta})$   holds if  both $(\bJ_{\psi, \beta, \leq})$ and    $(\bJ_{\psi, \beta, \geq})$  hold
 with  possibly  different constants $\kappa_i$ in the upper and lower bounds;

\item[(iv)]    $(\bJ_{\psi,  0_+, \leq})$   holds if
there are positive constants $\kappa_5$ and $\kappa_6$ so that
  \begin{equation*} 
  \left\{
\begin{split}
&J(x, y) \leq  \frac{\kappa_5 }{V_D(x,  \rho_D(x, y)) \psi_*  (\rho_D(x, y) ) }
 \quad \hbox{for   } (x, y) \in D\times D \setminus {\sf diag}, \\
&  \sup_{x\in D}\int_{\{y\in D\,\mid\, \rho_D(x,y)>1\}} \rho_D(x,y)^2 J(x, y)\,
\m({\d} y)
 \leq \kappa_6<\infty,
\end{split}
\right.
\end{equation*}
 where
 \begin{equation}
 \label{e:phi*}
 \psi_*(r):= \psi(r) \1_{\{r \leq 1\}} + r^2\1_{\{r>1\}} \quad \hbox{for } r\geq 0 ;
 \end{equation}
 \end{enumerate}
 }
\end{df}

 Define for $\beta \in \,[\,0, \,+\infty\,]$, $x\in \bar D$, $t>0$ and $r\geq 0$,
 \begin{equation*}
 p_{\psi, \beta } (t,x, r) :=
  \frac1{V_D(x,\sqrt{t})}\, \exp ( - r^2/t)+ \left(
  \frac1{V_D(x,\psi^{-1}(t))} \wedge
 \frac{t}{V(x, r) \psi (r)  \exp (r^\beta) } \right).
\end{equation*}
Define for $\beta \in ]\,0, \,1\,]$,
  $$
 H_{\psi,\beta} (t, x, r) :=
 \begin{cases} \displaystyle
 \frac{1}{V_D(x,\sqrt{t})} \wedge   p_{\psi, \beta } (t,x, r) 
 \quad & \hbox{if } t \in \,]\,0,\,1\,] ,  \smallskip \\
 \displaystyle   \frac{1}{V_D(x,  \sqrt{t} ) } \exp  \left( - \left( r^\beta \wedge  ({r^2}/{t})\right) \right)
    &\hbox{if } t \in\,]\,1,\,+\infty\,[;
    \end{cases}
 $$
 for $\beta\in \,]\,1,\,+\infty\,[$,
$$
 H_{\psi,\beta} (t,x,r) :=
 \begin{cases} \displaystyle
 \frac{1}{V_D(x,\sqrt{t})} \wedge p_{\psi, \beta } (t,x, r)  
  \quad & \hbox{if } t \in ]\,0,\,1\,] \hbox{ and }  r \leq 1,  \smallskip \\
   \displaystyle  \frac{t}
    {V_D(x, r )\psi(r)}
     \exp \left(  - \left(r(1+\log ^+ (r/t))^{(\beta-1)/\beta}\right)\wedge r^\beta \right)
  & \hbox{if } t \in ]\,0,\,1\,] \hbox{ and }    r> 1 , \smallskip \\
  \displaystyle   \frac{1}{V_D(x,\sqrt{t})}\exp \left( - \left(  r\, ( 1+\log^+ (r/t) )^{(\beta-1)/\beta} \right) \wedge (r^2/t)\right)
  &\hbox{if } t \in \,]\,1,\,+\infty\,[  ; \end{cases}
$$
where $\log^+(x):=\log (x\vee 1)$,
and $H_{\psi,\infty} (t,x,r) := \lim_{\beta \to \infty} H_{\psi,\beta} (t, x, r) $ for $\beta =\infty$, that is,
$$
 H_{\psi,\infty} (t,x,r) :=   \begin{cases}  \displaystyle
 \frac{1}{V_D(x,\sqrt{t})}  p_{\psi, \beta } (t,x, r) 
 \qquad & \hbox{if } t \in ]\,0,\,1\,] \hbox{ and }  r \le1,
  \smallskip \\
 \displaystyle    \frac{t}
    {V_D(x, r )\psi(r)}
     \exp \left(  - r(1+\log^+ (r/t) ) \right)
  & \hbox{if } t \in ]\,0,\,1\,]  \hbox{ and }
r  >1 ,  \smallskip \\
 \displaystyle    \frac{1}{V_D(x,\sqrt{t})}\exp \left( -  \left( r\, \big( 1+\log^+ (r/t)  \big)\right)  \wedge (r^2/t)\right)
  &\hbox{if } t \in \,]\,1,\,+\infty\,[  .
  \end{cases}
$$
We further define for $x\in \bar D$, $t>0$ and $r\geq 0$,
$$
 H_{\psi,0_+} (t,x,r) :=  \frac{1}{V_D(x,
\sqrt{t}
 ) } \wedge  p_{\psi_0, 0} (t,x, r),
$$
where $\psi_*$ is given by \eqref{e:phi*}.

Combining \cite[Theorems 1.5 and 1.6]{CKKW} and  our Theorem~\ref{thm:mainresult1}, Corollary~\ref{cor:mainresult} and 
Theorem~\ref{thm:shortstability}, we also obtain the stability result  on  estimates for fundamental solutions under Feynman-Kac perturbations for this process.

\begin{thm}\label{applthm;stableSJP3}
Suppose that {\bf X} satisfies {\bf (AC)} and it 
is transient. 
Let $u\in  {\F}_{\loc}\cap QC(D_{\partial})$ be a bounded finely continuous (nearly) Borel function on $D$. 
 Assume $\mu_1+N(e^{F_1}-1)\mu_H\in S_{N\!K_1}^1({\bf X})$, $\mu_{\<u\>}\in S_{N\!K_{\infty}}^1({\bf X})$ and   
$\mu_2+N(F_2)\mu_H\in S_{D_0}^1({\bf X})$. Then we have the following: 
\begin{enumerate}
\item\label{item:applstablemainresult1SJP3} 
Suppose $({\bf J}_{\psi, 0_+, \leq})$ holds. If  
$\lambda^{\Q}(\overline{\mu}_1) >0$, 
then there exists an integral kernel $p_t^A(x,y)$ of the Feynman-Kac semigroup $(P^A_t)_{t>0}$ 
satisfying 
$$
p_t^A(x,y) \le  c_1  H_{\psi,0_+} (t,x, c_2 \rho_D(x, y) )
   \quad \hbox{for all } x,y\in {\bar D} \hbox{ and } t>0.
$$

\item\label{item:applstablemainresult1*SJP3}
Suppose $J(x,y)$ satisfies  $({\bf J}_{\psi_1,  \beta_*, \leq })$ and $({\bf J}_{\psi_2, \beta^*, \geq })$
for some strictly increasing functions
 $\psi_1,\psi_2$
 satisfying $\psi_i(0)=0$, $\psi_i(1)=1$ and 
 $L(\beta^i_1,C^i_L)$ and $U(\beta^i_2,C^i_U)$ for some $0<\beta^i_1 \le \beta^i_2<2$ 
 for $i=1,2$,
 and for $\beta_*\leq \beta^*$ in $ \{0_+\}\, \cup\, ]\,0, \,+\infty\,] $   excluding   $\beta_*=\beta^*=0_+$. If  
$\lambda^{\Q}(\overline{\mu}_1) >0$, 
then there exists an integral kernel $p_t^A(x,y)$ of the Feynman-Kac semigroup $(P^A_t)_{t>0}$ in the strict sense   
satisfying
 \begin{equation}\label{e:1.23}
 c_1
 H_{\psi_2, \beta^*}   ( t, x, c_2 \rho_D(x, y)  )  \leq p_t^A(x,y) \leq c_3 H_{\psi_1, \beta_*} (t, x, c_4 \rho_D(x, y)),
  \end{equation}

\item\label{item:applstablemainresult2SJP3}  
If there exists an integral kernel $p_t^A(x,y)$ of the Feynman-Kac semigroup $(P^A_t)_{t>0}$ in the strict sense  satisfying 
\eqref{e:1.23}  then $\lambda^{\Q}(\overline{\mu}_1) >0$.
\end{enumerate}
\end{thm}
\begin{thm}\label{thm:shortstabilitySJP3}
Suppose that {\bf X} satisfies {\bf (AC)}. 
Let $u\in  {\F}_{\loc}\cap QC(D_{\partial})$ be a bounded finely continuous (nearly) Borel function on $D$. 
Suppose that  
$\mu_1+N(e^{F_1}-1)\mu_H\in S_{D\!K}^1({\bf X})$, 
$\mu_{\<u\>}\in S_{K}^1({\bf X})$ and 
$\mu_2+N(F_2)\mu_H\in  S_{D}^1({\bf X})$ hold. 

Suppose $J(x,y)$ satisfies  $({\bf J}_{\psi_1,  \beta_*, \leq })$ and $({\bf J}_{\psi_2, \beta^*, \geq })$
for some strictly increasing functions
 $\psi_1,\psi_2$
 satisfying $\psi_i(0)=0$, $\psi_i(1)=1$ and 
 $L(\beta^i_1,C^i_L)$ and $U(\beta^i_2,C^i_U)$ for some $0<\beta^i_1 \le \beta^i_2<2$ 
 for $i=1,2$,
 and for $\beta_*\leq \beta^*$ in $ \{0_+\}\, \cup\, ]\,0,\,+\infty\,] $   excluding   $\beta_*=\beta^*=0_+$. 
Then there exists an integral kernel 
$p_t^A(x,y)$ of the Feynman-Kac semigroup $(P^A_t)_{t>0}$ 
in the strict sense  
satisfying 
$$
 c_1e^{-kt}
 H_{\psi_2, \beta^*}   ( t, x, c_2 \rho_D(x, y)  )  \leq p_t^A(x,y) \leq c_3 e^{kt}H_{\psi_1, \beta_*} (t, x, c_4 \rho_D(x, y)),
$$
for some $k \ge 0$. 
\end{thm}

\subsection{Symmetric jump processes with tempered jumps at infinity in $\R^d$}\label{ex:relastable}

Let $(\E,\F)$ be a pure jump type regular Dirichlet form on $L^2(\R^d)$ defined by
\begin{align*}
\E(f,g)=\int_{\R^d\times \R^d}(f(x)-f(y))(g(x)-g(y))J(x,y){\rm d}x{\rm d}y, \quad f,g \in \F=\overline{C_0(\R^d)}^{\E_1^{1/2}},
\end{align*}
where $J(x,y)$ is a symmetric Borel function on $\R^d\times \R^d\setminus {\sf diag}$ satisfying {\bf UJS} condition, that is,
there is a constant $c>0$ such that
for a.e $x, y\in \R^d$,
$$
 J(x,y) \le \frac{c}{r^d}
\int_{B(x,r)} J(z,y)\, {\d} z \quad
\hbox { whenever }  r\le  |x-y|/2.
 $$
 Clearly, if $J(x,y)$  satisfies {\bf UJS} condition, so does $\wt{J}(x,y)$ satisfying \eqref{eq:comparisonjumpdensity}.

We also assume that the jump kernel $J$ has the following estimates;
\begin{align}
J(x,y) \asymp \frac{1}{|x-y|^d\phi(|x-y|)e^{|x-y|^\beta}},
\label{JestimatesAsymp}
\end{align}
where $\beta\in \,]\,0,\, +\infty\,]$ and $\phi:\,]\,0,\,+\infty\,[\,\to\,]\,0,\,+\infty\,[$ is a strictly increasing function satisfying $L(\alpha_1,c_L)$ and $U(\alpha_2, c_U)$ with some $0<\alpha_1 \le \alpha_2<2$.
When $\beta=\infty$, \eqref{JestimatesAsymp} is equivalent to 
\begin{align*}
J(x,y) \asymp \frac{1}{|x-y|^d\phi(|x-y|)} {\bf 1}_{\{|x-y| \le 1\}}.
\end{align*}
Let ${\bf X}=(X_t, \PP_x)$ be the symmetric Hunt process on $\R^d$ associated to $(\E,\F)$. 
 It is known in \cite{CKK:TranAMS2011, CKK:AMSCorregium} that ${\bf X}$ is a conservative Feller process and admits a locally jointly H\"older continuous transition density, a heat kernel $p_t(x,y)$ on $]\,0,\,+\infty\,[\times \R^d\times \R^d$, which has the following estimates: 
\begin{align}
p_t(x,y) \asymp_k q_{\beta} (t, |x-y|)\quad\text{ for }(t,x,y) \in \,
]\,0,\,+\infty\,[\times \R^d \times \R^d
\label{estimatestablelike2}
\end{align}
with $k=0$, where
 for $\beta\in \,]\,0,\,1\,]$,
$$
q_\beta (t, r):=
\begin{cases}
\frac{1}{\phi^{-1}(t)^d}\wedge  \frac{t}{r^d \phi (r) \Phi (r)}&\text{ if } t\le 1,\\
 t^{-d/2}\exp \left(-
r^{\beta} \wedge  ({r^2}/t) \right)&\text{ if  } t>1;
\end{cases}
$$
and for $\beta \in \,]\,1,\,+\infty\,]$,
$$
q_\beta (t, r):=
\begin{cases}
\frac{1}{\phi^{-1}(t)^d}\wedge  \frac{t}{r^d \phi (r) \Phi (r)}&\text{ if }   t\le 1, r<1,\\
t \exp\left(- ( r (1 + \log^+({r}/{t}))^{(\beta-1)/{\beta}})\wedge  r^\beta\right)  &\hbox{ if } t\le 1, r \ge 1,\\
t^{-d/2}\exp\left(- ( r (1 + \log^+({r}/{t}))^{(\beta-1)/{\beta}})\wedge   ({r^2}/t)\right)  &\hbox{ if } t> 1.
\end{cases}
$$
Note that {\bf X} is transient if and only if $d\geq3$. 
The rotationally symmetric relativistic $\alpha$-stable process ${\bf X}$ 
on $\R^d$  with mass $m>0$ (see \cite[Introduction]{CKK:TranAMS2011} for definition)
is a typical example of this subsection. 

By Theorem~\ref{thm:mainresult1}\eqref{item:mainresult1*} (see also the proof of \cite[Proposition~2.3]{Wada}), we can obtain the following: 

\begin{thm}\label{applthm;stable}
Suppose that {\bf X} satisfies {\bf (AC)}.
Suppose $d\geq3$. 
Let $u\in  {\F}_{\loc}\cap QC(E_{\partial})$ be a bounded finely continuous (nearly) Borel function on $E$. 
 Assume $\mu_1+N(e^{F_1}-1)\mu_H\in S_{N\!K_1}^1({\bf X})$, $\mu_{\<u\>}\in S_{N\!K_{\infty}}^1({\bf X})$ and   
$\mu_2+N(F_2)\mu_H\in S_{D_0}^1({\bf X})$. 
If $\lambda^{\Q}(\overline{\mu}_1) >0$, 
then there exists an integral kernel $p_t^A(x,y)$ of the Feynman-Kac semigroup $(P^A_t)_{t>0}$ in the strict sense  
satisfying 
\eqref{estimatestablelike2} with $k=0$.

\end{thm}
The conclusion of Theorem \ref{applthm;stable} extends 
\cite[Theorem~1.1]{Wada}, which treats the case 
that {\bf X} is the rotationally symmetric relativistic $\alpha$-stable process with 
$\mu_1 \in S_{C\!K_\infty}^1({\bf X}) \cap S_0^{(0)}({\bf X})$ and $u=\mu_2=F_1=F_2=0$. 

\bigskip

The next corollary follows from Corollary \ref{cor:mainresult} \eqref{item:cormainresult1*}.
\begin{cor}\label{applcorstable:mainresult1}
Suppose that {\bf X} satisfies {\bf (AC)}. 
Let $u\in  {\F}_{\loc}\cap QC(E_{\partial})$ be a bounded finely continuous (nearly) Borel function on $E$. 
 Assume $\mu_1+N(e^{F_1}-1)\mu_H\in S_{N\!K_1}^1({\bf X}^{(\alpha)})$, $\mu_{\<u\>}\in S_{N\!K_{\infty}}^1({\bf X}^{(\alpha)})$ and $\mu_2+N(F_2)\mu_H\in S_{D}^1({\bf X})$ for a fixed $\alpha>0$. 
If $\lambda^{\Q_{\alpha}}(\overline{\mu}_1) >0$, 
then there exists an integral kernel $p_t^A(x,y)$ of the Feynman-Kac semigroup $(P^A_t)_{t>0}$ in the strict sense satisfying 
\eqref{estimatestablelike2} 
with some constant $k\ge 0$ depending on
$\beta, \alpha_1,c_L, \alpha_2, c_U$.
\end{cor}
 
The following theorem can be proved based on Corollary~\ref{applcorstable:mainresult1} in the same way of the proof of Theorem~\ref{thm:shortstability}. We omit its proof.  
\begin{thm}\label{thm:shortstabilityRelStable}
Suppose that {\bf X} satisfies {\bf (AC)}. 
Let $u\in  {\F}_{\loc}\cap QC(E_{\partial})$ be a bounded finely continuous (nearly) Borel function on $E$. 
Suppose  that  
$\mu_1+N(e^{F_1}-1)\mu_H\in S_{E\!K}^1({\bf X})$, 
$\mu_{\<u\>}\in S_{K}^1({\bf X})$ and 
$\mu_2+N(F_2)\mu_H\in  S_{D}^1({\bf X})$ hold. Then there exists an integral kernel 
$p_t^A(x,y)$ of the Feynman-Kac semigroup $(P^A_t)_{t>0}$ 
in the strict sense  
satisfying 
\eqref{estimatestablelike2} 
for some constant $k:=k(\alpha) \ge 0$ depending on $\beta, \alpha_1,c_L, \alpha_2, c_U$.
\end{thm}

\emph{\bf Acknowledgment.}
{\rm The authors would like to thank Professors Takashi Kumagai and 
Masayoshi Takeda for their valuable comments to the first draft of this paper. } 

\providecommand{\bysame}{\leavevmode\hbox to3em{\hrulefill}\thinspace}
\providecommand{\MR}{\relax\ifhmode\unskip\space\fi MR }
\providecommand{\MRhref}[2]{%
  \href{http://www.ams.org/mathscinet-getitem?mr=#1}{#2}
}
\providecommand{\href}[2]{#2}

\end{document}